\documentclass[11pt, a4paper]{article}
\usepackage[cp1251]{inputenc}
\usepackage{amsmath} \usepackage{euscript}
\usepackage{graphicx}
\usepackage{longtable}
\usepackage{mathrsfs}
\usepackage{amssymb, amscd}
\begin{document}

\newcounter{bnomer} \newcounter{snomer}
\newcounter{bsnomer}
\setcounter{bnomer}{0}
\renewcommand{\thesnomer}{\thebnomer.\arabic{snomer}}
\renewcommand{\thebsnomer}{\thebnomer.\arabic{bsnomer}}
\renewcommand{\refname}{\begin{center}\large{\textbf{References}}\end{center}}

\setcounter{MaxMatrixCols}{14}

\newcommand\restr[2]{{
  \left.\kern-\nulldelimiterspace 
  #1 
  \right|_{#2} 
}}

\newcommand{\sect}[1]{%
\setcounter{snomer}{0}\setcounter{bsnomer}{0}
\refstepcounter{bnomer}
\par\bigskip\begin{center}\large{\textbf{\arabic{bnomer}. {#1}}}\end{center}}
\newcommand{\sst}[1]{%
\refstepcounter{bsnomer}
\par\bigskip\textbf{\arabic{bnomer}.\arabic{bsnomer}. {#1}}\par}
\newcommand{\defi}[1]{%
\refstepcounter{snomer}
\par\medskip\textbf{Definition \arabic{bnomer}.\arabic{snomer}. }{#1}\par\medskip}
\newcommand{\theo}[2]{%
\refstepcounter{snomer}
\par\textbf{Theorem \arabic{bnomer}.\arabic{snomer}. }{#2} {\emph{#1}}\hspace{\fill}$\square$\par}
\newcommand{\mtheop}[2]{%
\refstepcounter{snomer}
\par\textbf{Theorem \arabic{bnomer}.\arabic{snomer}. }{\emph{#1}}
\par\textsc{Proof}. {#2}\hspace{\fill}$\square$\par}
\newcommand{\mcorop}[2]{%
\refstepcounter{snomer}
\par\textbf{Corollary \arabic{bnomer}.\arabic{snomer}. }{\emph{#1}}
\par\textsc{Proof}. {#2}\hspace{\fill}$\square$\par}
\newcommand{\mtheo}[1]{%
\refstepcounter{snomer}
\par\medskip\textbf{Theorem \arabic{bnomer}.\arabic{snomer}. }{\emph{#1}}\par\medskip}
\newcommand{\theobn}[1]{%
\par\medskip\textbf{Theorem. }{\emph{#1}}\par\medskip}
\newcommand{\theoc}[2]{%
\refstepcounter{snomer}
\par\medskip\textbf{Theorem \arabic{bnomer}.\arabic{snomer}. }{#1} {\emph{#2}}\par\medskip}
\newcommand{\mlemm}[1]{%
\refstepcounter{snomer}
\par\medskip\textbf{Lemma \arabic{bnomer}.\arabic{snomer}. }{\emph{#1}}\par\medskip}
\newcommand{\mprop}[1]{%
\refstepcounter{snomer}
\par\medskip\textbf{Proposition \arabic{bnomer}.\arabic{snomer}. }{\emph{#1}}\par\medskip}
\newcommand{\theobp}[2]{%
\refstepcounter{snomer}
\par\textbf{Theorem \arabic{bnomer}.\arabic{snomer}. }{#2} {\emph{#1}}\par}
\newcommand{\theop}[2]{%
\refstepcounter{snomer}
\par\textbf{Theorem \arabic{bnomer}.\arabic{snomer}. }{\emph{#1}}
\par\textsc{Proof}. {#2}\hspace{\fill}$\square$\par}
\newcommand{\theosp}[2]{%
\refstepcounter{snomer}
\par\textbf{Theorem \arabic{bnomer}.\arabic{snomer}. }{\emph{#1}}
\par\textbf{Sketch of the proof}. {#2}\hspace{\fill}$\square$\par}
\newcommand{\exam}[1]{%
\refstepcounter{snomer}
\par\medskip\textbf{Example \arabic{bnomer}.\arabic{snomer}. }{#1}\par\medskip}
\newcommand{\deno}[1]{%
\refstepcounter{snomer}
\par\textbf{Notation \arabic{bnomer}.\arabic{snomer}. }{#1}\par}
\newcommand{\lemm}[1]{%
\refstepcounter{snomer}
\par\textbf{Lemma \arabic{bnomer}.\arabic{snomer}. }{\emph{#1}}\hspace{\fill}$\square$\par}
\newcommand{\lemmp}[2]{%
\refstepcounter{snomer}
\par\medskip\textbf{Lemma \arabic{bnomer}.\arabic{snomer}. }{\emph{#1}}
\par\textsc{Proof}. {#2}\hspace{\fill}$\square$\par\medskip}
\newcommand{\coro}[1]{%
\refstepcounter{snomer}
\par\textbf{Corollary \arabic{bnomer}.\arabic{snomer}. }{\emph{#1}}\hspace{\fill}$\square$\par}
\newcommand{\mcoro}[1]{%
\refstepcounter{snomer}
\par\textbf{Corollary \arabic{bnomer}.\arabic{snomer}. }{\emph{#1}}\par\medskip}
\newcommand{\corop}[2]{%
\refstepcounter{snomer}
\par\textbf{Corollary \arabic{bnomer}.\arabic{snomer}. }{\emph{#1}}
\par\textsc{Proof}. {#2}\hspace{\fill}$\square$\par}
\newcommand{\nota}[1]{%
\refstepcounter{snomer}
\par\medskip\textbf{Remark \arabic{bnomer}.\arabic{snomer}. }{#1}\par\medskip}
\newcommand{\propp}[2]{%
\refstepcounter{snomer}
\par\medskip\textbf{Proposition \arabic{bnomer}.\arabic{snomer}. }{\emph{#1}}
\par\textsc{Proof}. {#2}\hspace{\fill}$\square$\par\medskip}
\newcommand{\hypo}[1]{%
\refstepcounter{snomer}
\par\medskip\textbf{Conjecture \arabic{bnomer}.\arabic{snomer}. }{\emph{#1}}\par\medskip}
\newcommand{\prop}[1]{%
\refstepcounter{snomer}
\par\textbf{Proposition \arabic{bnomer}.\arabic{snomer}. }{\emph{#1}}\hspace{\fill}$\square$\par}

\newcommand{\proof}[2]{%
\par\medskip\textsc{Proof{#1}}. \hspace{-0.2cm}{#2}\hspace{\fill}$\square$\par\medskip}

\makeatletter
\def\iddots{\mathinner{\mkern1mu\raise\p@
\vbox{\kern7\p@\hbox{.}}\mkern2mu
\raise4\p@\hbox{.}\mkern2mu\raise7\p@\hbox{.}\mkern1mu}}
\makeatother

\newcommand{\okr}[2]{%
\refstepcounter{snomer}
\par\medskip\textbf{{#1} \arabic{bnomer}.\arabic{snomer}. }{\emph{#2}}\par\medskip}

\newcommand{\Ind}[3]{%
\mathrm{Ind}_{#1}^{#2}{#3}}
\newcommand{\Res}[3]{%
\mathrm{Res}_{#1}^{#2}{#3}}
\newcommand{\epsi}{\epsilon}
\newcommand{\tri}{\triangleleft}
\newcommand{\Supp}[1]{%
\mathrm{Supp}(#1)}

\newcommand{\lee}{\leqslant}
\newcommand{\gee}{\geqslant}
\newcommand{\reg}{\mathrm{reg}}
\newcommand{\Ann}{\mathrm{Ann}\,}
\newcommand{\Cent}[1]{\mathbin\mathrm{Cent}({#1})}
\newcommand{\PCent}[1]{\mathbin\mathrm{PCent}({#1})}
\newcommand{\Exp}[1]{\mathbin\mathrm{Exp}({#1})}
\newcommand{\empr}[2]{[-{#1},{#1}]\times[-{#2},{#2}]}
\newcommand{\sreg}{\mathrm{sreg}}
\newcommand{\ilm}{\varinjlim}
\newcommand{\plm}{\varprojlim}
\newcommand{\codim}{\mathrm{codim}\,}
\newcommand{\chara}{\mathrm{char}\,}
\newcommand{\rk}{\mathrm{rk}\,}
\newcommand{\chr}{\mathrm{ch}\,}
\newcommand{\Ker}{\mathrm{Ker}\,}
\newcommand{\id}{\mathrm{id}}
\newcommand{\Ad}{\mathrm{Ad}}
\newcommand{\col}{\mathrm{col}}
\newcommand{\row}{\mathrm{row}}
\newcommand{\low}{\mathrm{low}}
\newcommand{\pho}{\hphantom{\quad}\vphantom{\mid}}
\newcommand{\fho}[1]{\vphantom{\mid}\setbox0\hbox{00}\hbox to \wd0{\hss\ensuremath{#1}\hss}}
\newcommand{\wt}{\widetilde}
\newcommand{\wh}{\widehat}
\newcommand{\ad}[1]{\mathrm{ad}_{#1}}
\newcommand{\tr}{\mathrm{tr}\,}
\newcommand{\GL}{\mathrm{GL}}
\newcommand{\SL}{\mathrm{SL}}
\newcommand{\SO}{\mathrm{SO}}
\newcommand{\Or}{\mathrm{O}}
\newcommand{\Sp}{\mathrm{Sp}}
\newcommand{\Sa}{\mathrm{S}}
\newcommand{\Ua}{\mathrm{U}}
\newcommand{\Mat}{\mathrm{Mat}}
\newcommand{\Stab}{\mathrm{Stab}}
\newcommand{\htt}{\mathfrak{h}}
\newcommand{\spt}{\mathfrak{sp}}
\newcommand{\slt}{\mathfrak{sl}}
\newcommand{\sot}{\mathfrak{so}}

\newcommand{\vfi}{\varphi}
\newcommand{\vpi}{\varpi}
\newcommand{\teta}{\vartheta}
\newcommand{\Bfi}{\Phi}
\newcommand{\Fp}{\mathbb{F}}
\newcommand{\Rp}{\mathbb{R}}
\newcommand{\Zp}{\mathbb{Z}}
\newcommand{\Cp}{\mathbb{C}}
\newcommand{\Ap}{\mathbb{A}}
\newcommand{\Pp}{\mathbb{P}}
\newcommand{\Np}{\mathbb{N}}
\newcommand{\ut}{\mathfrak{u}}
\newcommand{\at}{\mathfrak{a}}
\newcommand{\hei}{\mathfrak{hei}}
\newcommand{\nt}{\mathfrak{n}}
\newcommand{\mt}{\mathfrak{m}}
\newcommand{\rt}{\mathfrak{r}}
\newcommand{\rad}{\mathfrak{rad}}
\newcommand{\bt}{\mathfrak{b}}
\newcommand{\unt}{\underline{\mathfrak{n}}}
\newcommand{\gt}{\mathfrak{g}}
\newcommand{\vt}{\mathfrak{v}}
\newcommand{\pt}{\mathfrak{p}}
\newcommand{\Xt}{\mathfrak{X}}
\newcommand{\Po}{\mathcal{P}}
\newcommand{\PV}{\mathcal{PV}}
\newcommand{\Uo}{\EuScript{U}}
\newcommand{\Fo}{\EuScript{F}}
\newcommand{\Do}{\EuScript{D}}
\newcommand{\Eo}{\EuScript{E}}
\newcommand{\Jo}{\EuScript{J}}
\newcommand{\Iu}{\mathcal{I}}
\newcommand{\Mo}{\mathcal{M}}
\newcommand{\Nu}{\mathcal{N}}
\newcommand{\Ro}{\mathcal{R}}
\newcommand{\Co}{\mathcal{C}}
\newcommand{\Lo}{\mathcal{L}}
\newcommand{\Ou}{\mathcal{O}}
\newcommand{\Uu}{\mathcal{U}}
\newcommand{\Au}{\mathcal{A}}
\newcommand{\Vu}{\mathcal{V}}
\newcommand{\Du}{\mathcal{D}}
\newcommand{\Bu}{\mathcal{B}}
\newcommand{\Sy}{\mathcal{Z}}
\newcommand{\Sb}{\mathcal{F}}
\newcommand{\Gr}{\mathcal{G}}
\newcommand{\rtc}[1]{C_{#1}^{\mathrm{red}}}

\newcommand{\JSpec}[1]{\mathrm{JSpec}\,{#1}}
\newcommand{\MSpec}[1]{\mathrm{MSpec}\,{#1}}
\newcommand{\PSpec}[1]{\mathrm{PSpec}\,{#1}}
\newcommand{\APbr}[1]{\mathrm{span}\{#1\}}
\newcommand{\APbre}[1]{\langle #1\rangle}
\newcommand{\APro}[1]{\setcounter{AP}{#1}\Roman{AP}}\newcommand{\apro}[1]{{\rm\setcounter{AP}{#1}\roman{AP}}}
\newcommand{\ot}{\xleftarrow[]{}}
\newcounter{AP}


\author{Mikhail Ignatyev\and Alexey Petukhov}
\date{}
\title{The orbit method for locally nilpotent infinite-dimensional\\ Lie algebras}\maketitle
\begin{abstract} Let $\mathfrak{n}$ be a locally nilpotent infinite-dimensional Lie algebra over $\mathbb{C}$. Let $\mathrm{U}(\mathfrak{n})$ and $\mathrm{S}(\mathfrak{n})$ be its universal enveloping algebra and its symmetric algebra respectively. Consider the Jacobson topology on the primitive spectrum of $\mathrm{U}(\mathfrak{n})$ and the Poisson topology on the primitive Poisson spectrum of $\mathrm{S}(\mathfrak{n})$. 
We provide a homeomorphism between the corresponding topological spaces (on the level of points, it gives a bijection between the primitive ideals of $\mathrm{U}(\mathfrak{n})$ and $\mathrm{S}(\mathfrak{n})$). We also show that all primitive ideals of $\mathrm{S}(\mathfrak{n})$ from an open set in a properly chosen topology are  generated by their intersections with the Poisson center. Under the assumption that $\mathfrak{n}$ is a nil-Dynkin Lie algebra, we give two criteria for primitive ideals $I(\lambda)\subset\mathrm{S}(\mathfrak{n})$ and $J(\lambda)\subset\mathrm{U}(\mathfrak{n})$, $\lambda\in\mathfrak{n}^*$, to be nonzero. Most of these results generalize the 
known facts about primitive and Poisson spectrum for finite-dimensional nilpotent Lie algebras (but note that for a finite-dimensional nilpotent Lie algebra all primitive ideals $I(\lambda)$, $J(\lambda)$ are nonzero).

 \end{abstract}
\vspace{-0.4cm}\tableofcontents\newpage

\sect{Introduction}\label{sect:intro}\addcontentsline{toc}{section}{\ref{sect:intro}. Introduction}
We work over the field $\Cp$ of complex numbers. By definition, a locally nilpotent Lie algebra is a direct limit of nested finite-dimensional nilpotent Lie algebras. In this paper we discuss primitive ideals (respectively, primitive Poisson ideals) of the universal enveloping algebra $\Ua(\nt)$ (respectively, of the symmetric algebra $\Sa(\nt)$) of such a Lie algebra~$\nt$.


The corresponding theory for the finite-dimensional Lie algebras was developed in the mid-end of 20th century by the efforts (and insights) of many brilliant mathematicians. We would like to mention explicitly J. Dixmier, M. Duflo, A. Joseph, A. Kirillov, B. Kostant, O. Mathieu, M. Vergne among them. Of course, in the infinite-dimensional situation, we use the power, the intuition and the details of proofs of the finite-dimensional setting, but the final results are rather different. 

It is known that similar questions for the finite-dimensional (nilpotent) Lie algebras boil down to some questions for the corresponding coadjoint representations. This idea works almost in the same way for the infinite-dimensional setting. Hence, the new features of the infinite-dimensional setting are coming from the new features of the corresponding coadjoint representation.

The most straightforward difference is that the dual space of a countable-dimensional space is uncountable-dimensional. Next, it turns out that many infinite-dimensional locally nilpotent Lie algebras are centerless (and this is completely opposite to the finite-dimensional setting)~\cite{DimitrovPenkov1}. Moreover, in many cases even the center of the corresponding universal enveloping algebra consists of the elements of the ground field~\cite{IgnatyevPenkov1} (we enhance these statements for a class of nil-Dynkin infinite-dimensional locally nilpotent Lie algebras which will be described later).

The most important tool in the representation theory of finite-dimensional nilpotent Lie algebras is an algebraic version of the Kirillov's orbit method. 
For a given finite-dimensional nilpotent Lie algebra~$\nt$, this method establishes a homeomorphism between the space $\JSpec{\Ua(\nt)}$ of primitive ideals of $\Ua(\nt)$ (endowed with the Jacobson topology) and the space of coadjoint orbits on the dual space~$\nt^*$. By definition, the latter space is homeomorphic to the space $\PSpec{\Sa(\nt)}$ of primitive Poisson ideals of $\Sa(\nt)$ (endowed with the Poisson topology, see Subsection~\ref{sst:Poisson}). In more details, to each linear form $\lambda\in\nt^*$ one can attach the primitive ideal $J(\lambda)$ of $\Ua(\nt)$ and the primitive Poisson ideal $I(\lambda)$ of $\Sa(\nt)$, see Subsection~\ref{sst:orbit_method}. 
It turns out that each primitive ideal of $\Ua(\nt)$ (respectively, each primitive Poisson ideal of $\Sa(\nt)$) has the form $J(\lambda)$ (respectively, $I(\lambda)$). Furthermore, the map $I(\lambda)\mapsto J(\lambda)$ provides a homeomorphism between the spaces $\PSpec{\Sa(\nt)}$ and $\JSpec{\Ua(\nt)}$.

Our first main result claims that the orbit method still works in the infinite-dimensional situation (see Theorem~\ref{theo:orbit_method_ifd} for more details).
\theobn{Let $\nt$ be a countable-dimensional locally nilpotent complex Lie algebra. Then
\begin{equation*}
\begin{split}
&\text{\textup{i)} each primitive ideal of $\Ua(\nt)$ equals $J(\lambda)$ for a certain $\lambda\in\nt^*$};\\
&\text{\textup{ii)} each primitive Poisson ideal of $\Sa(\nt)$ equals $I(\lambda)$ for a certain $\lambda\in\nt^*$};\\
&\text{\textup{iii)} the map $I(\lambda)\mapsto J(\lambda)$ is a homeomorphism between $\PSpec(\Sa(\nt))$ and $\JSpec(\Ua(\nt))$}.\\
\end{split}
\end{equation*}}
In the finite-dimensional setting, given linear forms $\lambda,~\mu\in\nt^*$, the ideals $J(\lambda)$ and $J(\mu)$ coincide if and only if $\lambda$ and $\mu$ belong to the same coadjoint orbit. 
We can't provide an analogue of this result for a general locally nilpotent Lie algebra.
Nevertheless, we can prove a similar result for a properly chosen group together with its action on $\nt^*$ for a certain subclass of locally nilpotent Lie algebras, see Proposition~\ref{Porbmsoc}. This subclass consists of Lie algebras $\nt$ which can be exhausted by its finite-dimensional nilpotent ideals; such a Lie algebra is called {\it socle}. For example, the countable-dimensional Heisenberg algebra is socle. 
We also establish a version of Dixmier--Moeglin equivalence for socle Lie algebras.

Now, we turn to the results on nil-Dynkin algebras. It is well known that the Dynkin diagrams of types $A$, $B$, $C$, $D$, $E$, $F$, $G$ provide a description of simple finite-dimensional Lie algebras. As a byproduct of this procedure the same diagrams give a very detailed description of maximal nilpotent subalgebras of these simple Lie algebras. The Dynkin diagrams of types $A$, $B$, $C$, $D$ has countable analogues and every such an analogue defines the infinite-dimensional Lie algebra together with its maximal locally nilpotent subalgebra~\cite{DimitrovPenkov1}. This construction defines a wide variety of infinite-dimensional locally nilpotent Lie algebras, see Section~\ref{sect:nildynkin}. We refer to all Lie algebras defined by the above construction as nil-Dynkin algebras.

If $\nt$ is a finite-dimensional nilpotent Lie algebra then almost all primitive ideals of $\Ua(\nt)$ are centrally generated, i.e., are generated as ideals by their intersections with the center $Z(\nt)$ of $\Ua(\nt)$, see, e.g., Theorem~\ref{theo:almost_all_cent_gen_fin_dim}. (Moreover, if $\nt$ is a maximal nilpotent subalgebra of a simple finite-dimensional Lie algebra then the centrally generated ideals of $\Ua(\nt)$ can be described explicitly \cite{IgnatyevPenkov1}.) 
Similarly, almost all primitive Poisson ideals of $\Sa(\nt)$ are generated as ideals by their intersections with the Poisson center $Y(\nt)$ of $\Sa(\nt)$, see Theorem~\ref{theo:almost_all_Poisson_cent_gen_fin_dim}. 
This result is a form of a ``generic reducibility'' of fibers of maps between algebraic varieties in characteristic 0.   
In our second main result we generalize the latter facts to the case of nil-Dynkin algebras, see Theorem~\ref{theo:all_Poisson_cent_gen_ifd}.

\theobn{Let $\nt$ be a nil-Dynkin algebra.There exists an open dense \textup(with respect to the countable-Zariski topology defined in Subsection~\textup{\ref{sst:pro_varieties}}\textup) subset of $\nt^*$ such that $I(\lambda)$ is generated as an ideal by its intersection with the Poisson center $Y(\nt)$ of~$\Sa(\nt)$ for each $\lambda$ from this subset.}

Using this theorem, we establish two criteria for a primitive Poisson ideal $I(\lambda)$ to be nonzero: one with a hint of linear algebra and another one with a hint of commutative algebra, see Theorem~\ref{theo:non_trivial_ifd}. We consider these criteria as our third main result. In~\cite{IgnatyevPenkov1}, the explicit description of $Z(\nt)$ and $Y(\nt)$ was given. Of course, if $Y(\nt)\neq\Cp$ then each $I(\lambda)$ is nonzero, so we focus on the case when $Y(\nt)=\Cp$. 
The first criterion says that there exists an explicitly described countable collection $\{ \tilde{}\Xi_k\}_{k}$ of countable collections of polynomials from $\Sa(\nt)$ such that $I(\lambda)\ne0$ if and only if there exists $k$ for which $\lambda(\xi)=0$ for all $\xi\in\Xi_k$. It might be interesting to work out a larger class of Lie algebras in which this fact holds. The second criterion is given in terms of ``minors'' in all cases and it can be considered as a statement of the form ``a certain infinite submatrix of matrix defined by $\lambda$ is of finite nonmaximal rank'', see Example~\ref{exam:nontrivial_ifd}.

Note that it can be easily deduced from \cite{IgnatyevPenkov1} that the center of $\nt$ is zero if and only if $Z(\nt)$ and $Y(\nt)$ equal $\Cp$. Note also that in this case almost all ideals $I(\lambda)$ are zero. This is completely opposite to the finite-dimensional case. Indeed, if $\nt$ is finite-dimensional then $I(\lambda)$ is the annihilator in $S(\nt)$ of the coadjoint orbit of $\lambda$. Hence, the condition $I(\lambda)=0$ means that the corresponding orbit is dense in $\nt^*$, but all coadjoint orbits are proper closed subvarieties of $\nt^*$.

This article can be considered as a part of the project researching the coadjoint representations of infinite-dimensional Lie algebras, see~\cite{PenkovPetukhov1} for the case of limits of simple Lie algebras,~\cite{IgnatyevPenkov1} for the case of
nil-Dynkin Lie algebras, \cite{PetukhovSierra1} for the case of Witt Lie algebra.

The paper is organised as follows. In Section~\ref{sect:preliminaries} we present preliminary facts and results (mostly, for finite-dimensional Lie algebras). In Subsection~\ref{sst:Jacobson} (respectively, \ref{sst:Poisson}), we recollect the definition of the Jacobson (respectively, Poisson) topology on the space of primitive ideals of an associative algebra (respectively, on the space of primitive Poisson ideals of a Poisson algebra) together with a few basic properties of these notions. In Subsection~\ref{sst:orbit_method}, we recall how the orbit method works for finite-dimensional nilpotent Lie algebras. In Subsection~\ref{sst:nilradicals_fin_dim}, we discuss an important example of nilpotent Lie algebras --- nilradicals of Borel subalgebras of simple finite-dimensional Lie algebras. In Subsections~\ref{sst:induction_bimod} and~\ref{sst:embeddings} we prove several auxiliary lemmas about nested finite-dimensional nilpotent Lie algebras. Subsection~\ref{sst:cent_gen_ideals_fin_dim} contains our proof of the fact that almost all primitive (Poisson) ideals are generated by their intersections with the (Poisson) center. 

Section~\ref{sect: 1mr} is devoted to the proof of the first main result. In Subsection~\ref{sst:pro_varieties} we recall the definition of pro-variety together with some basic properties of this notion. In particular, given a countable-dimensional vector space $V$, we define the countable-Zariski topology on its dual space $V^*$ and prove that $V^*$ is irreducible with respect to this topology. In Subsection~\ref{sst:locally_nilp_Lie_algs} we introduce the notion of locally nilpotent Lie algebra $\nt$ and establish a bijection between the radical Poisson ideals of $\Sa(\nt)$ and the radical two-sided ideals of $\Ua(\nt)$, see Proposition~\ref{prop:radical_U_n_S_n_ifd}. Further, given a linear form $\lambda\in\nt^*$, we define the ideal $J(\lambda)$ of $\Ua(\nt)$ and check that it is primitive, see Theorem~\ref{theo:primitive_ifdv}. The proof of the first main result (Theorem~\ref{theo:orbit_method_ifd}) given in Subsection~\ref{sst:prime_theorem} is based on an alternative description of $I(\lambda)$ presented in Subsection~\ref{sst:alternative_J_lambda}. 
In Subsections~\ref{sst:socle} and~\ref{sst: adjg} we introduce socle Lie algebras together with the corresponding ``adjoint'' (pro-)groups to prove that the primitive ideals are in one-to-one correspondence with the coadjoint orbits (Proposition~\ref{Porbmsoc}) and show that the corresponding orbit is closed (Lemma~\ref{Lcadjc}). 
We also give a version of Dixmier--Moeglin equivalence for radical ideals in this subsection (Proposition~\ref{Pdmogsa}). 

Section~\ref{sect:nildynkin} is devoted to nil-Dynkin algebras. In Subsection~\ref{sst:nilradicals_inf_dim} we define nil-Dynkin algebras and, given such an algebra $\nt$, recall the description of $Z(\nt)$ (or, equivalently, of $Y(\nt)$) from \cite{IgnatyevPenkov1}. Using this description, in Subsection~\ref{sst:cent_gen_ideals_ifd} we prove our second main result (Theorem~\ref{theo:all_Poisson_cent_gen_ifd}), which claims that almost all primitive Poisson ideals in $\Sa(\nt)$ are generated by their intersections with $Y(\nt)$. Finally, in Subsection~\ref{sst:critetion_non_trivial_ifd} we prove our third main result (Theorem~\ref{theo:non_trivial_ifd}), two criteria for $I(\lambda)$ to be nonzero discussed above. We also apply these criteria to several different nil-Dynkin algebras to give few pictures which can help to understand the pattern behind this result, see Example~\ref{exam:nontrivial_ifd}.

\textsc{Acknowledgements}. The first author was supported by the Foundation for the Advancement of Theoretical Physics and Mathematics ``BASIS'', grant no. 18--1--7--2--1. The second author was supported by RFBR, grant no. 20--01--00091--a. A part of this work was done during our stay at the Oberwolfach Research Institute for Mathematics (program ``Research in pairs'') and the stay of the first author at the Max Planck Institute for Mathematics in spring 2018. We thank these institutions for their hospitality. 
We thank Stephane Launois, Omar Leon Sanchez and Vladimir Zhgoon for useful discussions.

\sect{Preliminaries}\label{sect:preliminaries}\addcontentsline{toc}{section}{\ref{sect:preliminaries}. Preliminaries}
\sst{Jacobson topology and Jacobson spectrum}\label{sst:Jacobson}\addcontentsline{toc}{subsection}{\ref{sst:Jacobson}. Jacobson topology and Jacobson spectrum}

In this subsection, we briefly recall the notion of Jacobson topology and present ring-theoretical preliminaries needed for the sequel.
A detailed discussion can be found, e.g., in the classical books \cite{Dixmier1}, \cite{Joseph2}, \cite{MCR}.
Throughout the paper the ground field will be the field $\mathbb{C}$ of complex numbers (except Remark~\ref{nota:ideals_in_Un}). Let $A$ be an associative algebra (possibly, infinite-dimensional).
Denote by~${\rm Cent}(A)$ the center of $A$.

\defi{A (two-sided) ideal $J$ of $A$ is \emph{prime} if $J\neq A$ and, given two-sided ideals $J', J''$ of $A$ with $J', J''\not\subset J$, one has $J'J''\not\subset J$. An ideal $J$ is called \emph{completely prime} if there are no zero divisors in the quotient algebra~$A/J$. An ideal $J$ is called \emph{primitive} if it is the annihilator of a simple (left) $A$-module (equivalently, if $J$ is the largest ideal in a maximal left ideal of $A$).}

It is not hard to check that 
\begin{equation}\label{Emprpr}J\text{ is maximal }\implies J\text{ is primitive }\implies J\text{ is prime}.\end{equation}
Note also that $J\text{ is completely prime }\implies J\text{ is prime}$. 

One can attach to $A$ the topological space $\JSpec{A}$.\ The construction of $\JSpec{A}$ is as follows:
\begin{equation*}
\begin{split}
\text{i)}&\text{ the points of $\JSpec{A}$ are the primitive ideals of $A$};\\
\text{ii)}&\text{ to any set of elements $S\subset A$ we attach the subset $X_S$ of $\JSpec{A}$ by putting}
\end{split}
\end{equation*}
\begin{equation*}
X_S=\{J\in\JSpec{A}\mid S\subset J\}.
\end{equation*}
By definition, every closed set of $\JSpec{A}$ has the form $X_S$ for a certain $S\subset A$. We call the space $\JSpec{A}$ the \emph{Jacobson} (or, equivalently, \emph{primitive}) \emph{spectrum} of $A$. If $A$ is commutative, then the primitive ideals of $A$ are nothing but the maximal ideals of $A$, so $\JSpec{A}$ is the usual Zariski maximal spectrum.

Let $S$ be a subset of $A$. Put $$\sqrt S=\bigcap_{J\in X_S}J.$$
It is clear that $S\subset \sqrt S$, $X_S=X_{\sqrt S}$ and that $\sqrt S$ is a two-sided ideal of $A$.
Ideal $\sqrt S$ is called {\it the Jacobson radical} of $S$.

\defi{We say that an ideal $I$ of $A$ is {\it radical} if $\sqrt I=I$. Evidently, an ideal is radical if and only if it is the intersection of a family of primitive ideals.}
Obviously, $\sqrt{S}$ is a radical ideal for any $S\subset A$ and the closed subsets of $\JSpec{A}$ can be identified with the radical ideals of $A$. 
This identification reverses the inclusions.

If $A$ is a commutative algebra then we have
$$f\in\sqrt S\text{ if and only if }f^k\in (S)\text{ for some }k\in\Zp_{>0},$$
where $(S)$ is the ideal generated by $S$ (in other words, $f^k$ belongs to the ideal generated by $S$ for certain $k>0$). It is well known that, for an arbitrary finite- or countable-dimensional  associative algebra $A$ (not necessarily commutative), one has
\begin{equation}
f\in\sqrt S\text{ if and only if, for any $g\in A$, }(fg)^k\in (S)\text{ for some }k\in\Zp_{>0},\label{formula:radical_ideal_elements}
\end{equation}
where $(S)$ is the two-sided ideal generated by $S$, see, e.g., \cite[Lemma 4.4]{PenkovPetukhov1}. This means that, given $g\in A$, $(fg)^k$ belongs to the two-sided ideal generated by $S$ for certain $k>0$.
\okr{Lemma}{\label{Lintrad} Let $A$ be an associative algebra and $I, J$ be two-sided ideals of $A$. Then $$\sqrt{IJ}=\sqrt{I\cap J}.$$}
\proof{}{ Pick a primitive ideal $P$ of $A$. It is enough to show that the following conditions are equivalent:
\begin{equation*}
\begin{split}
&\text{\textup{i)} $I\subset P$ or $J\subset P$};\\
&\text{\textup{ii)} $I\cap J\subset P$};\\
&\text{\textup{iii)} $IJ\subset P$}.\\
\end{split}
\end{equation*}
It is clear that $(\text{i})\implies(\text{ii})\implies(\text{iii})$. The implication $(\text{i})\implies(\text{iii})$ follows from~(\ref{Emprpr}).}
The next lemma is a useful and straightforward generalization of Schur's lemma.
\okr{Lemma}{\label{Lschur} Let $A$ be a finite- or countable-dimensional associative algebra, and $J$ be a primitive ideal of $A$. Then the center of $A/J$ consists of scalars.}
\proof{}{ It is easy to deduce this statement from~\cite[Corollary~1.8, Chapter 9]{MCR}.}
\okr{Corollary}{\label{Lschurc} Let $A$ be a a finite- or countable-dimensional commutative associative algebra, and $M$ be a maximal ideal of $A$. Then $\dim_\Cp A/M=1$.}

\sst{Poisson topology and Poisson spectrum}\label{sst:Poisson}\addcontentsline{toc}{subsection}{\ref{sst:Poisson}. Poisson topology and Poisson spectrum}
Here we recall basic facts and notions for Poisson algebras, see, f.e., \cite{BrownGordon1}, \cite{Goodearl1} for the details. Let $\Au$ be a commutative (possibly, infinitely generated) algebra endowed with a polylinear skew-symmetric map
$$\{\cdot, \cdot\}\colon \Au\times \Au\to \Au.$$
The pair $(\Au, \{\cdot, \cdot\})$ is a {\it Poisson algebra} if $\{\cdot, \cdot\}$ satisfies the \emph{Jacobi identity}
$$\{x, \{y, z\}\}=\{\{x, y\}, z\}+\{y, \{x, z\}\}\text{ for all }x,y,z\in \Au,$$
and the \emph{Leibnitz rule}
$$\{x, yz\}=\{x, y\}z+y\{x, z\}\text{ for all }x,y,z\in \Au.$$
In this case $\{\cdot,\cdot\}$ is called the \emph{Poisson bracket} on $\Au$. Denote by ${\rm PCent}(\Au)=\{f\in \Au\mid \{f, \Au\}=0\}$ the {\it Poisson center} of $\Au$.
\defi{We say that an ideal $I$ of $\Au$ is {\it Poisson} if $\{\Au, I\}\subset I$.
A Poisson ideal $I$ of $\Au$ is called {\it primitive} if there exists a maximal ideal $I'$ of $\Au$ such that $I$ is the largest Poisson ideal in~$I'$. (For example, any maximal Poisson ideal is primitive.) }
It is well known that a primitive Poisson ideal is prime, see, f.e.,~\cite[Proposition~1.4]{SQ}. One can attach to $\Au$ the topological space $\PSpec{\Au}$.
The construction of ${\PSpec{\Au}}$ is similar to the construction of the Jacobson spectrum:
\begin{equation*}
\begin{split}
\text{i)}&\text{ the points of $\PSpec{\Au}$ are the primitive Poisson ideals of $\Au$};\\
\text{ii)}&\text{ to any set of elements $S\subset \Au$ we attach the subset $Z_S$ of $\PSpec{\Au}$ by putting}
\end{split}
\end{equation*}
\begin{equation*}
Z_S=\{I\in\PSpec{\Au}\mid S\subset I\}.
\end{equation*}
By definition, every closed set of $\PSpec{\Au}$ has the form $Z_S$ for a certain $S\subset \Au$.
We call the space $\PSpec{\Au}$ the \emph{Poisson} \emph{spectrum} of $\Au$ and the corresponding topology is called the \emph{Poisson topology}.
It is a subspace of the usual Zariski spectrum.

Let $S$ be a subset of $\Au$.
Put $$\sqrt[\Po] S=\bigcap_{I\in Z_S}I.$$
It is clear that $S\subset \sqrt[\Po] S$, $Z_S=Z_{\sqrt[\Po] S}$ and that $\sqrt[\Po] S$ is a Poisson ideal of $\Au$.
The ideal $\sqrt[\Po] S$ is called the {\it Poisson radical} of $S$.
We say that a Poisson ideal $I$ of $\Au$ is {\it radical} if $\sqrt[\Po] I=I$.
It turns out that, for any Poisson ideal $I$ of $\Au$, one has $\sqrt I=\sqrt[\Po] I$. It is clear that the closed subsets of $\PSpec{\Au}$ can be identified with the radical Poisson ideals of $\Au$. This identification reverses the inclusions.

\sst{The orbit method}\label{sst:orbit_method}\addcontentsline{toc}{subsection}{\ref{sst:orbit_method}. The orbit method}

Kirillov's orbit method appears in a wide variety of contexts in representation theory. In this subsection we describe the algebraic version of the orbit method for Lie algebras, see, e.g., the classical Dixmier's book \cite{Dixmier1} for the details. Let $\nt$ be a finite-dimensional Lie algebra. Its symmetric algebra and its universal enveloping algebra are denoted by $\Sa(\nt)$ and $\Ua(\nt)$ respectively. Observe that $\Sa(\nt)$ is a Poisson algebra with respect to the Poisson bracket defined by
$$\{x,y\}=[x,y],~x,y\in\nt.$$
We will now recall how the orbit method works for nilpotent $\nt$, i.e., we will establish a homeomorphism between topological spaces $\JSpec{\Ua(\nt)}$ and $\PSpec{\Sa(\nt)}$ in this case.

For the rest of the subsection we assume that $\nt$ is nilpotent. (In fact, the orbit method also works for solvable Lie algebras, but we will consider only nilpotent case.) It is well known that there exists a unique (up to isomorphism) unipotent algebraic group $N$ such that $\nt$ is the Lie algebra of $N$. We will write $N=\Exp\nt$. Let $\nt^*$ be the dual space of $\nt$. The group $N$ acts on $\nt$ by the adjoint action; the dual action of $N$ on $\nt^*$ is called \emph{coadjoint}. We will denote the result of this action by $g.\lambda$, $g\in N$, $\lambda\in\nt^*$. Given $\lambda\in\nt^*$, we denote by $N.\lambda$ its coadjoint orbit.

To each linear form $\lambda\in\nt^*$ one can assign a bilinear form $\beta_{\lambda}$ on $\nt$ by putting $$\beta_{\lambda}(x,y)=\lambda([x,y]).$$ A subalgebra $\pt\subseteq\nt$ is called a~\emph{polarization} of $\nt$ at $\lambda$ if it is a maximal $\beta_{\lambda}$-isotropic subspace. There is a nice construction of polarizations due to M. Vergne. Namely, let
\begin{equation}
\nt_1\subset\nt_2\subset\label{formula:Vergne_flag}\ldots\subset\nt_k=\nt,~\dim\nt=k,~\dim\nt_i=i\text{ for }1\leq i\leq k,
\end{equation}
be a complete flag of ideals of $\nt$ (clearly, $\nt$ admits such a flag). Denote by $\rt_i$ the kernel of the restriction of $\beta_{\lambda}$ to $\nt_i$, and put $$\pt=\sum_{i=1}^k\rt_i.$$ Then $\pt$ is a polarization of $\nt$ at $\lambda$ \cite{Vergne}. (Note that if all $\nt_i$'s in (\ref{formula:Vergne_flag}) are just subalgebras of $\nt$, then $\pt$ is a maximal $\beta_{\lambda}$-isotropic subspace of $\nt$, but, possibly, not a subalgebra.) Denote by $\Po(\lambda)$ the set of all polarizations of $\nt$ at $\lambda$, and pick $\pt\in\Po(\lambda)$.

Denote by $L_{\lambda}$ the one-dimensional representation of $\pt$ defined by $$\pt\ni x\mapsto\lambda(x)\in\Cp.$$
Further denote by~$V=V(\nt,\pt,\lambda)$ the induced representation of $\nt$, i.e., $$V=\Ua(\nt)\otimes_{U(\pt)} L_\lambda.$$

It turns out that $J(\lambda)=\Ann{V}$ (the annihilator of $V$ in $\Ua(\nt)$) depends only on $\lambda$, not on the choice of $\pt$. Moreover, if $\pt$ is obtained by Vergne's construction, then $V$ is simple, so $J(\lambda)$ is a primitive two-sided ideal of $\Ua(\nt)$. Further, $J(\lambda)=J(\mu)$ if and only if the coadjoint $N$-orbits of $\lambda$ and $\mu$ coincide. Finally, the \emph{Dixmier map} $$\Du\colon\nt^*\to{\rm JSpec}(\Ua(\nt)),~\lambda\mapsto J(\lambda),$$ induces a homeomorphism between the topological spaces $$\nt^*/N\approx{\rm JSpec}(\Ua(\nt)),$$ see \cite{Dixmier1}, \cite{BorhoGabrielRentschler}. Here the space of coadjoint orbits $\nt^*/N$ is endowed with the quotient topology derived from the Zariski topology on $\nt^*$.

On the other hand, radical Poisson ideals in $\Sa(\nt)$ are in one-to-one correspondence with Zariski closed $N$-stable subsets of $\nt^*$, where we interpret $\Sa(\nt)$ as $\Cp[\nt^*]$: to such a subset of $\nt^*$ we attach its annihilator in $\Sa(\nt)$. Furthermore, maximal Poisson ideals correspond to minimal closed $N$-stable subsets. Since the group $N$ is unipotent, all $N$-orbits on $\nt^*$ are closed, which means that each primitive Poisson ideal of $\Sa(\nt)$ is in fact maximal Poisson. It follows that $\PSpec{S(\nt)}$ and $\nt^*/N$ are homeomorphic.

More precisely, given $\lambda\in\nt^*$, one can consider the maximal ideal $I_{\lambda}$ of $\Sa(\nt)$ consisting of all polynomials from $\Sa(\nt)$ vanishing at $\lambda$. There is the unique largest Poisson ideal $I(\lambda)$ inside $I_{\lambda}$, see Subsection~\ref{sst:Poisson}. Clearly, $I(\lambda)$ is exactly the radical Poisson ideal of $\Sa(\nt)$ annihilating the coadjoint $N$-orbit of $\lambda$. In particular, $I(\lambda)=I(\mu)$ if and only if $\lambda$ and $\mu$ belong to the same $N$-orbit. Thus, the map $\lambda\mapsto I(\lambda)$ induces the above homeomorphism $\nt^*/N\approx\PSpec{\Sa(\nt)}$. Combining it with the homeomorphism between $\nt^*/N$ and $\JSpec{\Ua(\nt)}$, we obtain a required homeomorphism from $\PSpec{\Sa(\nt)}$ to $\JSpec(\Ua(\nt))$. Note that the latter homeomorphism preserves inclusions, because both of the homeomorphisms from $\nt^*/N$ to $\PSpec{\Sa(\nt)}$ and to $\JSpec{\Ua(\nt)}$ reverse inclusions. This implies that there is a one-to-one correspondence between radical ideals of~${\rm U}(\nt)$ and radical Poisson ideals of~${\rm S}(\nt)$.

\deno{\label{Nij} If $J$ is a radical two-sided ideal of~${\rm U}(\nt)$, we denote by $I(J)$ the corresponding radical Poisson ideal of $\Sa(\nt)$.
If $I$ is a radical Poisson ideal of~${\rm S}(\nt)$ we denote by $J(I)$ the corresponding radical two-sided ideal of $\Ua(\nt)$.}
Note that $I\subset\Sa(\nt)$ is prime if and only if $J(I)\subset\Ua(\nt)$ is prime~\cite[Proposition~6.3.5]{Dixmier1}.

\nota{In this remark the ground field can be an arbitrary field of zero characteristic. For $\Ua(\nt)$, one can say more about two-sided ideals. Namely,\label{nota:ideals_in_Un} recall that the \emph{Weyl algebra} $\Au_s$ of $2s$ variables is the unital associative algebra with generators $p_i$,~$q_i$ for $1\leq i\leq s$, and relations $[p_i,q_i]=1$, $[p_i,q_j]=0$ for $i\neq j$, $[p_i,p_j]=[q_i,q_j]=0$ for all $i,~j$. Now, if $J$ is a two-sided ideal of the enveloping algebra~$\Ua(\nt)$, then the following conditions are equivalent \cite[Proposition~4.7.4, Theorem 4.7.9]{Dixmier1}:
\begin{itemize}
\item $J$ is primitive;
\item $J$ is maximal;
\item the center of $\Ua(\nt)/J$ is trivial;
\item $U(\nt)/J$ is isomorphic to the Weyl algebra $\Au_s$.
\end{itemize}
(Here $2s$ is the dimension of the orbit of $\lambda$ given $J=J(\lambda)$.) An ideal in $\Ua(\nt)$ is prime if and only if it is completely prime~\cite[Theorem~3.7.2]{Dixmier1}. Furthermore, for any ideal $J$ there exist finitely many minimal pairwise distinct prime ideals $J_1$, $\ldots$, $J_r$ of $\Ua(\nt)$ containing $J$.~None of $J_i$ contains the intersection of the others, and $\sqrt{J}=J_1\cap\ldots\cap J_r$, see~\cite[Proposition 3.1.10]{Dixmier1}. It follows that any prime ideal is radical and the intersection of any family of radical ideals is radical.
The straightforward analogues of the above statements hold for the radical Poisson ideals~of~$\Sa(\nt)$.}

\nota{Here we collect \label{nota:polar_n_m} some facts from linear algebra about polarizations. First, note that, given a linear form $\lambda\in\nt$ and $\pt\in\Po(\lambda)$, one has
$$\dim\pt=\dim\nt-\rk\beta_{\lambda}/2=\rk\beta_{\lambda}/2+\dim\Ker\beta_{\lambda},$$
where $\Ker\beta_{\lambda}=\{x\in\nt\mid\beta_{\lambda}(x,\nt)=0\}$ is the kernel of $\beta_{\lambda}$, and $\rk\beta_{\lambda}$ is the rank of the form $\beta_{\lambda}$.

Next, assume that $\nt$ is a subalgebra of a nilpotent finite-dimensional Lie algebra $\mt$ of codimension one (then $\nt$ is in fact an ideal of $\mt$). Pick a linear form $\lambda\in\mt^*$. Then the following alternative occurs: either $\rk\restr{\beta_{\lambda}}{\nt}=\rk\beta_{\lambda}-2$ (in this case each polarization of $\nt$ at $\restr{\lambda}{\nt}$ is in fact a polarization of $\mt$ at $\lambda$), or $\rk\restr{\beta_{\lambda}}{\nt}=\rk\beta_{\lambda}$. In the latter case, for any $\pt\in\Po(\lambda)$, we have that $\pt\cap\nt\in\Po(\restr{\lambda}{\nt})$. Indeed, $\pt\cap\nt$ is a subalgebra of $\nt$ and a $\restr{\beta_{\lambda}}{\nt}$-isotropic subspace of $\nt$. Further, $\dim(\pt\cap\nt)\geq\dim\pt-1$, but the dimension of a maximal $\restr{\beta_{\lambda}}{\nt}$-isotropic subspace of $\nt$ equals $\dim\nt-\dfrac{1}{2}\rk\restr{\beta_{\lambda}}{\nt}=\dim\mt-1-\dfrac{1}{2}\rk\beta_{\lambda}=\dim\pt-1$. Hence $\pt\cap\nt$ is a maximal $\restr{\beta_{\lambda}}{\nt}$-isotropic subspace of $\nt$.}

\exam{(Heisenberg Lie algebra)\label{Ehla} Pick $n\ge1$. Let $\nt=\hei_n(\Cp)$ be the Lie algebra with generators $z$, $x_i$, $y_i$, $1\le i\le n$, and relations
$$[x_i,y_i]=z\text{ for all }i,~[x_i,y_j]=0\text{ for }i\neq j,~[x_i,z]=[y_i,z]=0\text{ for all }i.$$
We call $\nt$ the \emph{Heisenberg Lie algebra}. There are two classes of coadjoint orbits on $\nt^*$:

i) if $\lambda(z)=\alpha\neq0$ then $N.\lambda=\{\mu\in\nt^*\mid\mu(z)=\alpha\}$.

ii) if $\lambda(z)=0$ then $N.\lambda=\{\lambda\}$.

In a similar way, there are two classes of primitive ideals of ${\rm U}(\nt)$ and ${\rm S}(\nt)$:

i) every $\alpha\in\Cp^{\times}=\Cp\setminus\{0\}$ defines the two-sided (respectively, Poisson) ideal  of ${\rm U}(\nt)$ (respectively, of~${\rm S}(\nt)$) generated by $z-\alpha$. 
It is easy to verify that the quotient of $\Ua(\nt)$ by this ideal is a Weyl algebra and hence is simple (a similar argument is applicable to the Poisson side). Thus, this ideal is maximal and hence primitive (respectively, Poisson primitive).

ii) every $\lambda\in\nt^*$ with $\lambda(z)=0$ defines the ideal $J(\lambda)$ generated by $$x_i-\lambda(x_i),~y_i-\lambda(y_i),~1\le i\le n,~\text{~and~}~z.$$ 
It is easy to verify that the quotient by $J(\lambda)$ is isomorphic to $\mathbb C$. Thus $J(\lambda)$ is maximal and hence primitive.

The bijection between coadjoint orbits and primitive ideals is clear from these descriptions.}\newpage

\sst{Nilradicals of Borel subalgebras}\label{sst:nilradicals_fin_dim}\addcontentsline{toc}{subsection}{\ref{sst:nilradicals_fin_dim}. Nilradicals of Borel subalgebras}

In this subsection, we briefly recall definitions of classical finite-dimensional simple Lie algebras and fix notation for the nilradicals of their Borel subalgebras. These nilradicals (and their infinite-dimensional analogues defined in Subsection~\ref{sst:nilradicals_inf_dim}) provide one of the main examples of nilpotent Lie algebras for our purposes.

Pick $n\in\Zp_{>0}$. Let $\gt$ denote one of the Lie algebras $\slt_n(\Cp)$, $\sot_{2n}(\Cp)$, $\sot_{2n+1}(\Cp)$ or $\spt_{2n}(\Cp)$. The algebra $\sot_{2n}(\Cp)$ (respectively, $\sot_{2n+1}(\Cp)$ and $\spt_{2n}(\Cp)$) is realized as the sub\-al\-gebra of $\slt_{2n}(\Cp)$ (respectively, $\slt_{2n+1}(\Cp)$ and $\slt_{2n}(\Cp)$) consisting of all~$x$ such that $$\beta(u,xv)+\beta(xu,v)=0$$ for all $u,v$ in $\Cp^{2n}$ (respectively, in $\Cp^{2n+1}$ and $\Cp^{2n}$), where
\begin{equation*}
\beta(u,v)=\begin{cases}
\sum\nolimits_{i=1}^n(u_iv_{-i}+u_{-i}v_i)&\text{for }\sot_{2n}(\Cp),\\
u_0v_0+\sum\nolimits_{i=1}^n(u_iv_{-i}+u_{-i}v_i)&\text{for }\sot_{2n+1}(\Cp),\\
\sum\nolimits_{i=1}^n(u_iv_{-i}-u_{-i}v_i)&\text{for }\spt_{2n}(\Cp).
\end{cases}
\end{equation*}
Here for $\sot_{2n}(\Cp)$ (respectively, for $\sot_{2n+1}$ and $\spt_{2n}(\Cp))$ we denote by $e_1,\ldots,e_n,e_{-n},\ldots,e_{-1}$ (res\-pec\-tively, by $e_1,\ldots,e_n,e_0,e_{-n},\ldots,e_{-1}$ and $e_1,\ldots,e_n,e_{-n},\ldots,e_{-1}$) the standard basis of $\Cp^{2n}$ (res\-pec\-tively, of $\Cp^{2n+1}$ and $\Cp^{2n}$), and by $x_i$ the coordinate of a vector $x$ corresponding to $e_i$.

The set of all diagonal matrices from $\gt$ is a Cartan subalgebra of $\gt$; we denote it by $\htt$. Let $\Phi$ be the root system of $\gt$ with respect to $\htt$. Note that $\Phi$ is of type $A_{n-1}$ (respectively, $D_n$, $B_n$ and $C_n$) for $\slt_n(\Cp)$ (respectively, for $\sot_{2n}(\Cp)$, $\sot_{2n+1}(\Cp)$ and $\spt_{2n}(\Cp)$). The set of all upper-triangular matrices from $\gt$ is a Borel subalgebra of $\gt$ containing $\htt$; we denote it by $\bt$. Let $\Phi^+$ be the set of positive roots with respect to $\bt$. As usual, we identify $\Phi^+$ with the following subset of $\Rp^n$:
\begin{equation}\predisplaypenalty=0
\begin{split}
A_{n-1}^+&=\{\epsi_i-\epsi_j,~1\leq i<j\leq n\},\\
B_n^+&=\{\epsi_i-\epsi_j,~1\leq i<j\leq n\}\cup\{\epsi_i+\epsi_j,~1\leq i<j\leq n\}\cup\{\epsi_i,~1\leq i\leq n\},\\\label{formula:root_basis_nilradical_fin_dim}
C_n^+&=\{\epsi_i-\epsi_j,~1\leq i<j\leq n\}\cup\{\epsi_i+\epsi_j,~1\leq i<j\leq n\}\cup\{2\epsi_i,~1\leq i\leq n\},\\
D_n^+&=\{\epsi_i-\epsi_j,~1\leq i<j\leq n\}\cup\{\epsi_i+\epsi_j,~1\leq i<j\leq n\}.\\
\end{split}
\end{equation}
Here $\{\epsi_i\}_{i=1}^n$ is the standard basis of $\Rp^n$.

Denote by $\nt$ the algebra of all strictly upper-triangular matrices from $\gt$. Then $\nt$ has a basis consisting of root vectors $e_{\alpha}$, $\alpha\in\Phi^+$, where
\begin{equation*}\predisplaypenalty=0
\begin{split}
e_{\epsi_i}&=\sqrt{2}(e_{0,i}-e_{-i,0}),~e_{2\epsi_i}=e_{i,-i},\\
e_{\epsi_i-\epsi_j}&=\begin{cases}
e_{i,j}&\text{for }A_{n-1},\\
e_{i,j}-e_{-j,-i}&\text{for }B_n,~C_n\text{ and }D_n,
\end{cases}\\
e_{\epsi_i+\epsi_j}&=\begin{cases}
e_{i,-j}-e_{j,-i}&\text{for }B_n\text{ and }D_n,\\
e_{i,-j}+e_{j,-i}&\text{for }C_n,
\end{cases}
\end{split}
\end{equation*}
and $e_{i,j}$ are the usual elementary matrices. For $\sot_{2n}(\Cp)$ (respectively, for $\sot_{2n+1}(\Cp)$ and $\spt_{2n}(\Cp)$) we index the rows (from left to right) and the columns (from top to bottom) of matrices by the numbers $1,\ldots,n,-n,\ldots,-1$ (respectively, by the numbers $1,\ldots,n,0,-n,\ldots,-1$ and $1,\ldots,n,-n,\ldots,-1$). Note that $\gt=\htt\oplus\nt\oplus\nt_-$, where $\nt_-=\langle e_{-\alpha},~\alpha\in\Phi^+\rangle_{\Cp}$, and, by definition, $e_{-\alpha}=e_{\alpha}^T$. (The superscript~$T$ always indicates matrix transposition.) The set $\{e_{\alpha},~\alpha\in\Phi\}$ can be uniquely extended to a Chevalley basis of $\gt$. Clearly, $\nt$ is the nilradical of the Borel subalgebra $\bt$.

Since $\nt$ is nilpotent, the space $\JSpec{\Ua(\nt)}$ can be described in terms of the coadjoint orbits on $\nt^*$. But a classification of the coadjoint orbits on $\nt^*$ itself is an extremely hard problem. For instance, for $A_{n-1}$ a complete classification is known only for $n\leq8$ \cite{IgnatyevPanov1}. On the other hand, almost all orbits on $\nt^*$ has maximal possible dimension (such orbits are called \emph{regular}), and for $A_{n-1}$ all regular orbits were in fact described in the pioneering Kirillov's work \cite{Kirillov1} on the orbit method in 1962. Let us briefly recall this description.

\exam{Let $\Phi=A_{n-1}$. We put $m=[n/2]$ and\label{exam:regular_orbits_A_n} define the regular functions $\xi_i$, $1\leq i\leq m$, on $\nt^*$ by the following rule:
\begin{equation*}
\xi_i(\lambda)=\begin{vmatrix}
\lambda(e_{1,n-i+1})&\ldots&\lambda(e_{1,n-1})&\lambda(e_{1,n})\\
\lambda(e_{2,n-i+1})&\ldots&\lambda(e_{2,n-1})&\lambda(e_{2,n})\\
\vdots&\iddots&\vdots&\vdots\\
\lambda(e_{i,n-i+1})&\ldots&\lambda(e_{i,n-1})&\lambda(e_{i,n})\\
\end{vmatrix},~\lambda\in\nt^*.
\end{equation*}
It is easy to see that $\xi_i\in\PCent{\Sa(\nt)}$ for all $i$, so $\xi_i(\mu)=\xi_i(\lambda)=c_i$ for all $\mu\in N.\lambda$ (for certain $c_i\in\Cp$). It turns out that the coadjoint $N$-orbit of $\lambda$ has maximal possible dimension $$\dim N.\lambda=2((n-2)+(n-4)+\ldots)$$ if and only if $c_i\neq0$ for all $i$ for odd $n$, and $c_i\neq0$ for all $i<m$ for even $n$. In this case, $I(\lambda)$ is generated by $\xi_i-c_i$, $1\leq i\leq m$.}


\sst{Induction through bimodules}\label{sst:induction_bimod}\addcontentsline{toc}{subsection}{\ref{sst:induction_bimod}. Induction through bimodules}
In this subsection we introduce a series of bimodules which will allow us to construct infinite-dimensional simple modules. 

\okr{Theorem}{\label{Thm:bim}Let $\mt$ be a finite-dimensional nilpotent Lie algebra, $\nt$ be a subalgebra of $\mt$, $\lambda\in\mt^*$ be a linear form.
Then there exists a $({\rm U}(\mt)/J(\lambda)){-}({\rm U}(\nt)/J(\restr{\lambda}{\nt}))$-bimodule ${}_{\mt}F_{\nt}$ together with a ${\rm U}(\nt){-}{\rm U}(\nt)$ morphism $$\phi\colon {\rm U}(\nt)\to {}_{\mt}F_{\nt}$$ such that\\
\textup{i)} the functor $F\colon V\mapsto {}_{\mt}F_{\nt}\otimes_{{\rm U}(\nt)}V$ from the category of ${\rm U}(\nt)/J(\restr{\lambda}{\nt})$-modules to the category of ${\rm U}(\mt)/J(\lambda)$-modules is exact and sends simple modules to simple modules\textup{;}\\
\textup{ii)} the induced map $\phi\otimes\id\colon V\cong {\rm U}(\nt)\otimes_{{\rm U}(\nt)}V\to F(V)$ is injective for any $({\rm U}(\nt)/J(\restr{\lambda}{\nt}))$-module $V$.}
Our proof of Theorem~\ref{Thm:bim} is based on the following lemma.
\okr{Lemma}{\label{L:bim}The statement of Theorem~\textup{\ref{Thm:bim}} holds under the assumption that $\dim\mt=\dim\nt+1$.}
\vspace{-0.2cm}\proof{ of Theorem~\ref{Thm:bim}}{ 
Pick a sequence of Lie algebras $$\nt=\nt_0\subset\nt_2\subset...\subset\nt_s=\mt$$ such that $\dim \nt_i=\dim\nt_0+i.$ 
Let ${}_{\nt_{i+1}}F_{\nt_i}$ be the module defined by Lemma~\ref{L:bim} for $\nt_i, \nt_{i+1}, \restr{\lambda}{\nt_{i+1}}$ with $0\le i<s$. 
Set $${}_{\mt}F_{\nt}:=({}_{\nt_{s}}F_{\nt_{s-1}})\otimes_{{\rm U}(\nt_{s-1})}...\otimes_{{\rm U}(\nt_{1})}({}_{\nt_1}F_{\nt_0})$$
and set $\phi$ to be the composition of the maps $\phi_i$ for all the above pairs. Then it is clear that ${}_\mt F_\nt$ and~$\phi$ satisfy all the conditions of Theorem~\ref{Thm:bim}. 
}
\proof{ of Lemma~\ref{L:bim}}{Consider $\beta_\lambda$ together with the restriction $\restr{\beta_{\lambda}}{\nt}$ of $\beta_\lambda$ to $\nt$. 
Then, by Remark~\ref{nota:polar_n_m}, either (1) or (2) holds:

(1) ${\rm rk}(\restr{\beta_{\lambda}}{\nt})={\rm rk}\beta_\lambda;$

(2) ${\rm rk}(\restr{\beta_{\lambda}}{\nt})+2={\rm rk}\beta_\lambda$.\par
\noindent If (1) holds then \cite[Lemma 6.5.6]{Dixmier1} implies that
$${\rm U}(\mt)/J(\lambda)={\rm U}(\nt)/(J(\lambda)\cap{\rm U}(\nt))\cong{\rm U}(\nt)/J(\restr{\lambda}{\nt}).$$
The algebra ${\rm U}(\mt)/J(\lambda)$ is a ${\rm U}(\mt){-}{\rm U}(\mt)$-bimodule in a natural way. Thus ${\rm U}(\mt)/J(\lambda)$ is also a ${\rm U}(\mt){-}{\rm U}(\nt)$-bimodule.
Put $_{\mt}F_\nt={\rm U}(\mt)/J(\lambda)$ and set $\phi$ to be the natural map ${\rm  U}(\nt)\to{\rm U}(\mt)/J(\lambda)$.
It is clear that $_{\mt}F_\nt$ satisfies the required conditions.

Assume that (2) holds.
Then we put $_{\mt}F_{\nt}={\rm U}(\mt)\otimes_{{\rm U}(\nt)}({\rm U}(\nt)/J(\restr{\lambda}{\nt}))$. 
The natural map $${\rm U}(\nt)\to {\rm U}(\mt)\otimes_{{\rm U}(\nt)}({\rm U}(\nt)/J(\lambda|_{\nt})),~a\mapsto a\otimes 1,$$ defines the desired map $\phi$. 

The functor $F\colon V\mapsto {}_{\mt}F_{\nt}\otimes_{{\rm U}(\nt)}V$ is clearly exact and sends simple objects to simple objects by \cite[5.3]{Dixmier1}. 
As the next step we show that all modules in the image of $F$ are annihilated by $J(\lambda)$. 
Thanks to~\cite[5.2.6, 5.1.7]{Dixmier1} we have that the annihilators of all modules in the image of $F$ are all the same. 
We wish to show that this annihilator is $J(\lambda)$.

Indeed, let $\frak p$ be a polarization of $\nt$ at $\restr{\lambda}{\nt}$. 
Then $\dim \frak p+\frac12\rk(\restr{\beta_\lambda}{\nt})=\dim \nt$ and hence $$\dim \frak p+\frac12 \rk\beta_\lambda=\dim \mt.$$ 
Therefore $\frak p$ is a polarization of $\mt$ at $\beta_\lambda$.

Thus $M:={\rm U}(\nt)\otimes_{{\rm U}(\frak p)} L_{\lambda}$ is annihilated by $J(\restr{\lambda}{\nt})$ and $F(M)\cong {\rm U}(\mt)\otimes_{{\rm U}(\frak p)} L_{\lambda}$ (recall the definition of $L_{\lambda}$ from Subsection~\ref{sst:orbit_method}). 
The latter module is annihilated by $J(\lambda)$. 
}

\sst{Embeddings of nilpotent Lie algebras}\label{sst:embeddings}\addcontentsline{toc}{subsection}{\ref{sst:embeddings}. Embeddings of nilpotent Lie algebras}

In this subsection, we consider embeddings of finite-dimensional nilpotent Lie algebras in more details. This is needed for the subsequent consideration of locally nilpotent infinite-dimensional Lie algebras. Let $\mt$ be a (finite-dimensional) nilpotent Lie algebra and $\nt$ be a subalgebra of $\mt$. Recall the homeomorphism $\JSpec{\Ua(\mt)}\approx\PSpec{\Sa(\mt)}$ from Subsection~\ref{sst:orbit_method}: $$\Sa(\mt)\supset I\mapsto J(I)\subset\Ua(\mt),~\Ua(\mt)\supset J\mapsto I(J)\subset\Sa(\mt).$$

Let $I$ be a radical Poisson ideal of ${\rm S}(\mt)$.
The intersection ${\rm S}(\nt)\cap I$ is a radical Poisson ideal of~${\rm S}(\nt)$ and it is natural to expect that, in some sense, $$J(I\cap{\rm S}(\nt))\text{ ``$\approx$'' }J(I)\cap{\rm U}(\nt).$$
In this subsection we will prove that the answer to the question $$\text{What does ``$\approx$'' mean here?}$$ is as nice as possible. Precisely, we will prove the following theorem.

\okr{Theorem}{\label{theo:Um_Un} Let $\nt, \mt$ be as above. Then $J(I\cap{\rm S}(\nt)) = J(I)\cap{\rm U}(\nt)$.}

As a first step, we will prove the following lemma.

\lemmp{The ideal $J(I)\cap{\rm U}(\nt)$ is radical.}{The ideal $J(I)$ is radical and hence, by Remark~\ref{nota:ideals_in_Un}, there exist prime ideals $J_1,..., J_s$ of~${\rm U}(\mt)$ such that $$J(I)=J_1\cap\ldots\cap J_s.$$
Next, we have
$$J(I)\cap {\rm U}(\nt)=(J_1\cap{\rm U}(\nt))\cap\ldots\cap(J_s\cap{\rm U}(\nt)).$$
An ideal of ${\rm U}(\mt)$ is prime if and only if it is completely prime (see Remark~\ref{nota:ideals_in_Un}).
The intersection of a completely prime ideal with a subalgebra is completely prime.
An intersection of prime ideals of ${\rm U}(\nt)$ is a radical ideal. This completes the proof.}

Therefore we left to show that $I(J(I)\cap{\rm U}(\nt))=I\cap {\rm S}(\nt)$. We start with a particular case of Theorem~\ref{theo:Um_Un}. Namely, assume that $\dim\mt=\dim\nt+1$.

\propp{Let $\nt, \mt$ be\label{prop:Um_Un_primitive} as above and assume that $\dim\mt=\dim\nt+1$, $I=I(\lambda)$ for some $\lambda\in\mt^*$. Then $$J(I\cap{\rm S}(\nt)) = J(I)\cap{\rm U}(\nt).$$}
{Consider $\beta_\lambda$ together with the restriction $\restr{\beta_{\lambda}}{\nt}$ of $\beta_\lambda$ to $\nt$. Then, by Remark~\ref{nota:polar_n_m}, either (1) or (2) holds:

(1) ${\rm rk}(\restr{\beta_{\lambda}}{\nt})={\rm rk}\beta_\lambda;$

(2) ${\rm rk}(\restr{\beta_{\lambda}}{\nt})+2={\rm rk}\beta_\lambda$.\par
\noindent If (1) holds then \cite[Lemma 6.5.6]{Dixmier1} says that $J(\lambda)\cap{\rm U}(\nt)=J(\restr{\lambda}{\nt})$ and $I(\lambda)\cap{\rm S}(\nt)=I(\restr{\lambda}{\nt})$. This implies the required equality.

Assume that (2) holds. Set $M=\Exp\mt$. Since $\nt$ is an ideal of $\mt$, \cite[Lemma~6.5.1]{Dixmier1} states that
$$J(\lambda)\cap{\rm U}(\nt)=\bigcap_{g\in M}J(g.\restr{\lambda}{\nt}).$$
We left to show that
$$I(\lambda)\cap{\rm S}(\nt)=\bigcap_{g\in M}I(g.\restr{\lambda}{\nt}).$$
This equality is a very straightforward exercise in commutative algebra and therefore the proof is complete.
}

\corop{Let $\nt,~\mt$ be as above and $I$ be a radical\label{coro:Um_Un_radical} Poisson ideal of $\Sa(\mt)$. Assume that $\dim\mt=\dim\nt+1$. Then $$J(I\cap\Sa(\nt)) = J(I)\cap\Ua(\nt).$$}
{Note that $$I=\bigcap_{\lambda\in S}I(\lambda)$$
for some subset $S\subset\mt^*$.
Hence
$$J(I)=\bigcap_{\lambda\in S}J(\lambda).$$
This together with Proposition~\ref{prop:Um_Un_primitive} implies the desired result.}

\proof{ of Theorem~\ref{theo:Um_Un}}{ We recall that there exists a complete flag of subalgebras in $\mt$, i.e., a chain of nested subalgebras $\nt=\nt_1\subset\ldots\subset\nt_k=\mt$ such that $\dim\nt_i=i$ for all $i$ from 1 to $k=\dim\mt$. By Corollary~\ref{coro:Um_Un_radical} we have
$$J(I)\cap{\rm U}(\nt)=(J(I)\cap{\rm U}(\nt_{k-1}))\cap{\rm U}(\nt)=J(I\cap{\rm S}(\nt_{k-1}))\cap{\rm U}(\nt)=$$
$$=(J(I\cap{\rm S}(\nt_{k-1}))\cap{\rm U}(\nt_{k-2}))\cap{\rm U}(\nt)=J(I\cap \Sa(\nt_{k-2}))\cap{\rm U}(\nt)=...=J(I\cap{\rm S}(\nt)).$$ The result follows.}

\sst{Centrally generated ideals}\label{sst:cent_gen_ideals_fin_dim}\addcontentsline{toc}{subsection}{\ref{sst:cent_gen_ideals_fin_dim}. Centrally generated ideals}

This subsection is devoted to the special case of primitive ideals, namely, to the centrally generated ones. Recall that an ideal $J$ of an associative algebra $A$ is called \emph{centrally generated} if it is generated as an ideal by its intersection with the center of $A$. Let $\nt$ be a finite-dimensional nilpotent Lie algebra and $N=\Exp\nt$. Let $P$ be a prime ideal of $\Ua(\nt)$. First, we will prove that almost all primitive ideals of $\Ua(\nt)/P$ are centrally generated. We start from the following observation.

For brevity, denote by $Z(\nt; P)=\Cent{\Ua(\nt)/P}$ the center of the quotient algebra $\Ua(\nt)/P$. 
Since $\Ua(\nt)/P$ is a domain, $Z(\nt; P)$ is also a domain. Let $J$ be a primitive ideal of $\Ua(\nt)/P$, then $J\cap Z(\nt; P)$ is a maximal ideal of $Z(\nt; P)$, because the center of $\Ua(\nt)/J$ is trivial by Remark~\ref{nota:ideals_in_Un}. Note that, in general, $Z(\nt; P)$ is not finitely generated \cite[4.9.20]{Dixmier1}. Nevertheless, one can consider the space $\JSpec{Z(\nt; P)}$, which is in this case nothing but the usual Zariski maximal spectrum of $Z(\nt; P)$, because $Z(\nt; P)$ is commutative. We will denote this topological space by $\MSpec{Z(\nt; P)}$. (Since $Z(\nt; P)$ is a domain, it is an irreducible space.) In particular, to each $f\in Z(\nt; P)$ one can assign a dense open subset of $\MSpec{Z(\nt; P)}$ of the form $$D(f)=\{M\in\MSpec{Z(\nt; P)}\mid f\notin M\}.$$
The following proposition seems to be known to the specialists.

\propp{\label{prop:open_subset_Z_n_cent_gen} There exists an element $f\in Z(\nt; P)$ such that if $M\in D(f)$ then the ideal of~$\Ua(\nt)/P$ generated by $M$ is primitive.}{Denote by $Q$ the field of fractions of $Z(\nt; P)$. 
One can consider the Lie algebra $\nt_Q=\nt\otimes_{\Cp}Q$ over the field $Q$ and its universal enveloping algebra $\Ua(\nt_Q)\cong\Ua(\nt)\otimes_\Cp Q$, which is an associative algebra over the field $Q$.

On the other hand, $\Ua(\nt)/P\otimes_{Z(\nt; P)}Q$ is also an associative algebra over $Q$, and the natural embedding $\nt\hookrightarrow\Ua(\nt)$ induces the following map: 
$$\nt_Q\to\Ua(\nt)/P\otimes_{Z(\nt; P)}Q,~x\otimes\alpha\mapsto\wt x\otimes\alpha,~x\in\nt,~\alpha\in Q,$$ where $\wt x$ denotes the image of $x$ in the quotient algebra $\Ua(\nt)/P$. 
By the universal property of $\Ua(\nt_Q)$, the latter map can be extended to the unique morphism of associative algebras (over $Q$) $$\vfi\colon\Ua(\nt_Q)\to\Ua(\nt)/P\otimes_{Z(\nt; P)}Q,$$ which is clearly surjective (cf. \cite[4.1.3]{Dixmier1}). 
It is evident that $\Cent{\Ua(\nt)/P\otimes_{Z(\nt; P)}Q}=Q$, hence, by Remark~\ref{nota:ideals_in_Un}, the kernel of $\vfi$ is a primitive ideal of $\Ua(\nt_Q)/P$.

Denote this primitive ideal of $\Ua(\nt_Q)/P$ by $\Jo$. Again by Remark~\ref{nota:ideals_in_Un}, $\Ua(\nt)/P\otimes_{Z(\nt; P)}Q$ is isomorphic to a Weyl algebra $\Au_s$ over $Q$ for certain $s$ (as the quotient of $\Ua(\nt_Q)/P$ by the primitive ideal $\Jo$). 
Fix an isomorphism $$\psi\colon\Au_s\to\Ua(\nt)/P\otimes_{Z(\nt; P)}Q.$$ 
Let $p_i$, $q_i$, $1\leq i\leq s$, be the standard generators of $\Au_s$. 
Fix a $\Cp$-basis $e_i$, $1\leq i\leq n=\dim_{\Cp}\nt$, of the Lie algebra $\nt$, then all $\wt e_i\otimes1$, $1\leq i\leq s$, generate $\Ua(\nt)/P\otimes_{Z(\nt; P)}Q$ as a $Q$-algebra. 
Hence, $\psi(p_i)$ and $\psi(q_i)$ can be expressed as polynomials in $\wt e_i\otimes 1$ with coefficients in $Q$, and vice versa.

More precisely, given two $s$-tuples $A=(a_1,~\ldots,~a_s)$, $B=(b_1,~\ldots,~b_s)$ (respectively, an $n$-tuple $C=(c_1,~\ldots,~c_n)$) of nonnegative integers, we put
\begin{equation*}
p^Aq^B=\psi(p_1)^{a_1}\psi(q_1)^{b_1}\ldots\psi(p_s)^{a_s}\psi(q_s)^{b_s},~
e^C=\wt e_1^{c_1}\ldots\wt e_n^{c_n}\otimes1.
\end{equation*}
Then both the sets $\{p^Aq^B\}$ (where $A$ and $B$ run over all $s$-tuples independently) and $\{e^C\}$ (where $C$ runs over all $n$-tuples) generate $\Ua(\nt)/P\otimes_{Z(\nt; P)}Q$ as a vector space over $Q$. Hence, for all possible $i$, there are finite expressions
\begin{equation*}
\psi(p_i)=\sum\nolimits_{C}\alpha_{i,C}e^C,~\psi(q_i)=\sum\nolimits_{C}\beta_{i,C}e^C,~\wt e_i\otimes1=\sum\nolimits_{A,B}\gamma_{i,A,B}p^Aq^B
\end{equation*}
for certain $\alpha_{i,C}$, $\beta_{i,C}$, $\gamma_{i,A,B}\in Q$. Let $f\in Z(\nt)$ be such that $f\alpha_{i,C}$, $f\beta_{i,C}$, $f\gamma_{i,A,B}\in Z(\nt)$ for all possible indices (i.e., $f$ is a ``common denominator'' of these coefficients).

Next, pick a maximal ideal $M$ of $Z(\nt; P)$ for which $f\notin M$ and denote by $$\tau\colon Z(\nt; P)\to Z(\nt; P)/M\cong\Cp$$ the canonical projection. By definition of $f$ one can clearly define the complex numbers $\tau(\alpha_{i,C})$, $\tau(\beta_{i,C})$ and $\tau(\gamma_{i,A,B})$. The linear map $$e_i\mapsto\sum\nolimits_{A,B}\tau(\gamma_{i,A,B})p_1^{a_1}q_1^{b_1}\ldots p_s^{a_s}q_s^{b_s}$$ from $\nt$ to $\Au_s$ induces the associative algebra epimorphism $\Ua(\nt)\to\Au_s$. Its kernel $K$ is a primitive ideal of $\Ua(\nt)/P$ by Remark~\ref{nota:ideals_in_Un}. Note that $M\subset K$.

Now, denote by $J$ the ideal of $\Ua(\nt)/P$ generated by~$M$, then $J\subset K$ and, consequently, $J$ does not contain $f$. Consider the quotient algebra $A=U(\nt)/J$. 
Let $\pi$ be the canonical projection from $\Ua(\nt)$ to $A$. 
Note that $\pi(f)=c$ for certain $c\in\Cp^{\times}$. Consider the elements $\wh p_i$, $\wh q_i$, $1\leq i\leq s$, of $A$ defined by
\begin{equation*}
\wh p_i=c^{-1}\sum\nolimits_{C}\pi\left(f\alpha_{i,C}e^C\right),~\wh q_i=c^{-1}\sum\nolimits_{C}\pi\left(f\beta_{i,C}e^C\right).
\end{equation*}
It is clear that these elements satisfy the defining relations of the Weyl algebra in $2s$ variables. 
On the other hand, for an arbitrary $i$ from 1 to $n$,
\begin{equation*}
\pi(e_i)=c^{-1}\sum\nolimits_{A,B}\pi(f\gamma_{i,A,B})\wh p^A\wh q^B,
\end{equation*}
where $\wh p^A\wh q^B=\wh p_1^{a_1}\wh q_1^{b_1}\ldots\wh p_s^{a_s}\wh q_s^{b_s}$.

It follows that the elements $\wh p_i$ and $\wh q_i$, $1\leq i\leq s$, generate $A$ as an algebra, so $A$ is a quotient of the Weyl algebra $\Au_s$. 
But the Weyl algebra is simple, so $A$ is isomorphic to $\Au_s$. 
Thus, Remark~\ref{nota:ideals_in_Un} implies that the ideal $J$ is primitive, as required.}

The following theorem is the main result of this subsection. 
Recall the definition of the Dixmier map $\Du\colon\nt^*\to\JSpec{\Ua(\nt)}$ from Subsection~\ref{sst:orbit_method}. 
Note that if $V$ is an open $N$-stable subset of $\nt^*$ then $\Du(V)$ is an open subset of $\JSpec{\Ua(\nt)}$, because $\Du$ induces a homeomorphism between the spaces $\nt^*/N$ and $\JSpec{\Ua(\nt)}$. 
Note also that the space $\nt^*/N$ (and, consequently, $\JSpec{\Ua(\nt)}$) is irreducible as a surjective image of the irreducible space $\nt^*$.

\theop{There exists\label{theo:almost_all_cent_gen_fin_dim} an open $N$-stable subset of $\nt^*$ such that $J(\lambda)$ is centrally generated for each $\lambda$ from this subset. In other words\textup, there exists an open \textup(dense\textup) subset of $\JSpec{\Ua(\nt)}$ such that each primitive ideal from this subset is centrally generated.}{For brevity, denote the Poisson center $\PCent{\Sa(\nt)}$ of $\Sa(\nt)$ by $Y(\nt)$. It is well known that the restriction of the canonical \emph{symmetrization map} $$\sigma\colon\Sa(\nt)\to\Ua(\nt),~x^k\mapsto x^k,~x\in\nt,~k\in\Zp_{\geq0},$$ to $Y(\nt)$ is an algebra isomorphism between $Y(\nt)$ and $Z(\nt)$ \cite[Proposition 4.8.12]{Dixmier1}.

Let $f$ be an element from Proposition~\ref{prop:open_subset_Z_n_cent_gen} for $P=\{0\}$. Put $F=\sigma^{-1}(f)\in Y(\nt)$ and set $$V=\{\lambda\in\nt^*\mid F(\lambda)\neq0\}\subset\nt^*.$$ (Here we identify $\Sa(\nt)$ with the algebra $\Cp[\nt^*]$ of regular functions on $\nt^*$.) Clearly, $V$ is an open subset of $\nt^*$; $V$ is $N$-stable because $F$ belongs to the Poisson center of $\Sa(\nt)$ and so is constant on $N$-orbits. We will check that $J(\lambda)$ is centrally generated for all $\lambda\in V$.

Pick a linear form $\lambda\in V$ and recall the definition of the ideal $I(\lambda)$ of $\Sa(\nt)$ from Subsection~\ref{sst:Poisson}. Since $F(\lambda)\neq0$, we have $F\notin I(\lambda)\cap Y(\nt)$. By \cite[6.6.11]{Dixmier1}, $\sigma(I(\lambda)\cap Y(\nt))=J(\lambda)\cap Z(\nt)$. Thus, $f\notin J(\lambda)\cap Z(\nt)$. By Proposition~\ref{prop:open_subset_Z_n_cent_gen}, the ideal $J(\lambda)$ is centrally generated.}

The following theorem is a Poisson analogue of Theorem~\ref{theo:almost_all_cent_gen_fin_dim}.
\mtheop{There exists\label{theo:almost_all_Poisson_cent_gen_fin_dim} an open $N$-stable subset of $\nt^*$ such that $I(\lambda)\in\PSpec{\Sa(\nt)}$ is generated as an ideal by its intersection with the Poisson center $Y(\nt)$ of $\Sa(\nt)$ for each $\lambda$ from this subset.}
{Thanks to Rosenlicht Theorem there exists an open affine subset $U$ of $\nt^*$, an algebraic variety $Q$ and a morphism of algebraic varieties $\phi\colon U\to Q$ such that $U\to Q$ is the geometric quotient $U/N$. 
Without loss of generality we may assume that $Q$ is affine and provides an isomorphism between $\Cp[Q]$ and $\Cp[U]^N$. (Here and below we denote by $\cdot^N$ the set of $N$-invariant vectors.)

Next, recall that every morphism of algebraic varieties in characteristic 0 is smooth on an open subset. 
Thus, there exists $U'\subset U$ such that the restriction $\phi'=\restr{\phi}{U'}\colon U'\to Y$ is smooth. 
The condition of being smooth is local and hence we may replace $U'$ by $\bigcup_{g\in N}(g.U)$. 
This allows us to assume that $U'$ is $N$-equivariant. 
We assume further that $U=U'$, i.e., that $\phi$ is smooth on $U$.  

Next, denote by $I$ the defining ideal of $\nt^*\setminus U$. 
Thanks to Lie Theorem we have $I^N\ne\{0\}$.
Fix $f\in I^N\setminus\{0\}$ and denote by $\nt^*_f$ the principal open subset of $\nt^*$ defined by $f$, and by $Q_f$ the principal open subset of $Q$ defined by $f$
(recall that $Y(\nt)\subset\Cp[Q]^N$). 
Finally, note that $\Cp[U]$ is the localization $\Sa(\nt)[f^{-1}]$ of $\Sa(\nt)$ by~$f$ and $\Cp[Q]= (\Sa(\nt)[f^{-1}])^N=Y(\nt)[f^{-1}]$. 
From now on we identify $\Cp[Q]$ with $Y(\nt)[f^{-1}]\subset \Sa(\nt)[f^{-1}]$. 
Pick a maximal ideal $M$ of $\Cp[Q]$. 
Since $\phi$ is smooth, the ideal of $\Cp[U]$ generated by $M$ is prime and is the defining ideal of a single $N$-orbit $N.\lambda$. 
It is easy to argue that $\Sa(\nt)[f^{-1}]/M=\Sa(\nt)/(M\cap Y(\nt))$. 
This implies that $U$ is a required open subset of~$\nt^*$.
}
\exam{In the rest\label{exam:center_nilradical_fd} of the subsection we will consider as an example the case of maximal nilpotent subalgebras of simple finite-dimensional Lie algebras. Let $\gt$, $\htt$, $\bt$, $\nt$, $\Phi$ be as in Subsection~\ref{sst:nilradicals_fin_dim}. Denote by $\Bu$ the following subset of $\Phi^+$:
\begin{equation*}
\Bu=\begin{cases}\bigcup\nolimits_{1\leq i\leq[n/2]}\{\epsi_i-\epsi_{n-i+1}\}&\text{for }A_{n-1},\\
\bigcup\nolimits_{1\leq i\leq n/2}\{\epsi_{2i-1}-\epsi_{2i},~\epsi_{2i-1}+\epsi_{2i}\}&\text{for }B_n,~n\text{ even},\\
\bigcup\nolimits_{1\leq i\leq[n/2]}\{\epsi_{2i-1}-\epsi_{2i},~\epsi_{2i-1}+\epsi_{2i}\}\cup\{\epsi_n\}&\text{for }B_n,~n\text{ odd},\\
\bigcup\nolimits_{1\leq i\leq n}\{2\epsi_i\}&\text{for }C_n,\\
\bigcup\nolimits_{1\leq i\leq[n/2]}\{\epsi_{2i-1}-\epsi_{2i},~\epsi_{2i-1}+\epsi_{2i}\}&\text{for }D_n.\\
\end{cases}
\end{equation*}
Note that $\Bu$ is a maximal strongly orthogonal subset of $\Phi^+$, i.e., $\Bu$ is maximal with the property that if $\alpha,\beta\in\Bu$ then neither $\alpha-\beta$ nor $\alpha+\beta$ belongs to $\Phi^+$. The set $\Bu$ is called the \emph{Kostant cascade} of orthogonal roots in~$\Phi^+$.

The description of generators of $Z(\nt)$ (or, equivalently, of $Y(\nt)$) goes back to J. Dixmier, A. Joseph and B. Kostant \cite{Dixmier1}, \cite{Joseph1}, \cite{Kostant1}, \cite{Kostant2}. We can consider $\Zp\Phi$, the $\Zp$-linear span of $\Phi$, as a subgroup of the group $\Xt$ of rational multiplicative characters of $H$ by putting $\pm\epsi_i(h)=h_{i,i}^{\pm1}$, where $H=\exp(\htt)$ and $h_{i,i}$ is the $i$-th diagonal element of a matrix $h\in H$. Recall that a vector $\lambda\in\Rp^n$ is called a \emph{weight} of $H$ if $2(\alpha,\lambda)/(\alpha,\alpha)$ is an integer for any $\alpha\in\Phi^+$, where $(\cdot,\cdot)$ is the standard inner product on $\Rp^n$. A weight $\lambda$ is called \emph{dominant} if $2(\alpha,\lambda)/(\alpha,\alpha)\geq0$ for all $\alpha\in\Phi^+$. An element $a$ of an $H$-module is called an \emph{$H$-weight vector}, if there exists $\nu\in\Xt$ such that $h\cdot a=\nu(h)a$ for all $h\in H$. By \cite[Theorems 6, 7]{Kostant2}, every $H$-weight occurs in $Y(\nt)$ with multiplicity at most 1. Furthermore, there exist unique (up to scalars) prime polynomials $\xi_{\beta}\in Y(\nt)$, $\beta\in\Bu$, such that each $\xi_{\beta}$ is an $H$-weight polynomial of a dominant weight $\vpi_{\beta}$ belonging to the $\Zp$-linear span $\Zp\Bu$ of $\Bu$. A remarkable fact is that $\xi_{\beta},~\beta\in\Bu$, are algebraically independent generators of $Y(\nt)$, so $Y(\nt)$ and $Z(\nt)$ are polynomial rings. We call $\{\xi_{\beta},~\beta\in\Bu\}$ the set of \emph{canonical generators} of $Y(\nt)$. The explicit formulas for the weights $\vpi_{\beta}$'s can be found, e.g., in \cite[Theorem 2.12]{Panov1}.

Below we present explicit formulas for $\xi_{\beta}$ for classical root systems (see \cite[Subsection 2.1]{Ignatyev1} for the details). We will use these formulas in Subsection~\ref{sst:nilradicals_inf_dim} considering infinite-dimensional setting. If $\Phi=A_{n-1}$ then $\vpi_{\beta}=\epsi_1+\ldots+\epsi_i-\epsi_{n-i+1}-\ldots-\epsi_n$ for $\beta=\epsi_i-\epsi_{n-i+1}$, $1\leq i\leq[n/2]$, and
\begin{equation}
\xi_{\epsi_i-\epsi_{n-i+1}}=\begin{vmatrix}\label{formula:Delta_i_A_n}
e_{1,n-i+1}&\ldots&e_{1,n-1}&e_{1,n}\\
e_{2,n-i+1}&\ldots&e_{2,n-1}&e_{2,n}\\
\vdots&\iddots&\vdots&\vdots\\
e_{i,n-i+1}&\ldots&e_{i,n-1}&e_{i,n}\\
\end{vmatrix}=\begin{vmatrix}
e_{\epsi_1-\epsi_{n-i+1}}&\ldots&e_{\epsi_1-\epsi_{n-1}}&e_{\epsi_1-\epsi_n}\\
e_{\epsi_2-\epsi_{n-i+1}}&\ldots&e_{\epsi_2-\epsi_{n-1}}&e_{\epsi_2-\epsi_n}\\
\vdots&\iddots&\vdots&\vdots\\
e_{\epsi_i-\epsi_{n-i+1}}&\ldots&e_{\epsi_i-\epsi_{n-1}}&e_{\epsi_i-\epsi_n}\\
\end{vmatrix}
\end{equation}
(cf. Example~\ref{exam:regular_orbits_A_n}). For $\Phi=C_n$ and $\beta=2\epsi_i$, $1\leq i\leq n$, one has $\vpi_{\beta}=2\epsi_1+\ldots+2\epsi_i$ and
\begin{equation}
\xi_{\beta}=\begin{vmatrix}\label{formula:Delta_i_C_n}
e_{\epsi_1+\epsi_i}&\ldots&e_{\epsi_1+\epsi_3}&e_{\epsi_1+\epsi_2}&2e_{2\epsi_1}\\
e_{\epsi_2+\epsi_i}&\ldots&e_{\epsi_2+\epsi_3}&2e_{2\epsi_2}&e_{\epsi_1+\epsi_2}\\
e_{\epsi_3+\epsi_i}&\ldots&2e_{2\epsi_3}&e_{\epsi_2+\epsi_3}&e_{\epsi_1+\epsi_3}\\
\vdots&\iddots&\vdots&\vdots&\vdots\\
2e_{2\epsi_i}&\ldots&e_{\epsi_3+\epsi_i}&e_{\epsi_2+\epsi_i}&e_{\epsi_1+\epsi_i}\\
\end{vmatrix}.
\end{equation}
Finally, if $\Phi=B_n$ or $D_n$ and $\beta=\epsi_{2i-1}+\epsi_{2i}$, $1\leq i\leq[n/2]$, then $\vpi_{\beta}=\epsi_1+\ldots+\epsi_{2i}$ and
\begin{equation}
\xi_{\beta}^2=\pm\begin{vmatrix}
e_{\epsi_1+\epsi_{i+1}}&\ldots&e_{\epsi_1+\epsi_3}&e_{\epsi_1+\epsi_2}&0\\
e_{\epsi_2+\epsi_{i+1}}&\ldots&e_{\epsi_2+\epsi_3}&0&-e_{\epsi_1+\epsi_2}\\
e_{\epsi_3+\epsi_{i+1}}&\ldots&0&-e_{\epsi_2+\epsi_3}&-e_{\epsi_1+\epsi_3}\\
\vdots&\iddots&\vdots&\vdots&\vdots\\
0&\ldots&-e_{\epsi_3+\epsi_{i+1}}&-e_{\epsi_2+\epsi_{i+1}}&-e_{\epsi_1+\epsi_{i+1}}\\
\end{vmatrix}.\label{formula:Delta_i_B_n_D_n_pf}
\end{equation}
(After a suitable reordering of indices, the matrix in the right-hand side becomes skew-symmetric, so $\xi_{\beta}$ is nothing but its Pfaffian.) Our normalization is such that the term $e_{\epsi_1+\epsi_2}e_{\epsi_3+\epsi_4}\ldots e_{\epsi_i+\epsi_{i+1}}$ enters $\xi_{\beta}$ with coefficient 1. (We will not use the canonical generators $\xi_{\beta}$ for $\beta=\epsi_{2i-1}-\epsi_{2i}$ for these root systems.)}

Now, let $J$ be a primitive ideal of $\Ua(\nt)$ (for an arbitrary classical root system). Put ${\Delta_{\beta}=\sigma(\xi_{\beta})}$. Thanks to Lemma~\ref{Lschur} there are unique scalars $c_{\beta}\in\Cp$ such that $\Delta_{\beta}-c_{\beta}\in J$, $\beta\in\Bu$. Denote by $\Delta$ the set of simple roots in $\Phi^+$. The description of the centrally generated ideals given in \cite[Theorem 3.1]{IgnatyevPenkov1} and \cite[Theorem~2.4]{Ignatyev1} is as follows: $J$ is centrally generated if and only if $c_{\beta}\neq0$ for all $\beta\in\Bu\setminus\Delta$. In particular, an element $f$ from Proposition~\ref{prop:open_subset_Z_n_cent_gen} can be chosen to be of the form 
\begin{equation}
f=\prod_{\beta\in\Bu\setminus\Delta}\label{formula:cent_gen_fin_dim}\Delta_{\beta}.
\end{equation}

\sect{Locally nilpotent Lie algebras}\label{sect: 1mr}\addcontentsline{toc}{section}{\ref{sect: 1mr}. Locally nilpotent Lie algebras}

The key results of this section are Theorem~\ref{theo:orbit_method_ifd} and Proposition~\ref{prop:radical_U_n_S_n_ifd}, which establish a complete analogue of the orbit method in the infinite-dimensional setting. 

We would like to describe the content of this section in more details. 
Subsection~\ref{sst:pro_varieties} contains necessary definitions and properties of pro-varieties needed for the studying of the dual space of an infinite-dimensional Lie algebra. 
Of course, these properties are in some sense ``dual'' to the properties of ind-varieties. 
The only ind-varieties we use are countable-dimensional vector spaces, so we do not discuss the general theory of ind-varieties in this paper.

In Subsection~\ref{sst:locally_nilp_Lie_algs}, we define the class of infinite-dimensional Lie algebras we are interested in (namely, locally nilpotent Lie algebras). 
Then we establish an inclusion-preserving bijection between radical ideals of $\Ua(\nt)$ and radical Poisson ideals of $\Sa(\nt)$ for such an algebra $\nt$, see Proposition~\ref{prop:radical_U_n_S_n_ifd}. 
We also prove that this bijection sends prime ideals to prime Poisson ideals and vice versa. 
Next, given a linear form $\lambda\in\nt^*$, we construct the ideal $J(\lambda)$ of $\Ua(\nt)$ and check that this ideal is primitive, see Theorem~\ref{theo:primitive_ifdv}. 
After that we present an alternative description of the ideal $I(\lambda)$ and, using this description, prove  (Subsection~\ref{sst:prime_theorem},  Theorem~\ref{main_theo:orbit_method}) that each primitive ideal of $\Ua(\nt)$ has the form $J(\lambda)$ for some $\lambda\in\nt^*$. 
This shows that a bijection of Proposition~\ref{prop:radical_U_n_S_n_ifd} gives rise to a homeomorphism between the space of primitive ideals of $\Ua(\nt)$ and the space of primitive Poisson ideals of $\Sa(\nt)$; we consider this as the first main result, see Theorem~\ref{theo:orbit_method_ifd}.

As an example, we discuss in Subsections~\ref{sst:socle} and~\ref{sst: adjg} a special class of locally nilpotent Lie algebras, so-called socle Lie algebras, where the homeomorphism between $\JSpec{\Ua(\nt)}$ and $\PSpec{\Sa(\nt)}$ mentioned above can be obtained via the coadjoint action of a certain pro-group on the dual space $\nt^*$ (similarly to the finite-dimensional case).

\sst{Pro-varieties}\label{sst:pro_varieties}\addcontentsline{toc}{subsection}{\ref{sst:pro_varieties}. Pro-varieties}
In this subsection, we briefly recall basic facts about pro-varieties which are needed for the sequel. By definition, a \emph{pro-variety} is the projective limit $X=\plm X_n$ of a chain of morphisms of algebraic varieties
\begin{equation}
\label{formula:pro_variety}
X_1\stackrel{\vfi_1}\ot X_2\stackrel{\vfi_2}\ot\ldots\stackrel{\vfi_{n-1}}{\ot}
X_n\stackrel{\vfi_n}\ot X_{n+1}\stackrel{\vfi_{n+1}}\ot\ldots.
\end{equation} Obviously, the projective limit of a chain (\ref{formula:pro_variety}) does not change if we replace the sequence $\{X_n\}_{n\geq1}$ by a subsequence $\{X_{i_n}\}_{n\geq 1}$, and the morphisms $\vfi_n$ by the compositions $$\wt\vfi_{i_n}=\vfi_{i_n}\circ\vfi_{i_n+1}\circ\ldots\circ\vfi_{i_{n+1}-1}.$$
In what follows we only consider chains (\ref{formula:pro_variety}) where the morphisms $\vfi_n$ are dominant.



The following lemma is very natural and seems to be well known.
\lemmp{ \label{Lprone}If $X_n\ne\varnothing$ for all $n\ge1$ then $X$ has at least one point.}{Let $X=\varprojlim X_n$ be a pro-variety. Then every $X_n$ can be covered by a finite collection of affine charts. Every such a chart is a union of several irreducible affine varieties. Thus without loss of generality we can assume that $X_n$ is affine and irreducible for all $n$.

Set $R_n:=\mathbb C[X_n]$ to be the algebra of regular functions on the affine variety $X_n$. This defines a sequence $R_1\to R_2\to\ldots$ of monomorphisms which is dual to the sequence $X_1\leftarrow X_2\leftarrow\ldots$. Put $R:=\lim\limits_{n\to\infty}R_n$. Pick a maximal ideal $M$ of $R$. It is clear that $R$ is at most countable-dimensional and hence by Corollary~\ref{Lschurc} we have $R/M\cong\mathbb C$. 
A maximal ideal of $R$ with this property defines a point of $X$ in a straightforward way.}

\exam{ Let $V$ be a countable-dimensional complex vector space with a fixed basis $E=\{e_1,~e_2,~\ldots\}$. Set $V_n=\langle e_1,~\ldots,~e_n\rangle_{\Cp}$, then $V=\ilm V_n$. Let $V^*$ and $V_n^*$ be the dual spaces of $V$ and $V_n$, $n\geq1$, respectively. Then $V^*$ is the projective limit $V^*=\plm V_n^*$ where morphisms $V_n^*\ot V_{n+1}^*$ are nothing but the restrictions of linear functions from $V_{n+1}^*$ to $V_n$. Thus, $V^*$ is a pro-variety.}

In this special case we will consider two different topologies on $V^*$. First, let $\Sa(V)$ be the symmetric algebra of the space $V$. The points of $\JSpec{\Sa(V)}=\MSpec{\Sa(V)}$ of $\Sa(V)$ can be identified with $V^*$. This introduces the first topology on $V^*$; we call the topology induced by this construction the {\it $\Sa(V)$-Zariski topology} on $V^*$.

Second, we say that a subset $Z$ of $V^*$ is \emph{closed in countable-Zariski topology} if $Z$ is a union of countably many $\Sa(V)$-Zariski closed subsets of $V$. One can immediately see that it is again a topology on $V^*$, which is finer than the $\Sa(V)$-Zariski topology. Some properties of this topology are ``strange'' even in the finite-dimensional case: for example, $\mathbb{Q}$ is a closed subset of $\Cp$ in this topology. On the other hand, we have the following property needed in Subsection~\ref{sst:cent_gen_ideals_ifd}.

\propp{\label{prop:Zariski_countable_irr_v2} Let $V$ be a countable-dimensional $\mathbb C$-vector space. Then $V^*$ is irreducible with respect to the countable-Zariski topology.}{Assume to the contrary that $V^*$ can be represented as a union of two proper nonempty countable-Zariski closed subsets. This implies that there exist nonempty proper $\Sa(V)$-Zariski closed subsets $Z_n$ of $V^*$, $n\geq1$, such that $V^*=\bigcup_nZ_n$. We may assume without loss of generality that each $Z_n$ has the form $$Z_n=\{\lambda\in V^*\mid f_n(\lambda)=0\}$$ for certain $f_n\in\Sa(V)$.

Consider the localisation $\Sa(V)[f_n^{-1},~n\geq1]$ of $\Sa(V)$ by $f_n,~n\geq1$. Pick a maximal ideal $M$ of $\Sa(V)[f_n^{-1},~n\geq1]$. It is evident that $\Sa(V)[f_n^{-1},~n\geq1]$ is countable-dimensional and hence $M$ is of codimension 1 in $\Sa(V)$. Let $\lambda$ be the point of $V^*$ corresponding to $M$. By definition we have that $f_n(\lambda)\neq0$ for all $n\geq1$. Therefore, $\lambda\notin\bigcup_n Z_n$ and, consequently, $V^*\neq\bigcup_n Z_n.$}
Note that the proof of Proposition~\ref{prop:Zariski_countable_irr_v2} is very similar to the proof of Lemma~\ref{Lprone}
\sst{Locally nilpotent Lie algebras}\label{sst:locally_nilp_Lie_algs}\addcontentsline{toc}{subsection}{\ref{sst:locally_nilp_Lie_algs}. Locally nilpotent Lie algebras}

Now we will introduce the main definition of the paper.

\defi{Let $\nt$ be a countable-dimensional Lie algebra expressed as an inductive limit of its nested finite-dimensional nilpotent subalgebras
\begin{equation*}
\nt_1\subset\nt_2\subset\ldots\subset\nt_k\subset\ldots.\label{formula:locally_nilpotent}
\end{equation*}
Then $\nt$ is called \emph{locally nilpotent}.}
From now on, assume that $\nt$ is locally nilpotent. 
The definitions of $\Ua(\nt)$ and $\Sa(\nt)$ coincide with the definitions of these algebras in the finite-dimensional case. 
It is clear that $\Ua(\nt)$ and $\Sa(\nt)$ are also countably-dimensional, and $\Sa(\nt)$ is a Poisson algebra. Below we will prove that, for the Lie algebra $\nt$, we still have a bijection between radical ideals of $\Ua(\nt)$ and radical Poisson ideals of $\Sa(\nt)$, see Proposition~\ref{prop:radical_U_n_S_n_ifd}. Furthermore, it defines a homeomorphism between $\JSpec{\Ua(\nt)}$ and $\PSpec{\Sa(\nt)}$ thanks to Theorem~\ref{theo:orbit_method_ifd}.

We need to set up the notation. 
Pick a Poisson ideal $I$ of $\Sa(\nt)$ and an ideal $J$ of ${\rm U}(\nt)$.
Put $$J_k=J\cap \Ua(\nt_k),~I_k=I\cap \Sa(\nt_k)$$ for $k\geq1$.
Clearly, each $I_k$ is a Poisson ideal of $\Sa(\nt_k)$ and each $J_k$ is a two-sided ideal of ${\rm U}(\nt_k)$.
Moreover, $I$ is the inductive limit $I=\ilm I_k$, and $J$ is the inductive limit of the respective sequence of the ideals~$J_k$, $k\geq1$.
\lemmp{The ideal \label{lemm:radical_intersections_finitary} $I$ \textup(respectively\textup, $J$\textup) is radical if and only if each $I_k$ \textup(respectively\textup, each $J_k$\textup) is radical.
}{For $I$, this is a simple exercise in commutative algebra, so we will proceed for $J$. If all $J_k$'s are radical then $J$ is radical by~(\ref{formula:radical_ideal_elements}).
Assume now that $J$ is radical.
Formula (\ref{formula:radical_ideal_elements}) implies that
$$J\cap\Ua(\nt_k)=\sqrt J\cap\Ua(\nt_k)=\bigcap_{l\ge k}\sqrt{J\cap \Ua(\nt_l)}\cap\Ua(\nt_k).$$
It follows from Remark~\ref{nota:ideals_in_Un} that $\sqrt{J\cap \Ua(\nt_l)}$ is an intersection of several completely prime ideals of~${\rm U}(\nt_l)$.
Thus $\sqrt{J\cap \Ua(\nt_l)}\cap\Ua(\nt_k)$ is an intersection of several completely prime ideals of ${\rm U}(\nt_k)$.
Hence $\sqrt{(J\cap \Ua(\nt_l)}\cap\Ua(\nt_k)$ is an intersection of several prime ideals.
Therefore $J\cap\Ua(\nt_k)$ is a radical ideal. The result follows.}

\propp{\textup{i)}\label{prop:radical_U_n_S_n_ifd} Let $I$ be a radical Poisson ideal of ${\rm S}(\nt)$. Set $J(I)=\bigcup_kJ(I_k)$. Then $J(I)$ is a radical ideal of $\Ua(\nt)$. \textup{ii)} Let $J$ be a radical ideal of ${\rm U}(\nt)$. Set $I(J)=\bigcup_kI(J_k)$. Then $I(J)$ is a radical Poisson ideal of $\Sa(\nt)$. \textup{iii)} The maps $I\mapsto J(I),~J\mapsto I(J)$ provide an inclusion-preserving bijection between the radical Poisson ideals of $\Sa(\nt)$ and the radical ideals of $\Ua(\nt)$.}
{All the parts are implied by a combination of Lemma~\ref{lemm:radical_intersections_finitary} and Theorem~\ref{theo:Um_Un}.}
\corop{\label{Crpcp} A radical Poisson ideal $I\subset\Sa(\nt)$ is prime if and only if $J(I)$ is prime. Moreover\textup, if $J(I)$ is prime then $J(I)$ is completely prime.}
{ Fix a radical Poisson ideal $I\subset\Sa(\nt)$. 
Assume $I$ is prime. 
Then $I_k=I\cap\Sa(\nt_k)$ is a prime ideal of $\Sa(\nt_k)$ for all $k\ge1$, and hence $J(I_k)$ is prime for all $k\ge1$, see~\cite[Proposition 6.3.5]{Dixmier1}. Hence $J(I_k)$ is completely prime for all $k\ge1$, see Remark~\ref{nota:ideals_in_Un}.
This implies that $J(I)$ is prime and completely prime.

Assume that $I$ is not prime. Then there exist $a, b\in\Sa(\nt)$ with $ab\in I$ and $a, b\notin I.$ 
To proceed we use notation of~\cite{PetukhovSierra1}. 
Set
$$I_b:=(I:a):=\{f\in\Sa(\nt)\mid af\in I\},\hspace{10pt}I_a:=(I\hat+a):=\bigcap_{c\in(I:a)}(I:c).$$ Thanks to~\cite[Lemma~2.1]{PetukhovSierra1} both $I_a$ and $I_b$ are radical Poisson ideals. Moreover,~\cite[Lemma~2.1]{PetukhovSierra1} implies that $b\in I_b$ (therefore $I_b\not\subset I$), $a\in I_a$ (therefore $I_a\not\subset I$), and by~\cite[Lemma 2.3]{PetukhovSierra1} we have $I_a\cap I_b\subset \sqrt I=I.$ 
From Proposition~\ref{prop:radical_U_n_S_n_ifd} we have
\begin{equation}\label{Eidjab}J(I_a)\not\subset J(I),~J(I_b)\not\subset J(I),~J(I_a)\cap J(I_b)=J(I_a\cap I_b)\subset J(I).\end{equation}
Thanks to Lemma~\ref{Lintrad} we have $\sqrt{J(I_a)\cap J(I_b)}=\sqrt{J(I_a)J(I_b)}$ and hence $J(I_a)J(I_b)\subset J(I)$. This together with~(\ref{Eidjab}) implies that $J(I)$ is not prime.
}
\defi{ Pick $\lambda\in\nt^*$. Assign to $\lambda$ the primitive Poisson ideal $I(\lambda)$ of~$\Sa(\nt)$, which is by definition the largest Poisson ideal in $I_{\lambda}$, where $I_{\lambda}$ is the kernel of the evaluation map $$\Sa(\nt)\to\Cp\colon f\mapsto f(\lambda).$$ 
Corollary~\ref{Lschurc} implies that $I(\lambda),~\lambda\in\nt^*$, are all the primitive Poisson ideals of ${\rm S}(\nt)$. Set $J(\lambda):=J(I(\lambda))$.
}
\nota{ For a Noetherian associative algebra it is known that a prime ideal is radical. We do not know whether or not this holds in our setting.}

It is natural to expect that $J(\lambda),~\lambda\in\nt^*$, are all the primitive ideals of $\Ua(\nt)$.
We first show that each $J(\lambda)$, $\lambda\in\nt^*$, is primitive.
The proof is based on Theorem~\ref{Thm:bim}.
\theop{\label{theo:primitive_ifdv} For every $\lambda\in\nt^*$ the ideal $J(\lambda)$ is primitive.}
{ Without loss of generality assume that $\dim\nt_1=1$. 
The restriction $\restr{\lambda}{\nt_1}$ is a character of~$\nt_1$ and hence it defines a one-dimensional $\nt_1$-module $M_1$. 
Set $M_{i+1}={}_{\nt_{i+1}}F_{\nt_i}\otimes_{\Ua(\nt_i)}M_i$ for all $i\ge1$ where ${}_{\nt_{i+1}}F_{\nt_i}$ is a $\Ua(\nt_{i+1}){-}\Ua(\nt_i)$ bimodule defined by Theorem~\ref{Thm:bim}. 
Theorem~\ref{Thm:bim} also implies that $M_i$ is a simple $\nt_i$-module for all $i\ge1$. 
Then the map $\phi$ from Theorem~\ref{Thm:bim} defines the $\nt_i$-embedding $M_i\to M_{i+1}$ for all $i\ge1$. 
Let $M=\ilm M_i$. 
It is clear that $M$ is a limit of simple $\nt_i$-modules with ${\rm Ann}_{ \Ua(\nt)}M = J(\lambda)$. This immediately implies that $M$ is a simple $\nt$-module and $J(\lambda)$ is primitive. }

\sst{An alternative description of $I(\lambda)$}\label{sst:alternative_J_lambda}\addcontentsline{toc}{subsection}{\ref{sst:alternative_J_lambda}. An alternative description of $I(\lambda)$}
In this subsection we provide an alternative characterisation of ideals of the form $I(\lambda)$. Such a characterisation will be used to show that every primitive ideal of $\Ua(\nt)$ is of the form $J(\lambda)$ for a locally nilpotent Lie algebra $\nt$, see Theorem~\ref{main_theo:orbit_method}. 

Let $\nt$ be a locally nilpotent Lie algebra together with an exhaustion $\nt_1\subset\nt_2\subset\ldots$ of $\nt$ by its finite-dimensional nilpotent subalgebras.
Pick a radical Poisson ideal $I\subset\Sa(\nt)$ of $\Sa(\nt)$, and let $J:=J(I)$ be the corresponding radical two-sided ideal of $\Ua(\nt)$. 
As above, set $J_n:=J\cap\Ua(\nt_n), I_n:=I\cap\Sa(\nt_n)$.
For every ideal $I'\subset\Sa(\nt_n)$ we let ${\rm Var}(I')\subset\nt_n^*$ to be the set of common zeros of $I'$ in $\nt_n^*$, and put $V_n:={\rm Var}(I_n)$.
Recall that $N_l=\Exp{\nt_l}$ denotes the unipotent group attached to $\nt_l$. 
Set $\phi_{l\to n}$ to be the canonical $N_n$-equivariant  projection $\nt_l^*\to\nt_n^*$. 
The description of $V_n$ for $I=I(\lambda)$ is given in the following lemma. 
\lemmp{\label{Lvaril} Pick $\lambda\in\nt^*$ and assume $I=I(\lambda)$. Further\textup, set $O(\lambda; l)=N_l.\restr{\lambda}{\nt_l}$ to be the coadjoint $N_l$-orbit of the linear form $\restr{\lambda}{\nt_l}\in\nt_l^*$.
Then\\
\indent$(1)$ $\phi_{l\to n}O(\lambda; l)\subset\phi_{l+1\to n}O(\lambda; l+1)$ for all $l\ge n$\textup;\\
\indent$(2)$ $\overline{\phi_{l\to n}O(\lambda; l)}=\overline{\phi_{l+1\to n}O(\lambda; l+1)}$ for all $l\gg n$\textup;\\
\indent$(3)$ $V_n=\bigcup_{l\ge n}\overline{\phi_{l\to n}O(\lambda; l)}$ and $V_n=\overline{\phi_{l\to n}O(\lambda; l)}$ for all $l\gg n$.
}{By definition, $\phi_{l+1\to l}(\restr{\lambda}{\nt_{l+1}})=\restr{\lambda}{\nt_l}$ and hence $O(\lambda; l)\subset \phi_{l+1\to l}O(\lambda; l+1).$ 
Therefore $\phi_{l\to n}O(\lambda; l)\subset\phi_{l+1\to n}O(\lambda; l)$ for all $l\ge n$. This proves (1). 

Each variety $\overline{\phi_{l\to n}O(\lambda; l)}$ is irreducible because it is the closure of an image of the irreducible (i.e., connected) group $N_l$. Therefore the sequence of varieties $\overline{\phi_{l\to n}O(\lambda; l)}$ stabilizes for $l\gg n$ because all these varieties are contained in $\nt_n^*$. This proves (2). 

We left to show (3). Thanks to (1) and (2) it is enough to show that $V_n=\overline{\phi_{l\to n}O(\lambda; l)}$ for all $l\gg n$. 
It is equivalent to the condition $I_n=I(\restr{\lambda}{\nt_l})\cap\Sa(\nt_n)$ for all $l\gg n$. 
(Recall that $I(\restr{\lambda}{\nt_l})$ is the annihilator of $O(\lambda; l)$ in $\Sa(\nt_l)$.) 
It is enough to show that $I_n\subset I(\restr{\lambda}{\nt_l})\cap\Sa(\nt_n)$ and $I(\restr{\lambda}{\nt_l})\cap\Sa(\nt_n)\subset I_n$ for all $l\gg n$. 
The first inclusion is trivial because $I(\restr{\lambda}{\nt_l})$ is the largest Poisson ideal of $\Sa(\nt_l)$ contained in the maximal ideal attached to $\restr{\lambda}{\nt_l}$. 
For the second inclusion set $\tilde I_n:=\bigcap_{l\ge n}I(\restr{\lambda}{\nt_l})$ and note that thanks to step (2) we have $\tilde I_n=I(\restr{\lambda}{\nt_l})\cap\Sa(\nt_n)$ for all $l\gg n$. 
This implies that there exists a Poisson ideal $\tilde I$ of $\Sa(\nt)$ such that  $\tilde I\cap\Sa(\nt_n)=\tilde I_n$ for all $n\ge1$; then we have $\tilde I\subset I=I(\lambda)$ and thus $\tilde I_n\subset I_n$ for all $n\ge 1$.
}
The map $\phi_{l\to n}$ is (by definition) a moment map with respect to the action of $N_n$ on $V_l$;
a well-known feature of a moment map implies the following~\cite[the~last~paragraph~of~\S 2.5]{Vinberg}.  

\lemm{\label{Lmmap} Let $O_l$ be a $N_l$-coadjoint orbit of $\nt_l^*$. Set $d_n(O_l)$ to be the maximal dimension of an $N_n$-orbit on $O_l$. Then $\dim\phi_{l\to n}(O_l)=d_n(O_l)$.} 
Let $I$ again be an arbitrary Poisson ideal of $\Sa(\nt)$.
Set $d_{l\to n}=d_n(V_l)$ to be the maximal dimension of an $N_n$-orbit on $V_l$. 
Fix $n\ge1$ and consider the sequence $d_{n\to n},~d_{n+1\to n},~d_{n+2\to n},~\ldots$. It is clear that this sequence is nondecreasing and that it is bounded by $\dim N_n=\dim \nt_n$. 
This implies that this sequence stabilizes from some point and we denote the stable value of this sequence by $d_n(I)$. 

Now we have enough tools to provide an alternative characterization of ideals of the form $I(\lambda)$. 
\propp{\label{Pppoicr} Assume $I$ is prime. Then $I=I(\lambda)$ for some $\lambda$ if and only if $d_n(I)=\dim V_n$ for all $n\ge 1$.}
{It is clear from Lemmas~\ref{Lvaril},~\ref{Lmmap} that if $I=I(\lambda)$ for some $\lambda\in\nt^*$ then $d_n(I)=\dim V_n$ for all $n\ge 1$. Thus we left to check the opposite statement. 

From now on we assume that $\dim V_n=d_n(I)$ for all $n\ge 1$. 
This implies that for every $n\ge 1$ there exists $l\ge n$ such that $d_n(V_l)=\dim V_n$. 
In more details, this means that there exists an $N_l$-coadjoint orbit $O_l\subset\nt^*_l$ such that $d_n(O_l)=\dim V_n$. 

General arguments imply that there exists a nonempty open $N_n$-stable subset $V_l^\circ$ of $V_l$ such that the dimensions of all $N_n$-orbits from $V_l^\circ$ equal $d_n(V_l)$. 
The complement to $V_l^\circ$ in $V_l$ can be described as the zero set of a finite collection of polynomials. 
The union of all such finite collections of polynomials for all suitable pairs $l, n$ is at most countable and we denote them $f_1, f_2,\ldots$. 

Consider the localization $(\Sa(\nt)/I)[f_1^{-1}, f_2^{-1},\ldots]$ of $\Sa(\nt)/I$ and a maximal ideal $$M\subset (\Sa(\nt)/I)[f_1^{-1}, f_2^{-1},\ldots]$$ inside it. 
Corollary~\ref{Lschurc} implies that $M$ has codimension 1 in $(\Sa(\nt)/I)[f_1^{-1}, f_2^{-1},\ldots]$. 
Hence $M$ defines a homomorphism of rings $$(\Sa(\nt)/I)[f_1^{-1}, f_2^{-1},\ldots]\to\mathbb C.$$
Denote by $\lambda\in\nt^*$ the linear form defined by $M$. By definition we have $I\subset I(\lambda).$ 

We left to show that $I(\lambda)\subset I$, i.e., that $I=I(\lambda)$. 
The definition of $\lambda$ implies that 
\begin{equation}\label{Eqdv}d_n(I(\lambda))=\dim V_n.\end{equation} 
On the other hand both $V_n$ and ${\rm Var}(I(\lambda)_n)$ are irreducible and ${\rm Var}(I(\lambda)_n)\subset V_n$. Together with (\ref{Eqdv}) this implies the desired result.}

\okr{Corollary}{ \label{Cprtpro} Let $I$ be a prime Poisson ideal of $\Sa(\nt)$. Assume that\textup, for every $n\ge1$\textup, there exist $l\ge n$ and a coadjoint $N_l$-orbit $O_l\subset V_l\subset\nt_l^*$ such that the canonical map $O_l\to V_n$ is dominant. Then $I=I(\lambda)$ for some $\lambda\in\nt^*$.}
\proof{}{ Thanks to Lemmas~\ref{Lvaril} and~\ref{Lmmap} we have $d_n(V_l)\ge d_n(O_l)=\dim V_n$. 
Hence $$d_n(I)=\dim V_n.$$ 
This together with Proposition~\ref{Pppoicr} implies the desired result.}
\sst{The orbit method for locally nilpotent Lie algebras}\label{sst:prime_theorem}\addcontentsline{toc}{subsection}{\ref{sst:prime_theorem}. The orbit method for locally nilpotent Lie algebras}
In this subsection we will freely use notions, definitions and conventions related to skew fields of associative (noncommutative) algebras. A very basic introduction to this subject is given in~\cite{Dixmier1}, and there are quite a lot of books in the area of skew fields with no focus on Lie algebras, see, f.e.,~\cite{Cn1}, \cite{Cn2}, \cite{Cn3}.

Let $\nt=\ilm \nt_n$ be a locally nilpotent Lie algebra and let $J$ be a radical two-sided ideal of $\Ua(\nt)$. 
Set $J_n:=J\cap\Ua(\nt_n)$. 
Whenever $J$ is prime it is completely prime (Corollary~\ref{Crpcp}) and we set $Q(\nt; J)$ to be the limit of the quotient skew fields $Q(\nt_n; J_n)$ of $\Ua(\nt_n)/J_n$.

In this subsection we prove that there is a bijection between the primitive ideals of $\Ua(\nt)$ and the primitive ideals of $\Sa(\nt)$ which extends to the homeomorphism between $\JSpec{\Ua(\nt)}$ and $\PSpec{\Sa(\nt)}$, cf. Corollary~\ref{Crpcp}. 
In fact, the only thing which remains to be checked is the following theorem.
\theop{\label{main_theo:orbit_method}The following conditions are equivalent\textup:
\begin{equation*}
\begin{split}
&\text{\textup{i)} $J$ equals $J(\lambda)$ for some $\lambda\in\gt^*$};\\
&\text{\textup{ii)} $J$ is primitive};\\
&\text{\textup{iii)} $J$ is prime and radical and $Q(\nt;J)$ has trivial center}.\\
\end{split}
\end{equation*}
}{ Thanks to Theorem~\ref{theo:primitive_ifdv}, (i) implies (ii). It is also clear that (ii) implies that $J$ is prime and radical. The last condition of (iii) is implied by a minor modification of~\cite[Lemma~4.1.6,~Proposition~4.1.7]{Dixmier1}.
(The only fact about $\Ua(\nt)$ and $J$ needed for \cite[Lemma~4.1.6]{Dixmier1} is the existence of~$Q(\nt; J)$.) Therefore we left to show that (iii) implies (i). The proof of this implication is given by Proposition~\ref{Pprt0} below.}

As an immediate corollary we obtain the following result (which we consider as the first main results of the paper).

\mtheo{Let $\nt$ be\label{theo:orbit_method_ifd} a countable-dimensional locally nilpotent complex Lie algebra. Then
\begin{equation*}
\begin{split}
&\text{\textup{i)} each primitive ideal of $\Ua(\nt)$ equals $J(\lambda)$ for a certain $\lambda\in\nt^*$};\\
&\text{\textup{ii)} each primitive Poisson ideal of $\Sa(\nt)$ equals $I(\lambda)$ for a certain $\lambda\in\nt^*$};\\
&\text{\textup{iii)} the map $I(\lambda)\mapsto J(\lambda)$ is a homeomorphism between $\PSpec(\Sa(\nt))$ and $\JSpec(\Ua(\nt))$}.\\
\end{split}
\end{equation*}}
\okr{Proposition}{\label{Pprt0}Assume $J$ is prime and radical in $\Ua(\nt)$. If $Q(\nt; J)$ has trivial center then $J=J(\lambda)$ for some $\lambda\in\nt^*$.}
From now on we assume that $J$ is prime and radical ideal of $\Ua(\nt)$ and that $Q(\nt; J)$ has trivial center. 
To prove Proposition~\ref{Pprt0} we need more notation. 
Set $R_n:=\Ua(\nt_n)/J_n$. 
Corollary~\ref{Crpcp} implies that $J_n$ is a prime ideal of $\Ua(\nt_n)$ and hence $R_n$ is a (noncommutative) domain for all $n$. 
Denote by~$Z_n$ the center of $R_n$ and by $Q_n$ the quotient field of $Z_n$. 
It  is clear that
$$\Ua(\nt)/J\cong \ilm R_n$$
and that the center of $R:=\Ua(\nt)/J$ equals $\bigcap_nZ_n$. 
Recall that $I(J_n)$ denotes the prime Poisson ideal of $\Sa(\nt_n)$ corresponding to $J_n$, see Notation~\ref{Nij}. Further, $V_n$ stands for the set of common zeros of $I(J_n)$ in $\nt_n^*$ and ${\rm Var}(I')\subset\nt_n^*$ is the set of common zeros of an ideal $I'\subset \Sa(\nt_n^*)$. 

Recall the notion of linear disjointness, see, e.g., \cite[Chapter 5, \S 2.5]{Bou}.
It is easy to deduce Proposition~\ref{Pprt0} from Corollary~\ref{Cprtpro} and the following propositions.
\okr{Proposition}{\label{Pprtalg} Assume $J$ is prime and radical in $\Ua(\nt)$. If $Q(\nt; J)$ has trivial center then\textup, for every $n\ge 1$, there exists $l> n$ such that $Z_n$ and $Z_l$ are linearly disjoint. }
\okr{Proposition}{\label{Pprtinj} Consider $n\ge1$ and $l\ge n$. If $Z_n$ and $Z_l$ are linearly disjoint then there exists a primitive ideal $P$ of $R_l$ such that the canonical map $R_n\to R_l/P$ is injective.}
\okr{Proposition}{\label{Pprtdom} Consider $n\ge1$ and $l\ge n$ together with a primitive ideal $P$ of $R_l$. If the canonical map $R_n\to R_l/P$ is injective then $V_l$ contains a coadjoint $N_l$-orbit $O_l$ such that the canonical map $O_l\to V_n$ is dominant.}
Proposition~\ref{Pprtdom} is relatively simple and quite straightforward and we give the proof of it first. Propositions~\ref{Pprtalg} and~\ref{Pprtinj} are more involved and the rest of this section is devoted to their proofs. 
\proof{ of Proposition~\ref{Pprtdom}}{ Let $\wt P$ be the full preimage of $P$ under the canonical map $$\Ua(\nt_l)\to R_l.$$ 
The assumptions of Proposition~\ref{Pprtdom} imply that $J_n=\wt P\cap\Ua(\nt_n).$ 
This together with Theorem~\ref{theo:Um_Un} implies the desired condition for $O_l={\rm Var}(I(P_l))$. 
}

Now we proceed to the proofs of Propositions~\ref{Pprtalg} and~\ref{Pprtinj}.
\lemmp{We have $Z_n\cap Z_l=\Cp$ for all $l\gg n$.}{The algebras $R_n, n\ge1,$ are domains and therefore $Z_n, n\ge 1,$ are commutative domains. 
For a commutative domain $A$ we use notions of the transcendence degree $A\mapsto {\rm trdeg}(A)$ and the Gelfand--Kirillov dimension $A\mapsto{\rm GKdim}(A)$, see~\cite{KrauseLenagan1} for more details. We have
$${\rm GKdim}(Z_n\cap Z_i)={\rm trdeg}(Z_n\cap Z_i)\ge {\rm trdeg} (Z_n\cap Z_j)={\rm GKdim}(Z_n\cap Z_j)$$
whenever $i<j$. 
The transcendence degree is always a nonnegative integer and hence ${\rm trdeg}(Z_n\cap Z_l)$ is equal to the same $d\in\mathbb Z_{\ge0}$ for all $l\gg n$. 

The desired statement is equivalent to the condition $d=0$. This can be deduced from the following lemma.
\okr{Lemma}{\label{LWnormal} Let $\nt_0$ be a finite-dimensional nilpotent Lie algebra and $J_0$ be a prime ideal of~$\Ua(\nt_0)$. Denote by $Z(\nt_0; J_0)$ the center of $\Ua(\nt_0)/J_0$ and by $QZ(\nt_0; J_0)$ the quotient field of $Z(\nt_0; J_0)$. 
Then the following conditions are equivalent\textup:
\begin{equation*}
\begin{split}
&\text{\textup{i)} $f\in Z(\nt_0; J_0)$};\\
&\text{\textup{ii)} $f\in\Ua(\nt_0)/J_0$ is algebraic over $QZ(\nt_0; J_0)$};\\
\end{split}
\end{equation*}}
\proof{}{ According to~\cite[Proposition~4.7.1, Lemma~6.4.5]{Dixmier1} we have that $$(\Ua(\nt_0)/J_0)\otimes_{Z(\nt_0; J_0)}QZ(\nt_0; J_0)$$ is isomorphic to the Weyl algebra over $QZ(\nt_0; J_0)$. Thus it is enough to check a similar statement for the Weyl algebra. See~\cite[p. 295]{Goodearl2} for the latter fact.}
Pick $i$ with ${\rm trdeg}(Z_i\cap Z_n)=d$ and $i\ge n$. 
Then ${\rm trdeg}(Z_i\cap Z_n)={\rm trdeg}(Z_{j}\cap Z_n)$ for all $j>i$. 
Fix $j>i$. 
The above condition implies that every element of $Z_i\cap Z_n$ is algebraic over $Z_j$. 
Hence Lemma~\ref{LWnormal} implies that $Z_i\cap Z_n\subset Z_j$ for all $j\ge i$. Thus all elements of $Z_i\cap Z_n$ belong to the center of $\Ua(\nt)/J$ and therefore $Z_n\cap Z_i=\Cp$.
}

\proof{ of Proposition~\ref{Pprtinj}}{ 
Pick a $\Cp$-basis $z_1, z_2, \ldots$ of $Z_n$ with $z_1=1$ and extend it to a $\Cp$-basis $c_1, c_2, \ldots$ of $R_l$. 
Thanks to the assumption that $Z_n$ and $Z_l$ are linearly disjoint, we have that $z_1, z_2, \ldots$ are $Q_l$-linearly independent. 
Hence we can choose a subset $\{c_1', c_2', \ldots\}$ of the set $\{c_1, c_2, \ldots\}$ such that all $z_i$ are among $c_i'$  and $\{c_1', c_2', \ldots\}$ is a $Q_l$-basis of $R_l\otimes_{Z_l}Q_l$. 
For every $i, j\ge 1$ we have 
\begin{equation}\label{Emult}c_i'c_j'=\sum_{k}\alpha_{ijk}c_k', c_i=\sum_k \beta_{ik} c_k'\end{equation}
for the unique constants $\alpha_{ijk}, \beta_{ik}\in Q_l$ (here the right-hand sums are always finite). 
These constants are fractions of elements of $Z_l$ and we denote the corresponding denominators by $\gamma_1, \gamma_2, \ldots$. 
Pick $f$ as in Proposition~\ref{prop:open_subset_Z_n_cent_gen}.  
We add $f$ to the sequence of constant $\gamma_1, \gamma_2, \ldots$. 
It is easy to verify that $\gamma_1, \gamma_2, \ldots$ considered as a set is at most countable. 
Pick a maximal ideal $M$ of $Z_l$ such that $$\gamma_1, \gamma_2, \ldots\not\in M,\hspace{10pt}Z_l/M\cong\Cp$$ (such an ideal exists thanks to Corollary~\ref{Lschurc}). 
We wish to show that $P:=MR_l$ is a maximal two-sided ideal of $R_l$ (and hence primitive).

Denote by $\psi$ the homomorphism $Z_l\to Z_l/M\cong\Cp$. 
Thanks to the definition of $\psi$, we can extend $\psi$ to $\alpha_{ijk}$ and $\beta_{ik}$. 
Set $\bar\alpha_{ijk}:=\psi(\alpha_{ijk}), \bar\beta_{ik}:=\psi(\beta_{ik})$.
Then we define algebra $\bar R_l$ as a $\Cp$-vector space with basis $\bar c_1', \bar c'_2, \ldots$ and the multiplication law completely analogoues to~(\ref{Emult}). 
This algebra is associative because the multiplication law~(\ref{Emult}) is associative. 
We also can define a linear map $\wt\psi\colon R_l\to \bar R_l$ by formulas $$c_i\mapsto\sum_k\bar\beta_{ik}\bar c'_k,$$
and it is easy to check that this linear map is a morphism of associative algebras such that $$\wt\psi(M)=0, \wt\psi(1)=\wt\psi(z_1)=1.$$ 
Hence $M$ is a subset of the kernel of $\wt\psi$. 
The latter kernel is nonzero and thanks to Proposition~\ref{prop:open_subset_Z_n_cent_gen} we have that this kernel is generated by $M$ (i.e., is equal to $P$) and is primitive.
We have that $\bar c_i'=\wt\psi(c_i')$ form a basis of $\bar R_l$ and hence the restriction $\restr{\wt\psi}{Z_n}$ is injective. 
Thus $P\cap Z_n=\varnothing$ and $P$ satisfies all the desired properties. 
}

\okr{Lemma}{\label{Ltrdes} Let $\mathbb K_1\supset \mathbb K_2\supset\ldots$ be a sequence of fields such that, for all $n\ge1$, $\mathbb K_n$ is algebraically closed in $\mathbb K_1$. 
Pick $n\ge1$ and $x_1,\ldots, x_d\in\mathbb K_n$ such that $x_1, \ldots, x_d$ are algebraically dependent over $\mathbb K_l$ for all $l$. Then $x_1, \ldots, x_d$ are algebraically dependent over $\mathbb K=\bigcap_{n\ge1}\mathbb K_n$. }
\proof{}{ 
If $x_1, \ldots, x_{d-1}$ are algebraically dependent over $\mathbb K_l$ for all $l$ then we can replace $x_1, \ldots, x_d$ by $x_1, \ldots, x_{d-1}$. 
Thus we can assume that $x_1, \ldots, x_{d-1}$ are algebraically independent over $\mathbb K_l$ for all $l\ge l'$ for a positive integer $l'\ge n$. 

Hence $x_d$ is algebraic over $\mathbb K_l(x_1, \ldots, x_{d-1})$ for all $l\ge l'$.
Denote by $p_l$ the minimal monic polynomial of $x_d$ over $\mathbb K_l(x_1, \ldots, x_{d-1})$ for all $l\ge l'$. 
Denote by ${\rm Rts}_{l}$ the set of roots of $p_l$ in $\overline{\mathbb K_l(x_1, \ldots, x_{d-1})}$. 
It is clear that $p_l\mid p_{l+1}$ and hence ${\rm Rts}_{l}\subset {\rm Rts}_{l+1}$. 
This implies that $${\rm Rts}_l\subset\overline{\mathbb K_{l+1}(x_1, \ldots, x_{d-1})}.$$ 
Therefore the coefficients of $p_l$ belong to $\mathbb K_l(x_1, \ldots, x_{d-1})\cap \overline{\mathbb K_{l+1}(x_1, \ldots, x_{d-1})}$. 
This together with Lemma~\ref{Lvopros} implies that the coefficients of $p_l$ belong to $(\mathbb K_l\cap\overline{\mathbb K}_{l+1})(x_1, \ldots, x_{d-1})$. 
The assumption that $\mathbb K_{l+1}$ is algebraically closed in $\mathbb K_1$ implies that $\mathbb K_l\cap\overline{\mathbb K}_{l+1}=\mathbb K_{l+1}$. 
Therefore the coefficients of $p_l$ belong to $\mathbb K(x_1, \ldots, x_{d-1})$. This implies the desired result.}

\okr{Lemma}{\label{Lvopros}Let $\mathbb K_1\subset\mathbb K_2$ be a nested pair of fields. 
Fix $d'\ge0$ and a collection $x_1, \ldots, x_{d'}$ of algebraically independent over $\mathbb K_2$ variables. 
Identify $\overline{\mathbb K_1(x_1, \ldots, x_{d'})}$ with its image under the canonical embedding of $\overline{\mathbb K_1(x_1, \ldots, x_{d'})}$ into the algebraic closure $\overline{\mathbb K_2(x_1, \ldots, x_{d'})}$ of $\mathbb K_2(x_1, \ldots, x_{d'})$. 
Then $\overline{\mathbb K_1(x_1, \ldots, x_{d'})}\cap\mathbb K_2(x_1, \ldots, x_{d'})=(\overline{\mathbb K}_1\cap\mathbb K_2)(x_1, \ldots, x_{d'})$.}
\proof{}{ We will explicitly prove the statement of Lemma~\ref{Lvopros} under the assumption that $d'=1$ and it is easy to deduce the general case from this one by induction. 

Fix $q\in\mathbb K_2(x_1)\cap\overline{\mathbb K_1(x_1)}$. 
It is enough to show that $q\in(\mathbb K_2\cap\overline{\mathbb K}_1)(x_1)$. 
If $q=0$ then the latter statement is trivial. 
Assume $q\ne0$. 
Then $q$ can be expressed as $f/g$ where $f$, $g$ are relatively prime polynomials over $\mathbb K_2$ and $g$ is monic. 
The fact that $q$ is algebraic over $\mathbb K_1(x_1)$ implies that there exists a finite set $S$ such that $q(a)$ is well-defined and $q(a)\in\overline{\mathbb K}_1$ for all $a\in\overline{\mathbb K}_1\setminus S$. 

Consider the set of all pairs $(f', g')\in\overline{\mathbb K}_1[x]\oplus\overline{\mathbb K}_1[x]$ such that $f'(a)=q(a)g'(a)$ for all $a\in\overline{\mathbb K}_1\backslash S$. 
It is clear that the set of these pairs form a $\overline{\mathbb K}_1[x]$-submodule of $\overline{\mathbb K}_1[x]\oplus\overline{\mathbb K}_1[x]$; we denote this submodule~$M$. 
Let $M'$ be the projection of $M$ on the second summand. 
It is clear that $(f, g)\in M$ and thus that $M'$ is a nonzero ideal of $\overline{\mathbb K}_1[x]$ and hence it is generated by a single monic polynomyal $m_g\in\overline{\mathbb K}_1[x]$. 
Pick $m_f\in\overline{\mathbb K}_1[x]$ such that $(m_f, m_g)\in M$. 

We have $(f, g)\in M$ and thus there exists $l\in\overline{\mathbb K}_1[x]$ such that $g=m_gl$. 
Next, we have $$f(a)=q(a)m_g(a)l(a)=m_f(a)l(a)$$ for all $a\in\overline{\mathbb K}_1\setminus S$. 
This implies that $f-m_fl$ has at infinitely many roots and hence $f=m_fl$. 
Recall that $f$ and $g$ are relatively prime and hence $l\in\overline{\mathbb K}_1$. 
Next, both $g$ and $m_g$ are monic and hence $l=1$. 
This implies that $f, g\in(\mathbb K_2\cap\overline{\mathbb K}_1)[x]$ and therefore $q\in(\overline{\mathbb K}_1\cap\mathbb K_2)(x_1)$.}
The last ingredient in the proof of Proposition~\ref{Pprtalg} uses the skew field $Q(\nt; J)$ of $\Ua(\nt)/J$.
Set $T_l$ to be the subfield of $Q(\nt; J)$ generated by $Z_l, Z_{l+1},\ldots$ for all $l\ge 1$. 
It is clear that $T_l$ is commutative for all $l$. 
Further, put $\wt T_l$ to be the algebraic closure of $T_l$ in $T_1$. We will need the following lemma.
\okr{Lemma}{\label{Ltilt} If $Q(\nt; J)$ has trivial center then $\bigcap_l \wt T_l=\Cp$.}
\proof{}{ Assume to the contrary that $\bigcap_l \wt T_l\ne \Cp$ and fix $f\in\bigcap_l \wt T_l\setminus\Cp$. 
Also fix $a, b\in\Ua(\nt)/J$ with $f=a^{-1}b$. 
There exists $l'$ such that $a, b\in\Ua(\nt_{l'})/J_{l'}$. 
Consider $l\ge l'$. 
Recall that $Q_l$ is the quotient field of $Z_l$ and that $\Ua(\nt_l)/J_l\otimes_{Z_l}Q_l$ is isomorphic to the Weyl algebra $\Au_n(Q_l)$ over the field $Q_l$ for a positive integer $n\ge0$, see~\cite[Proposition~4.7.1, Lemma~6.4.5]{Dixmier1}. Thus we have a map
$$\Au_n(T_l)\cong \Au_n(Q_l)\otimes_{Q_l}T_l\to Q(\nt; J).$$
The first algebra is simple and therefore we can identify it with its image. 
It is clear that $f$ belongs to the quotient skew field of $\Au_n(T_l)$. 
Also, $f$ is algebraic over $T_l$ and hence $f\in T_l$~\cite[p. 295]{Goodearl2}. 
Therefore $f\in T_l$ for all $l\ge l'$.

This implies that $f$ commutes with the image of $\nt_l$ for all $l\ge l'$ and hence $f$ belongs to the center of $Q(\nt; J)$, i.e., $f\in\Cp$. 
This contradicts our assumption.}
\proof{ of Proposition~\ref{Pprtalg}}{ Assume to the contrary that $Z_n$ and $Z_l$ are not linearly disjoint for all $l>n$. 
Pick a transcendence basis $x_1, x_2, \ldots$ of $Z_n$. 
Thanks to the inequality $${\rm GKdim}(Z_n)\le {\rm GKdim}(\Ua(\nt_n))=\dim\nt_n,$$ this set is finite and we denote by $d$ the number of its elements. 
Then Lemma~\ref{Ltrdes} applied to $\mathbb K_l=\tilde T_l$ implies that $x_1, \ldots, x_d$ are algebraically dependent over $\bigcap_l\wt T_l=\Cp$, see Lemma~\ref{Ltilt}. This contradicts our assumption.
}
\newpage\sst{Socle Lie algebras}\label{sst:socle}\addcontentsline{toc}{subsection}{\ref{sst:socle}. Socle Lie algebras}
In this subsection we consider in details a special class of locally nilpotent Lie algebras, namely, socle Lie algebras. 
This class is much ``nicer'' then the class of all locally nilpotent Lie algebras. 
In particular, we can show in this case that there is a bijection between the set of orbits of a properly chosen (pro-)group and $\PSpec(\Sa(\nt))\approx\JSpec(\Ua(\nt))$. 
Also, in this case all such orbits are $\Sa(\nt)$-Zariski closed 
and a version of Dixmier--Moeglin equivalence holds, see Proposition~\ref{Pdmogsa}. 
\defi{A locally nilpotent Lie algebra $\nt$ is called \emph{socle} if it admits an exhaustion by its finite-dimensional nilpotent ideals.}
To state the above mentioned equivalence in a proper form we need a notion of a locally closed ideal both in Poisson and noncommutative setting. 
\defi{Let $I$ (respectively, $J$) be a Poisson (respectively, two-sided) ideal of $\Sa(\nt)$ (respectively, of $\Ua(\nt)$). We denote by ${\rm cl}(I)$ (respectively, by ${\rm cl}(J)$) the intersection of all radical Poisson (respectively, radical two-sided) ideals strictly containing $I$. We say that $I$ (respectively, $J$) is {\it locally closed} if $I\ne{\rm cl}(I)$ (respectively, if $J\ne{\rm cl}(J)$). For the topological meaning of locally closed ideals see~\cite[Lemma~II.7.7]{BGo}.}

We consider the following proposition as a Dixmier--Moeglin equivalence for the Poisson ideals of~$\Sa(\nt)$.
\okr{Proposition}{\textup{(cf.{~
\cite[Proposition~4.8.5]{Dixmier1}})}\label{Pdmogsa} Let $\nt$ be a socle Lie algebra and $I\subset\Sa(\nt)$ be a Poisson ideal. Then the following conditions are equivalent\textup:
\begin{equation*}
\begin{split}
&\text{\textup{i)} $I$ is maximal in the class of Poisson ideals\textup;}\\
&\text{\textup{ii)} $I$ is primitive};\\
&\text{\textup{iii)} the Poisson center of $\Sa(\nt)/I$ equals $\Cp$\textup;}\\
&\text{\textup{iv)} $I$ is prime and the Poisson center of the field of fractions of $\Sa(\nt)/I$ coincides with $\Cp$\textup;}\\
&\text{\textup{v)} $I$ is prime and locally closed.}\\
\end{split}
\end{equation*}}

\proof{}{ 
 Given a module $V$ of an arbitrary locally nilpotent Lie algebra~$\nt$, on which $\nt$ acts locally nilpotently, we denote by $V^{\nt}$ the submodule of $\nt$-invariants, i.e., $$V^{\nt}=\{v\in V\mid x(v)=0\text{ for all }x\in\nt\},$$ where by $x(v)$ we denote the image of $v$ under the linear operator corresponding to $x$. We note that $V^{\nt}\neq\{0\}$ for a finite-dimensional $V$. Indeed, the image of $\nt$ in ${\rm End}(V)$ is finite-dimensional and hence there exists a finite-dimensional subalgebra $\nt_k$ of $\nt$ which maps onto that image. By the Engel theorem we immediately obtain $V^{\nt}=V^{\nt_k}\ne 0$.

The implications $\text{\textup{(i)}}\implies\text{\textup{(ii)}},~ \text{\textup{(iv)}}\implies\text{\textup{(iii)}}$ are trivial. Next, $ \text{\textup{(i)}}$ implies $\text{\textup{(v)}}$ by~\cite[Proposition~1.4]{SQ}. 
Also $\text{\textup{(v)}}$ implies $\text{\textup{(ii)}}$ by~\cite[Proposition~1.7]{SQ}. 
Further, a minor modification of the proof of~\cite[Theorem~3.2]{BLSR} shows that $\text{\textup{(ii)}}\implies\text{\textup{(iv)}}$. 
Thus it is enough to show that $\text{\textup{(iii)}}$ implies $\text{\textup{(i)}}$ in our case.

Indeed, assume to the contrary that  there exists a Poisson ideal $I'\subset\Sa(\nt)$ with $I\subsetneq I'$ then the image $\bar I'$ of $I'$ in $\Sa(\nt)/I$ is a nontrivial Poisson ideal. 
Fix $0\ne f\in \bar I'$. 
Then there exists a finite-dimensional $\nt$-submodule $M_f$ of $\Sa(\nt)/I$ with $f\in M_f$. 
The condition that $\nt$ is a limit of nilpotent Lie algebras implies that $M_f^\nt\ne 0$. 
This contradicts $\text{\textup{(iii)}}$.
}

\nota{Fix an exhaustion of $\nt$ by its finite-dimensional ideals $\nt_k$ and a linear form $\lambda\in\nt^*$. 
We may assume without loss of generality that $\nt_1\subset\nt_2\subset\ldots$ and $\dim\nt_k=k$ for all $k\geq1$. 
Denote by~$\lambda_k$ the restriction of $\lambda$ to the ideal $\nt_k$ for $k\geq1$. 
Denote by $\rt_k\subset\nt_k$ the kernel of the bilinear from~$\beta_{\lambda_k}$, see Section~\ref{sst:orbit_method}: $\beta_{\lambda_k}(x,y)=\lambda_k([x,y])$, $x,~y\in\nt_k$. 
Put $$\pt=\sum_{k\geq1}\rt_k,~\pt_k=\sum_{i=1}^k\rt_i=\pt\cap\nt_k\text{ for each }k\geq1.$$
Then each $\pt_k$ is a polarization at $\lambda_k$, see Subsection~\ref{sst:orbit_method}.
Moreover, set $V_k=V(\nt_k,\pt_k,\lambda_k)$, then each $V_k$ is a simple $\nt_k$-module, and one can easily see that the natural map $V_k\to V_{k+1}$ is injective, and the corresponding simple $\nt$-module $V=\ilm V_k$ coincides with the $\nt$-module $V(\nt,\pt,\lambda)$ induced from the one-dimensional representation $x\mapsto\lambda(x)$ of $\pt$.}

\sst{The adjoint (pro-)group for a socle Lie algebras}\label{sst: adjg}\addcontentsline{toc}{subsection}{\ref{sst: adjg}. The adjoint (pro-)group for a socle Lie algebras}
Our next goal is to involve coadjoint orbits into the picture. First, one can consider the groups the $N_k$ such that $N_k=\Exp{\nt_k}$ and, for each $k$, the embedding $\nt_k\subset\nt_{k+1}$ induces the embedding $N_k\subset\ N_{k+1}$. Hence, we have the inductive limit $N=\ilm N_k$, and the group $N$ acts naturally on $\nt$ and on $\nt^*$. The point is that the group $N$ is ``too small" for our purposes.

\exam{(The countable-dimensional Heisenberg Lie algebra, cf. Example~\ref{Ehla}.) Let $\nt=\hei_\infty(\Cp)$ be the Lie algebra with generators $z$, $x_i$, $y_i$, $i\ge 1$, and relations
$$[x_i,y_i]=z\text{ for all }i,~[x_i,y_j]=0\text{ for }i\neq j,~[x_i,z]=[y_i,z]=0\text{ for all }i.$$
We call $\nt$ the \emph{countable-dimensional Heisenberg Lie algebra}. 
Clearly, $\nt$ is socle, because the subalgebra~$\nt_k$ generated by $z$ and $x_i$, $y_i$ for $1\leq i\leq k$ is an ideal.

There are two classes of primitive ideals of ${\rm U}(\nt)$ and ${\rm S}(\nt)$.

i) Every 
$\alpha\in\Cp^{\times}=\Cp\setminus\{0\}$ defines the two-sided (respectively, Poisson) ideal $J_\alpha$ (respectively, $I_\alpha$) of ${\rm U}(\nt)$ (respectively, ${\rm S}(\nt)$) generated by $z-\alpha$. 
It is easy to verify that
$$J_\alpha=J(\lambda)\iff I_\alpha=I(\lambda)\iff\lambda(z)=\alpha.$$

One can check that the quotient of $\Ua(\nt)$ by this ideal is a limit of Weyl algebras and hence is simple (a similar argument is applicable to the Poisson side). Thus, this ideal is maximal and hence primitive (respectively, Poisson primitive).

ii) Every $\lambda\in\nt^*$ with $\lambda(z)=0$ defines the ideal $J(\lambda)$ generated by $$x_i-\lambda(x_i),~y_i-\lambda(y_i),~i\ge 1,~\text{~and~}~z.$$ 
It is easy to verify that the quotient by $J(\lambda)$ is isomorphic to $\mathbb C$. Thus $J(\lambda)$ is maximal and hence primitive. 
Similar facts hold true in the Poisson case. 
It is easy to verify that, for $\mu\in\nt^*$,
$$J(\lambda)=J(\mu)\iff I(\lambda)=I(\mu)\iff\lambda=\mu.$$

The bijection between primitive Poisson ideals of $\Sa(\nt)$ and primitive ideals of $\Ua(\nt)$ is clear from these descriptions.

The adjoint group $N$ of $\hei_\infty(\Cp)$ (which is the injective limit of the adjoint groups of $\nt_k$) together with the action of $N$ on $\nt^*$ can be defined in a natural way. 
The group $N$ fixes all linear forms of type~(ii) and this matches the finite-dimensional case. 

On the other hand, consider a linear form $\lambda\in\nt^*$ with $\lambda(z)\ne0$, i.e., of type (i). 
Then the $N$-orbit~$N.\lambda$ of $\lambda$ consists of linear forms $\mu\in\nt^*$ such that 
\begin{center}$\lambda(z)=\mu(z)$ and $\lambda(x_i)=\mu(x_i), \lambda(y_i)=\mu(y_i)$ for all $i\gg0$\label{Exhla},\end{center}
i.e., the action of $N$ can change only finitely many coefficients $\lambda(x_i), \lambda(y_i)$.
Recall that $I(\mu)=I(\lambda)$ if and only if $\mu(z)=\lambda(z)$. 
This implies that $N.\lambda$ is much smaller than the set of linear forms $\mu$ for which $I(\lambda)=I(\mu)$. 

To fix this issue we construct out of $N$ a larger pro-group $\wh N\subset{GL}(\nt^*)$ such that the action of $\wh N$ can change all the coefficients $\lambda(x_i), \lambda(y_i)$ (using essentially the same formula for the adjoint action). 
Thus it is quite clear how to define such a group for Heisenberg Lie algebra and the construction suitable for every socle Lie algebra is given next.}

To define the pro-coadjoint action, we need some more notation. Since each $N_l$ acts on each $\nt_k^*$, we can denote by $\rho_k^l$ the morphism from $N_l$ to $\GL(\nt_k^*)$ which sends an element~$g$ to the corresponding linear operator on $\nt_k^*$ for $l\geq k$. Note that $\rho_k^l(N_l)\subset\rho_k^{l+1}(N_{l+1})$ for each $l\geq k$. Now, we put
\begin{equation*}
\wh N_k=\bigcup_{l\geq k}\rho_k^l(N_l).
\end{equation*}
Since each $\rho_k^l(N_l)$ is irreducible, there exists the minimal number $l(k)$ such that $$\wh N_k=\rho_k^{l(k)}(N_{l(k)}).$$ Obviously, there is the projection $\pi_{k+1}\colon\wh N_{k+1}\to\wh N_k$ for each $k$.
\defi{We set $\wh N:=\varprojlim \wh N_k$.}
Note that $\wh N$ is assembled out of the pieces of $N_k$s and therefore $\wh N$ acts on $\nt$ and thus on $\nt^*$. It is easy to verify that this construction does not depend on a choice of the sequence of ideals $\nt_1\subset\nt_2\subset\ldots$. 
Thus we consider $\wh N$ as ``the true adjoint group of $\nt$''. 
We claim that $\hat N$ is well-tailored for the orbit method and supplement this claim by Proposition~\ref{Porbmsoc} and Lemma~\ref{Lcadjc}.
\propp{\label{Porbmsoc}Pick $\lambda, \mu\in\nt^*$. Then $I(\lambda)=I(\mu)$ iff $\wh N.\lambda=\wh N.\mu$.}{First, we need to check that $I(g.\lambda)=I(\lambda)$ for each $g\in\wh N$. Thanks to Lemma~\ref{Lvaril} we have that $I(\lambda)$ can be constructed out of $I(\lambda_k)\subset\Sa(\nt_k)$ with $\lambda_k=\lambda|_{\nt_k}$. It is clear that $I(\lambda_k)=I((g.\lambda)_k)$.

It remains to show that if $I(\lambda)=I(\mu)$ then there exists $g\in \wh N$ such that $\lambda=g.\mu$.
As above, set $\lambda_k:=\restr{\lambda}{\nt_k}, \mu_k:=\restr{\mu}{\nt_k}$ for all $k\ge1$. Next, put $\wh N_k(\lambda_k\to\mu_k):=\{g\in\wh N_k\mid g.\lambda_k=\mu_k\}.$ It is clear that, for all $l\ge k\ge 1$, $\phi_{l, k}\colon\wh N_l\to \wh N_k$ maps $\wh N_l(\lambda_l\to\mu_l)$ to $\wh N_k(\lambda_k\to\mu_k)$, where $\phi_{l, k}=\pi_l\circ\pi_{l-1}\circ\ldots\circ\pi_{k+1}$. 
This implies that $\wh N_k(\lambda_k\to\mu_k)$ together with maps $\phi_{k, l}$ form a pro-variety; we denote this pro-variety by $\wh N(\lambda\to \mu)$. 
Clearly, each $\wh N_k(\lambda_k\to\mu_k)$ is nonempty and thus by Lemma~\ref{Lprone} the pro-variety $\wh N(\lambda\to\mu)$ has at least one point $g$. 
Evidently, $g\in\wh N$ and $g.\lambda=\mu$.}
We hope that an analogue of $\wh N$ can be constructed for every Lie algebra of a larger class of Lie algebras (say locally nilpotent Lie algebras) but we was not able to work out a relevant construction together with a relevant proof. 
Recall the definition of $\Sa(\nt)$-Zariski topology from Subsection~\ref{sst:pro_varieties}. We conclude this section by the following lemma.
\lemmp{\label{Lcadjc}Let $\Omega$ be a coadjoint $\wh N$-orbit on $\nt^*$. Then $\Omega$ is a closed \textup(in the $\Sa(\nt)$-Zariski topology\textup) pro-subvariety of $\nt^*$\textup, i.e.\textup, $\Omega$ equals to the common set of zeros in $\nt^*$ of $I(\lambda)$.}{Pick a linear form $\mu\in\Omega$. Assume to the contrary that $\Omega$ is not closed and pick $\mu'\in\overline{\Omega}\setminus \Omega$. 
Then we have $I(\mu)\subset I(\mu')$. 
Thanks to Proposition~\ref{Pdmogsa} we have that $I(\mu)$ is a maximal Poisson ideal of $\Sa(\nt)$ and hence $I(\mu)=I(\mu')$. 
Thus by Proposition~\ref{Porbmsoc} we have $\mu'\in\Omega$. 
This contradicts our assumption.}

\newpage\sect{Nil-Dynkin algebras}\label{sect:nildynkin}\addcontentsline{toc}{section}{\ref{sect:nildynkin}. Nil-Dynkin algebras}

This section is devoted to an important class of locally nilpotent Lie algebras, which we call\break nil-Dynkin algebras and consider as the most natural infinite-dimensional analogues of maximal nilpotent subalgebras of simple finite-dimensional Lie algebras from Subsection~\ref{sst:nilradicals_fin_dim}. For such algebras, the results of the previous section can be substantially improved.

More precisely, in Subsection~\ref{sst:nilradicals_inf_dim} we recall the definition of simple finitary infinite-dimensional Lie algebras and present a classification of their (splitting) Borel subalgebras due to I. Dimitrov and I.~Penkov. 
The splitting Borel subalgebras can be parametrized by linear orders on certain countable sets. 
To the contrary with the finite-dimensional setting, there are uncountably many nonisomorphic splitting Borel subalgebras. 
Next, we present the description of the center of $\Ua(\nt)$ of a nil-Dynkin algebra $\nt$ obtained by I. Penkov and the first author. This description is based on the infinite analogues of Kostant cascades and is similar to the finite-dimensional one from Example~\ref{exam:center_nilradical_fd}.

First, we prove an analogue of Theorem~\ref{theo:almost_all_Poisson_cent_gen_fin_dim} for nil-Dynkin algebras using the description of the center of $\Ua(\nt)$, see Theorem~\ref{theo:all_Poisson_cent_gen_ifd} in Subsection~\ref{sst:cent_gen_ideals_ifd}. Next, in Subsection~\ref{sst:critetion_non_trivial_ifd}, we deduce from this analogue a criterion for an ideal $I(\lambda)$ (or, equivalently, $J(\lambda)$) to be nontrivial. 
This criterion is given in terms of ranks of a certain ``submatrix'' of $\lambda$, see Theorem~\ref{theo:non_trivial_ifd} and Example~\ref{exam:nontrivial_ifd}.

\sst{Splitting Borel subalgebras of simple finitary Lie algebras}\label{sst:nilradicals_inf_dim}\addcontentsline{toc}{subsection}{\ref{sst:nilradicals_inf_dim}. Splitting Borel subalgebras of simple finitary Lie algebras}

Let $\slt_{\infty}(\Cp)$, $\sot_{\infty}(\Cp)$, $\spt_{\infty}(\Cp)$ be the three simple complex finitary countable-dimensional Lie algebras as classified by A.~Baranov \cite{Baranov1}. Each of them can be described as follows (see for example \cite{DimitrovPenkov2}). Consider an infinite chain of inclusions $$\gt_1\subset\gt_2\subset\ldots\subset\gt_n\subset\ldots$$ of simple Lie algebras, where all $\gt_n$ are of the same type $A$, $B$, $C$ or $D$, and $\rk{\gt_n}=n$ for types $B$, $C$, $D$, while $\rk\gt_n=n-1$ for type $A$. We fix such a chain in the following way. Let $\gt_n$ be the algebra denoted by $\slt_n(\Cp)$, $\sot_{2n+1}(\Cp)$, $\spt_{2n}(\Cp)$ or $\sot_{2n}(\Cp)$ in Subsection~\ref{sst:nilradicals_fin_dim}. The inclusions are trivial: they send the $(i,j)$th entry of a matrix to the $(i,j)$th entry of its image. 

Then the union $\gt=\bigcup\gt_n$ is isomorphic to $\slt_{\infty}(\Cp)$, $\sot_{\infty}(\Cp)$ or $\spt_{\infty}(\Cp)$. We choose the nested Cartan subalgebras $\htt_n\subset\gt_n$, $\htt_n\subset\htt_{n+1}$, as in Subsection~\ref{sst:nilradicals_fin_dim}: $\htt_n$ is the set of all diagonal matrices from $\gt_n$. Then each root space of $\gt_n$ is mapped to exactly one root space of $\gt_{n+1}$. The union $\htt=\bigcup\htt_n$ acts semisimply on $\gt$, and it is by definition a \emph{splitting Cartan subalgebra} of~$\gt$ \cite[Section 3]{DPSn}. We have a root decomposition $\gt=\htt\oplus\bigoplus_{\alpha\in\Phi}\gt^{\alpha}$ where $\Phi\subset\htt^*$ is \emph{the root system of} $\gt$ with respect to $\htt$ and $\gt^{\alpha}$ are the \emph{root spaces}. The root system $\Phi$ can be thought as the union of the root systems of $\gt_n$ and equals one of the following infinite root systems:
\begin{equation*}
\begin{split}
A_{\infty}&=\pm\{\epsi_i-\epsi_j,~i,j\in\Zp_{>0},~i<j\},\\
B_{\infty}&=\pm\{\epsi_i-\epsi_j,~i,j\in\Zp_{>0},~i<j\}\cup\pm\{\epsi_i+\epsi_j,~i,j\in\Zp_{>0},~i<j\}\cup\pm\{\epsi_i,~i\in\Zp_{>0}\},\\
C_{\infty}&=\pm\{\epsi_i-\epsi_j,~i,j\in\Zp_{>0},~i<j\}\cup\pm\{\epsi_i+\epsi_j,~i,j\in\Zp_{>0},~i<j\}\cup\pm\{2\epsi_i,~i\in\Zp_{>0}\},\\
D_{\infty}&=\pm\{\epsi_i-\epsi_j,~i,j\in\Zp_{>0},~i<j\}\cup\pm\{\epsi_i+\epsi_j,~i,j\in\Zp_{>0},~i<j\}.
\end{split}
\end{equation*}

Here we first present a classification of splitting Borel subalgebras of $\gt$, which is due to I. Penkov and I. Dimitrov \cite{DimitrovPenkov1}. A \emph{splitting Borel subalgebra} of $\gt$ is a subalgebra $\bt$ such that, for every~$n$, $\bt_n=\bt\cap\gt_n$ is a Borel subalgebra of $\gt_n$. It is well known that any splitting Borel subalgebra is conjugate via $\mathrm{Aut}\,\gt$ to a splitting Borel subalgebra containing $\htt$. Therefore, in what follows we restrict ourselves to considering only such Borel subalgebras $\bt$.

Recall \cite{DimitrovPenkov1} that a linear order on $\{0\}\cup\{\pm\epsi_i\}$ is $\Zp_2$-\emph{linear} if multiplication by $-1$ reverses the order. By \cite[Proposition 3]{DimitrovPenkov1}, there exists a bijection between splitting Borel subalgebras of $\gt$ containing $\htt$ and certain linearly ordered sets as follows.
\begin{equation*}\predisplaypenalty=0
\begin{split}
&\text{For }A_{\infty}\text{: linear orders on }\{\epsi_i\};\\
&\text{for }B_{\infty}\text{ and }C_{\infty}\text{: }\Zp_2\text{-linear orders on }\{0\}\cup\{\pm\epsi_i\};\\
&\text{for }D_{\infty}\text{: }\Zp_2\text{-linear orders on }\{0\}\cup\{\pm\epsi_i\}\text{ with the property that}\\
&\text{a minimal positive element (if it exists) belongs to }\{\epsi_i\}.
\end{split}
\end{equation*}
In the sequel we denote these linear orders by $\preceq$. To write down the above bijection, denote $\teta_i=\epsi_i$, if $\epsi_i\succ0$, and $\teta_i=-\epsi_i$, if $\epsi_i\prec0$ (for $A_{\infty}$, $\teta_i=\epsi_i$ for all $i$). Then put $\bt=\htt\oplus\nt$, where $\nt=\bigoplus\limits_{\alpha\in\Phi^+}\gt^{\alpha}$ and
\begin{equation*}\predisplaypenalty=0
\begin{split}
A_{\infty}^+&=\{\teta_i-\teta_j,~i,j\in\Zp_{>0},~\teta_i\succ\teta_j\},\\
B_{\infty}^+&=\{\teta_i-\teta_j,~i,j\in\Zp_{>0},~\teta_i\succ\teta_j\}\cup\{\teta_i+\teta_j,~i,j\in\Zp_{>0},
~\teta_i\succ\teta_j\}\cup\{\teta_i,~i\in\Zp_{>0}\},\\
C_{\infty}^+&=\{\teta_i-\teta_j,~i,j\in\Zp_{>0},~\teta_i\succ\teta_j\}\cup\{\teta_i+\teta_j,~i,j\in\Zp_{>0},~
\teta_i\succ\teta_j\}\cup\{2\teta_i,~i\in\Zp_{>0}\},\\
D_{\infty}^+&=\{\teta_i-\teta_j,~i,j\in\Zp_{>0},~\teta_i\succ\teta_j\}\cup\{\teta_i+\teta_j,~i,j\in\Zp_{>0},~\teta_i\succ\teta_j\}.\\
\end{split}.
\end{equation*}
Actually, we are interested only in isomorphism classes of these subalgebras. Consequently, we will assume without loss of generality that $\teta_i=\epsi_i$, i.e., that $\epsi_i\succ0$ for all $i$ and all $\Phi$. This is possible because each $\Zp_2$-linear order on $\{0\}\cup\{\pm\epsi_i\}$ is clearly isomorphic to a $\Zp_2$-linear order with the property $\epsi_i\succ0$ for all $i$.

\defi{The Lie algebras $\nt$ defined above are called \emph{nil-Dynkin Lie algebras}.}

Our next goal is to recall the description of the center $Z(\nt)$ of the enveloping algebra $U(\nt)$ from~\cite{IgnatyevPenkov1}. Fix $\nt$, i.e., fix an order $\preceq$ as above. Define the subset $\Nu\subseteq\Zp_{>0}$ by setting $\Nu=\bigcup_{k\geq0}\Nu_k$, where $\Nu_0=\varnothing$ and $\Nu_k$ for $k\geq1$ is defined inductively in the following table.
\begin{center}
\begin{tabular}{|l|l|}
\hline
$\Phi$&$\Nu_k$\\
\hline\hline
$A_{\infty}$&$\Nu_{k-1}\cup\{i,j\}$ if there exists a maximal element $\epsi_i$\\
&and a minimal element $\epsi_j$ of $\{\epsi_s,~s\in\Zp_{>0}\setminus\Nu_{k-1}\}$,\\
&$\Nu_{k-1}$ otherwise\\
\hline
$C_{\infty}$&$\Nu_{k-1}\cup\{i\}$ if there exists a maximal element $\epsi_i$ of $\{\epsi_s,~s\in\Zp_{>0}\setminus\Nu_{k-1}\}$,\\
&$\Nu_{k-1}$ otherwise\\
\hline
$B_{\infty}$,&$\Nu_{k-1}\cup\{i,j\}$ if there exists a maximal element $\epsi_i$ of $\{\epsi_s,~s\in\Zp_{>0}\setminus\Nu_{k-1}\}$\\
$D_{\infty}$&and a maximal element $\epsi_j$ of $\{\epsi_s,~s\in\Zp_{>0}\setminus\left(\Nu_{k-1}\cup\{i\}\right)\}$,\\
&$\Nu_{k-1}$ otherwise\\
\hline
\end{tabular}
\end{center}

\exam{Let $\Phi=A_{\infty}$. If $\epsi_1\succ\epsi_3\succ\ldots\succ\epsi_4\succ\epsi_2$, then $\Nu=\Zp_{>0}$. If $\epsi_1\succ\epsi_2\succ\epsi_3\succ\ldots$, then $\Nu=\varnothing$. On the other hand, if $\Phi\neq A_{\infty}$ and $\epsi_1\succ\epsi_2\succ\ldots\succ0\succ\ldots\succ-\epsi_2\succ-\epsi_1$, then $\Nu=\Zp_{>0}$.}

Now we can define the (possibly infinite) Kostant cascade corresponding to $\nt$. Namely, to each $\Nu_k$ such that $\Nu_{k-1}\subsetneq\Nu_k$, we assign the following root:
\begin{equation*}
\beta_k=\begin{cases}\epsi_i-\epsi_j,&\text{if $\Phi=A_{\infty}$ and $\Nu_k\setminus\Nu_{k-1}=\{i,j\}$, $\epsi_i\succ\epsi_j$,}\\
\epsi_i+\epsi_j,&\text{if $\Phi=B_{\infty}$ or $D_{\infty}$ and $\Nu_k\setminus\Nu_{k-1}=\{i,j\}$, $\epsi_i\succ\epsi_j$,}\\
2\epsi_i,&\text{if $\Phi=C_{\infty}$ and $\Nu_k\setminus\Nu_{k-1}=\{i\}$.}\\
\end{cases}
\end{equation*}

\defi{Put $\Bu=\{\beta_k,~k\geq1,~\Nu_{k-1}\subsetneq\Nu_k\}$. The subset $\Bu$ is called the \emph{Kostant cascade} corresponding to $\nt$.} 
Note that $\Bu$ is a strongly orthogonal subset of $\Phi^+$; however it is not necessarily maximal with this property (cf. Example~\ref{exam:center_nilradical_fd}).
For instance, if $\Phi=A_{\infty}$ and $\epsi_1\succ\epsi_3\succ\ldots\succ\epsi_4\succ\epsi_2$, then\break $\Bu=\{\epsi_1-\epsi_2,~\epsi_3-\epsi_4,~\epsi_5-\epsi_6,~\ldots\}$; if $\Phi\neq A_{\infty}$ and $\epsi_1\succ\epsi_2\succ\ldots\succ0\succ\ldots\succ-\epsi_2\succ-\epsi_1$, then
\begin{equation*}
\Bu=\begin{cases}\{\epsi_1+\epsi_2,~\epsi_3+\epsi_4,~\epsi_5+\epsi_6,~\ldots\}&\text{for $B_{\infty}$ and $D_{\infty}$},\\
\{2\epsi_1,~2\epsi_2,~2\epsi_3,~\ldots\}&\text{ for $C_{\infty}$}.
\end{cases}
\end{equation*}

To each finite nonempty subset $M\subset\Zp_{>0}$, one can assign a root subsystem $\Phi_M$ of $\Phi$ and a subalgebra~$\nt_M$ of $\nt$ by putting
\begin{equation*}
\begin{split}
\Phi_M&=\Phi\cap\langle\epsi_i,~i\in M\rangle_{\Rp},\\
\nt_M&=\bigoplus_{\alpha\in\Phi_M^+}\gt^{\alpha},~\Phi_M^+=\Phi_M\cap\Phi^+.
\end{split}
\end{equation*}
Then the subsystem $\Phi_M$ is isomorphic to the root system $\Phi_n$ of $\gt_n$ for $n=|M|$; we fix the isomorphism  $j_M\colon\Phi_n\to\Phi_M$ induced by $\epsi_i\mapsto\epsi_{a_i}$, where $M=\{a_1,\ldots,a_n\}$, $\epsi_{a_1}\succ\ldots\succ\epsi_{a_n}$. 
Put $\unt_n$ to be the maximal nilpotent subalgebra of $\gt_n$ considered in Subsection~\ref{sst:nilradicals_fin_dim}; then $\nt_M\cong\unt_n$. 

Note also that $\nt=\ilm\nt_M$. Here, for $M\subseteq M'$, the monomorphism ${i_{M,M'}\colon\nt_M\hookrightarrow\nt_{M'}}$ is just the inclusion. Further, our chain of embeddings of $\gt_n$ automatically inherits the basis $\{e_{\alpha},~\alpha\in\Phi^+\}$ of~$\nt$ consisting of root vectors $e_{\alpha}$ defined in Subsection~\ref{sst:nilradicals_fin_dim}. 
The linear map $\phi_M\colon\unt_n\to\nt_M$, $e_{\alpha}\mapsto e_{j_M(\alpha)}$, $M\subset\Zp_{>0},~n=|M|$, is the isomorphism of Lie algebras.
Consider a finite subset $M'\subset\Zp_{>0}$ with $M\subset M'$ and $n':=|M'|$. The inclusion $i_{M,M'}\colon\nt_M\hookrightarrow\nt_{M'}$ induces the unique embedding of Lie algebras $\kappa_{M, M'}\colon\unt_n\hookrightarrow\unt_{n+1}$ such that the following diagram is commutative.
$$
\begin{CD}
\unt_n@>{\kappa_{M, M'}}>>\unt_{n'}\\
@V{\phi_M}VV @VV{\phi_{M'}}V\\
\nt_M@>{i_{M, M'}}>>\nt_{M'}
\end{CD}
$$

We are now ready to write down a set of generators of $Z(\nt)$. 
To each root $\beta$ from the Kostant cascade $\Bu$ we will attach an element of $Z(\nt)$ and altogether these elements will be the free generators of~$Z(\nt)$. 
Pick $k\le |\Bu|$ and $\beta=\beta_k$. 
Let $M$ be a finite subset of $\Zp_{>0}$ such that $\Nu_k\subseteq M$. The isomorphism $\phi_M$ gives rise to an isomorphism $\phi_M\colon U(\unt_n)\to U(\nt_M)$, $n=|M|$. 
Put
$$\Delta_{\beta}:=\phi_M(\sigma(\xi_{j_M^{-1}(\beta)})),$$ 
see Example~\ref{exam:center_nilradical_fd} for the notation. 
Observe that $\Delta_{\beta}\in U(\nt_M)$ is given by one of formulas (\ref{formula:Delta_i_A_n}), (\ref{formula:Delta_i_C_n}) or (\ref{formula:Delta_i_B_n_D_n_pf}) with $e_{j_M(\alpha)}$ instead of $e_{\alpha}$. It is important that $\Delta_{\beta}\in U(\nt_M)$ depends only on $\beta$ but not on $M$. Moreover, it follows from the finite-dimensional theory that $\Delta_{\beta}$ belongs to the center of~$U(\nt)$. The next theorem was proved in \cite[Theorem 2.6]{IgnatyevPenkov1}.

\mtheo{If \label{theo:center_ifd} $\Phi=A_{\infty}$\textup{,} $C_{\infty}$\textup{,} $B_{\infty}$ or $D_{\infty}$\textup{,} then $\Delta_{\beta}$\textup{,} $\beta\in\Bu$\textup{,} are algebraically independent and generate $Z(\nt)$ as an algebra. In particular\textup{,} $Z(\nt)$ is a polynomial ring in $|\Bu|$ variables.}

\sst{Centrally generated Poisson ideals}\label{sst:cent_gen_ideals_ifd}\addcontentsline{toc}{subsection}{\ref{sst:cent_gen_ideals_ifd}. Centrally generated Poisson ideals}

One can easily deduce from Theorem~\ref{theo:center_ifd} that the Poisson center $Y(\nt)$ of $\Sa(\nt)$ is generated as an algebra by $\sigma^{-1}(\Delta_{\beta})$ for $\beta\in\Bu$. This follows immediately from the fact that, given a finite subset $M\subset\Zp_{>0}$, one has $$Y(\nt)\cap\Sa(\nt_M)=\bigcap_{M'\supseteq M}Y(\nt_{M'}),~Z(\nt)\cap\Ua(\nt_M)=\bigcap_{M'\supseteq M}Z(\nt_{M'}),$$ and the fact that the restriction of $\sigma$ to $Y(\nt_M')$ is an algebra isomorphism onto $Z(\nt_M')$.

In Subsection~\ref{sst:cent_gen_ideals_fin_dim} we have showed that, for a finite-dimensional nilpotent Lie algebra, almost all primitive ideals in $\Ua(\nt)$ are centrally generated.  
Similarly, all primitive Poisson ideals in $\Sa(\nt)$ are generated as ideals by their intersections with the Poisson center of $\Sa(\nt)$, see Theorem~\ref{theo:almost_all_Poisson_cent_gen_fin_dim}. 
Below we will prove an analog of this result for the class of nil-Dynkin algebras.

\theop{There exists\label{theo:all_Poisson_cent_gen_ifd} an open dense \textup(with respect to the countable-Zariski topology\textup)\break subset of $\nt^*$ such that $I(\lambda)$ is generated as an ideal by its intersection with the Poisson center $Y(\nt)$ of~$\Sa(\nt)$ for each $\lambda$ from this subset.}{Pick a root $\beta\in\Phi^+$ (possibly, not from $\Bu$) and a finite subset $M\subset\Zp_{>0}$ for which $\beta\in\Phi_M$. Let $|M|=n$, then we have the isomorphisms $j_M\colon\Phi_n\to\Phi_M$ and $\phi_M\colon\unt_n\to\nt_M$. Assume in addition that $j_M^{-1}(\beta)$ belongs to the Kostant cascade of $\Phi_n$ and denote by $\xi_{\beta}^M$ the image in $\Sa(\nt_M)$ of the corresponding canonical generator of $Y(\unt_n)$ under the isomorphism $\Sa(\unt_n)\to\Sa(\nt_M)$ induced by~$\phi_M$. 
From now on we will consider $\xi_{\beta}^M$ as an element of $\Sa(\nt)$. Note that this element depends on~$M$, not only on $\beta$, but for a given $\beta$ there exist only countably many possibilities for $\xi_{\beta}^M$. 
For example, if $\Phi=A_{\infty}$, $\beta=\epsi_i-\epsi_j$ for some $\epsi_i\succ\epsi_j$ and $M=\{i_1,~\ldots,~i_k,~i,~j,~j_k,~\ldots,~j_1\}$, where $\epsi_{i_1}\succ\ldots\succ\epsi_{i_k}\succ\epsi_i\succ\epsi_j\succ\epsi_{j_k}\succ\ldots\succ\epsi_{j_1}$, then one has
\begin{equation*}
\xi_{\beta}^M=\begin{vmatrix}
e_{i_1,j}&e_{i_1,j_k}&\ldots&e_{i_1,j_1}\\
e_{i_2,j}&e_{i_2,j_k}&\ldots&e_{i_2,j_1}\\
\vdots&\vdots&\iddots&\vdots\\
e_{i_k,j}&e_{i_k,j_k}&\ldots&e_{i_k,j_1}\\
e_{i,j}&e_{i,j_k}&\ldots&e_{i,j_1}\\
\end{vmatrix}=\begin{vmatrix}
e_{\epsi_{i_1}-\epsi_j}&e_{\epsi_{i_1}-\epsi_{j_k}}&\ldots&e_{\epsi_{i_1}-\epsi_{j_1}}\\
e_{\epsi_{i_2}-\epsi_j}&e_{\epsi_{i_2}-\epsi_{j_k}}&\ldots&e_{\epsi_{i_2}-\epsi_{j_1}}\\
\vdots&\vdots&\iddots&\vdots\\
e_{\epsi_{i_k}-\epsi_j}&e_{\epsi_{i_k}-\epsi_{j_k}}&\ldots&e_{\epsi_{i_k}-\epsi_{j_1}}\\
e_{\epsi_i-\epsi_j}&e_{\epsi_i-\epsi_{j_k}}&\ldots&e_{\epsi_i-\epsi_{j_1}}\\
\end{vmatrix}.
\end{equation*}
For other root systems, $\xi_{\beta}^M$ can be defined by formulas (\ref{formula:Delta_i_C_n}) and (\ref{formula:Delta_i_B_n_D_n_pf}), see Example~\ref{exam:center_nilradical_fd}.

Denote by $A$ the subset of $\nt^*$ consisting of all $\lambda$ such that $\xi_{\beta}^M(\lambda)\neq0$ for all $\beta\in\Phi^+$ and all $M\subset\Zp_{>0}$ for which $\xi_{\beta}^M$ exists. Evidently, $A$ is open (and so dense by Proposition~\ref{prop:Zariski_countable_irr_v2}) in the countable-Zariski topology on $\nt^*$. Hence it is enough to check that $I(\lambda)$ is generated by $I(\lambda)\cap Y(\nt)$ for each $\lambda\in A$.

Pick a linear form $\lambda\in A$, an element $t'\in I(\lambda)$ and a finite subset $M\subset\Zp_{>0}$ such that $t'\in\Sa(\nt_M)$. 
If $n=|M|$, we denote by $t$ the preimage of $t'$ in $\Sa(\unt_n)$ under the isomorphism $\Sa(\unt_n)\to\Sa(\nt_M)$ induced by $\phi_M\colon\unt_n\to\nt_M$. Denote by $\Bu_M$ the preimage in $\Phi_n$ of $\Bu\cap\Phi_M$ under the isomorphism $j_M\colon\Phi_n\to\Phi_M$. 
Denote also by $\Bu_n$ the Kostant cascade of $\Phi_n$. 
Without loss of generality we can enlarge $M$ in such a way that $\Bu\cap\Phi_M=\{\beta_1,~\ldots,~\beta_k\}$ for some $k$ (if not, we will just add the required indices to~$M$). 
Then, clearly, $\Bu_M\subset\Bu_n$. Denote by $\xi_{\beta}$, $\beta\in\Bu_n$, the canonical generators of $Y(\unt_n)$. 
Since $\lambda\in A$, we have $\xi_{\beta}(\phi_M^*(\restr{\lambda}{\nt_M}))=c_{\beta}\neq0$ for all $\beta\in\Bu_n$. According to (\ref{formula:cent_gen_fin_dim}) and Theorem~\ref{theo:almost_all_Poisson_cent_gen_fin_dim}, $t$ belongs to the (Poisson) ideal of $\Sa(\unt_n)$ generated by $\xi_{\beta}-c_{\beta}$ for $\beta\in\Bu_n$. It is enough to prove that $t$~belongs to the ideal of $\Sa(\unt_n)$ generated by $\xi_{\beta}-c_{\beta}$ for $\beta\in\Bu_M$. (Note that, for $\beta\in\Bu_M$, the image of $\xi_{\beta}$ in $\Sa(\nt)$ depends only on $\beta$, not on the choice of a sufficiently large subset $M$.)

Let $\beta_0$ be the minimal root from $\Bu_n\setminus\Bu_M$ (if $\Bu_n\setminus\Bu_M$ is empty then there is nothing to prove). Note that
\begin{equation*}
j_M(\beta_0)=\begin{cases}\epsi_{i_0}-\epsi_{j_0}&\text{for }\Phi=A_{\infty},\\
2\epsi_{i_0}&\text{for }\Phi=C_{\infty},\\
\epsi_{i_0}+\epsi_{j_0}&\text{for }\Phi=B_{\infty}\text{ or }D_{\infty}\\
\end{cases}
\end{equation*}
for certain $i_0,~j_0\in\Zp_{>0}$. We will proceed case by case.

\textsc{Case} $\Phi=A_{\infty}$. Here $j_M^{-1}(\beta_s)=\epsi_s-\epsi_{n-s+1}$ and $\beta_s=\epsi_{i_s}-\epsi_{j_s}$ for certain $i_s,~j_s$, $1\leq s\leq k$. Since $\beta_0\notin\Bu_M$, at least one of the following two conditions is satisfied: $\epsi_{i_0}$ is not the maximal element (with respect to the partial order $\preceq$) of the set $\left\{\epsi_i,~i\in\Zp_{>0}\setminus\bigcup_{s=1}^k\{i_s,~j_s\}\right\}$, or $\epsi_{j_0}$ is not the minimal element element of this set. Suppose that $\epsi_{i_0}$ is not maximal (the case when $\epsi_{j_0}$ is not minimal can be considered similarly). Let $i'\in\Zp_{>0}$ be such that $\epsi_{i_k}\succ\epsi_{i'}\succ\epsi_{i_0}$ (note that $i'\notin M$). Put $M'=M\cup\{i'\}$, then $\Phi_M\subset\Phi_{M'}$. We can fix the embedding $\restr{j_{M'}^{-1}}{\Phi_M}\circ j_M\colon\Phi_n\hookrightarrow\Phi_{n+1}$. This induces the embeddings $\kappa_{M,M'}\colon\unt_n\hookrightarrow\unt_{n+1}$ and $\Sa(\unt_n)\hookrightarrow\Sa(\unt_{n+1})$, hence we can consider $t$ as an element $\Sa(\unt_{n+1})$. We denote also by $\Bu_{n+1}$ the Kostant cascade of $\Phi_{n+1}$, then $\Bu_n\cap\Bu_{n+1}=\Bu_M$.

For example, let $n=8$ and the root vectors from $\unt_{n+1}\setminus\unt_n$ be marked gray in the picture below. Then $\Bu_{n+1}\cap\Bu_n=\Bu_M=\{\epsi_1-\epsi_9,~\epsi_2-\epsi_8\}$ (the corresponding root vectors are marked by $\otimes$), $\Bu_{n+1}\setminus\Bu_n=\{\epsi_3-\epsi_7,~\epsi_4-\epsi_6\}$ (marked by $\bullet$) and $\Bu_n\setminus\Bu_{n+1}=\{\epsi_4-\epsi_7,~\epsi_5-\epsi_6\}$ (marked by $\square$).
\begin{center}
\includegraphics{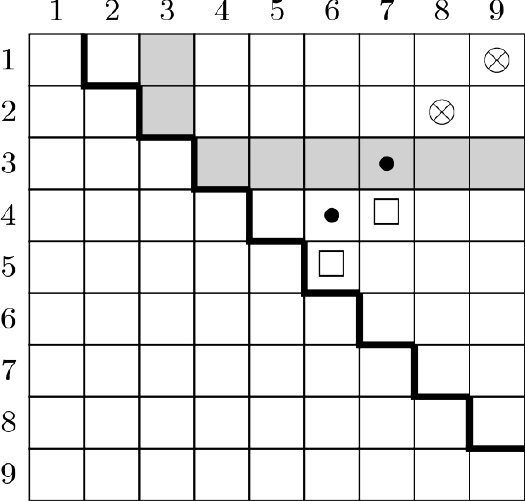}
\end{center}

Denote by $\lambda'$ the image of $\restr{\lambda}{\nt_{M'}}$ under the isomorphism $\phi_{M'}^*\colon\nt_{M'}^*\to\unt_{n+1}^*$. Put $N'=\Exp{\unt_{n+1}}$. One can consider $\Sa(\unt_{n+1})$ as the algebra $\Cp[\unt_{n+1}^*]$ of polynomial functions on $\unt_{n+1}^*$. Since $t'\in I(\lambda)$, one has $t(g.\lambda')=0$ for all $g\in N'$, where $g.\lambda'$ denotes the result of the coadjoint action. It follows, e.g., from \cite[Theorems 4, 8]{Kostant2} and \cite[Theorem 3.1]{IgnatyevPenkov1} (see also Example~\ref{exam:regular_orbits_A_n}) that $\lambda\in A$ implies that the coadjoint $N'$-orbit $N'.\lambda'$ of the linear form $\lambda'$ contains the unique linear form $\mu\in\unt_{n+1}^*$ such that, given $\beta\in\Phi_{n+1}^+$, $\mu(e_{\beta})\neq0$ if and only if $\beta\in\Bu_{n+1}$.

We will need the exponential map 
$$\exp\colon\unt_{n+1}\to N',~x\mapsto\sum_{i=0}^{\infty}\dfrac{x^i}{i!}.$$
The coadjoint action is defined by the formula $$((\exp x).f)(y)=\sum_{i=0}^{\infty}\dfrac{f((-\ad{x})^i(y))}{i!},~x,y\in\unt_{n+1},~f\in \unt_{n+1}^*,$$ where $\ad{x}(y)=[x,y]$. 
Now, put $m=|\Bu_n|$ and pick numbers $g_i\in\Cp^{\times}$ for $k+1\leq i\leq m$. 
Set
\begin{equation*}
g=\exp(g_me_{\epsi_{n-m+1}-\epsi_{n-m+2}})\ldots\exp(g_{k+1}e_{\epsi_{n-(k+1)+1}-\epsi_{n-(k+1)+2}})\in N'.
\end{equation*}
Put also $\nu=g.\mu$, then one can easily deduce that, for all $\alpha\in\Phi_{n+1}^+$,
\begin{equation}
\nu(e_{\alpha})=\begin{cases}g_i\mu(\epsi_i-\epsi_{n-i+1}),&\text{if }\alpha=\epsi_i-\epsi_{n-i+2}\text{ for }i>k,\\
\mu(\alpha)&\text{otherwise}.
\end{cases}\label{formula:nu_A}
\end{equation}
Note that $\epsi_i-\epsi_{n-i+2}$, $k+1\leq i\leq m$, are exactly the roots from $\Bu_n\setminus\Bu_M$ considered as a subset of the root system $\Phi_{n+1}$.

We see that $\xi_{\beta_i}(\nu)=\xi_{\beta_i}(\mu)=\xi_{\beta_i}(\lambda')=c_{\beta_i}$ for $1\leq i\leq k$, i.e., for $\beta_i\in\Bu_M$. 
At the same time, since $g_i$ are arbitrary nonzero numbers, formula (\ref{formula:Delta_i_A_n}) implies that the values of all other canonical generators of $Y(\unt_n)$ on $\nu$ can be arbitrary nonzero. 
Now, denote by $R$ the subset of $\unt_n^*$ consisting of all linear forms $f\in\unt_n^*$ such that the values of all canonical generators of $\Sa(\unt_n)$ on $f$ are nonzero. By \cite[Theorems 4, 8]{Kostant2} and \cite[Theorem 3.1]{IgnatyevPenkov1}, $R$ is an open dense $\Exp{\unt_n}$-stable subset of $\unt_n^*$. Furthermore, each coadjoint $\Exp{\unt_n}$-orbit in $R$ contains exactly one element of the form $\sum_{\beta\in\Bu_n}r_{\beta}e_{\beta}^*$, $r_{\beta}\neq0$, where $\{e_{\alpha}^*,~\alpha\in\Phi_n^+\}$ is the dual basis for the basis $\{e_{\alpha},~\alpha\in\Phi_n^+\}$ of $\unt_n$. Now, set $$Z=\{f\in\unt_n^*\mid\xi_{\beta}(f)=c_{\beta}\text{ for all }\beta\in\Bu_M\}.$$
Obviously, $Z$ is a closed $\Exp{\unt_n}$-stable subset of $\nt_n^*$. It is clear from (\ref{formula:Delta_i_A_n}) that $Z$ is isomorphic to an affine space (and so is irreducible), $\nu\in Z$ and that when $g_i$'s run over $\Cp^{\times}$ independently, $\nu$ runs over a complete system or representatives of $\Exp{\unt_n}$-orbits in $Z\cap R$. Since $t'\in I(\lambda)$, one has $t(\nu)=t(\mu)=t(\lambda')=0$. It follows that $\restr{t}{Z\cap R}\equiv0$ and thus $t$ vanishes on $Z$. But the latter condition means exactly that $t$ belongs to the ideal of $Z$ in $\Sa(\unt_n)$, i.e., to the ideal of $\Sa(\unt_n)$ generated by $\xi_{\beta}-c_{\beta}$, $\beta\in\Bu_M$, as required.

\textsc{Case} $\Phi=C_{\infty}$. The idea is similar to the case $A_{\infty}$, but the technical details differ slightly. Here $j_M^{-1}(\beta_s)=2\epsi_s$ and $\beta_s=2\epsi_{i_s}$ for certain $i_s$, $1\leq s\leq k$. Since $\beta_0\notin\Bu_M$, $\epsi_{i_0}$ is not maximal (with respect to the partial order $\preceq$) in the set $\left\{\epsi_i,~i\in\Zp_{>0}\setminus\{i_1,~\ldots,~i_k\}\right\}$. Let $i'\in\Zp_{>0}$ be such that $\epsi_{i_k}\succ\epsi_{i'}\succ\epsi_{i_0}$ (note that $i'\notin M$). Now we define $M'$, $\Phi_{M'}$, $\unt_{n+1}$, $\Bu_{n+1}$, the group $N'$, the linear form $\lambda'$ and the embeddings $\kappa_{M,M'}\colon\unt_n\hookrightarrow\unt_{n+1}$, $\Sa(\unt_n)\hookrightarrow\Sa(\unt_{n+1})$ as in the previous case. Note that $\Bu_{n+1}=\Bu_n\cup\{2\epsi_{k+1}\}$. Again, \cite[Theorem~3.1]{IgnatyevPenkov1} shows that $N'.\lambda'$ contains a linear form $\mu$ such that, given $\beta\in\Phi_{n+1}^+$, $\mu(e_{\beta})\neq0$ if and only if $\beta\in\Bu_{n+1}$.

Now, we pick numbers $g_i\in\Cp^{\times}$, $k+2\leq i\leq n+1$, and put $$g=\prod_{i=k+2}^{n+1}\exp(g_ie_{\epsi_{k+1}-\epsi_i})\in N'$$ ($i$~runs from $k+2$ to $n+1$), $\nu=g.\mu\in\nt_{n+1}^*$. Using the formula for the coadjoint action, one can easily check that, given $\alpha\in\Phi_{n+1}^+$,
\begin{equation*}
\nu(e_{\alpha})=\begin{cases}g_i\mu(2\epsi_{k+1}),&\text{if }\alpha=\epsi_{k+1}+\epsi_i,~k+1\leq i\leq n+1,\\
g_ig_j\mu(2\epsi_{k+1}),&\text{if }\alpha=\epsi_i+\epsi_j,~k+2\leq i<j\leq n+1,\\
g_i^2\mu(2\epsi_{k+1})+\mu(2\epsi_i),&\text{if }\alpha=2\epsi_i,~k+1\leq i\leq n+1,\\
\mu(\alpha)&\text{otherwise}
\end{cases}
\end{equation*}
(cf. (\ref{formula:nu_A})). Hence, using formula (\ref{formula:Delta_i_C_n}), one can calculate by the induction on $i$, $k+1\leq i\leq n$, that the value of the $i$th canonical generator of $Y(\unt_n)$ on $\nu$ is equal up to a sign to
$$\prod_{s=1}^k\mu(2\epsi_s)\left(\prod_{j=k+2}^{i+1}\mu(2\epsi_j)+\mu(2\epsi_{k+1})\sum_{j=k+2}^{i+1}g_j^2\prod_{k+2\leq r\leq i+1,~r\neq j}\mu(2\epsi_r)\right).$$
Since $g_i$ are arbitrary nonzero numbers, we see that the values of all these canonical generators of $Y(\unt_n)$ on $\nu$ can be arbitrary nonzero.
Thus, one can complete the proof arguing literally as in the last paragraph of the case $A_{\infty}$ (with formula (\ref{formula:Delta_i_C_n}) instead of (\ref{formula:Delta_i_A_n})).

\textsc{Case} $\Phi=D_{\infty}$. Here $j_M^{-1}(\beta_s)=\epsi_{2s-1}+\epsi_{2s}$ and $\beta_s=\epsi_{i_s}+\epsi_{j_s}$ for certain $i_s$, $j_s$, $1\leq s\leq k$. Since $\beta_0\notin\Bu_M$, there exists $i'\in\Zp_{>0}\setminus\{i_1,~\ldots,~i_k\}$ such that either $\epsi_{i_k}\succ\epsi_{i'}\succ\epsi_{i_0}$ or $\epsi_{i_0}\succ\epsi_{i'}\succ\epsi_{j_0}$ (note that $i'\notin M$). Now we define $M'$, $\Phi_{M'}$, $\unt_{n+1}$, $\Bu_{n+1}$, the group $N'$, the linear form $\lambda'$ and the embeddings $\kappa_{M,M'}\colon\unt_n\hookrightarrow\unt_{n+1}$, $\Sa(\unt_n)\hookrightarrow\Sa(\unt_{n+1})$ as in the case of $A_{\infty}$. Note that $\Bu_{n+1}\cap\Bu_n=\Bu_M$. It follows from \cite[Theorem~2.4]{Ignatyev1} that $N'.\lambda'$ contains a linear form $\mu$ such that, given $\beta\in\Phi_{n+1}^+$, $\mu(e_{\beta})\neq0$ if and only if $\beta\in\Bu_{n+1}$.

Now, we pick numbers $g_i\in\Cp^{\times}$, $k+1\leq i\leq [n/2]=|\Bu_n|$, and put $$g=g'\prod_{i=k+2}^{[n/2]}\exp(g_ie_{\epsi_{2i-1}-\epsi_{2i+1}})\in N'$$ ($i$~runs from $[n/2]$ to $k+1$), $\nu=g.\mu\in\nt_{n+1}^*$, where
\begin{equation*}
g'=\begin{cases}
\exp(g_{k+1}e_{\epsi_{2k+1}-\epsi_{2k+1}}),&\text{if }\epsi_{i_k}\succ\epsi_{i'}\succ\epsi_{i_0},\\
\exp(g_{k+1}e_{\epsi_{2k}-\epsi_{2k+1}}),&\text{if }\epsi_{i_0}\succ\epsi_{i'}\succ\epsi_{j_0}.\\
\end{cases}
\end{equation*}
Using the formula for the coadjoint action, one can easily check that, given $\alpha\in\Phi_{n+1}^+$,
\begin{equation*}
\nu(e_{\alpha})=\begin{cases}g_i\mu(\epsi_{2i-1}+\epsi_{2i}),&\text{if }\alpha=\epsi_{2i}+\epsi_{2i+1},~k+2\leq i\leq [n/2],\\
g_{k+1}\mu(\epsi_{2k-1}+\epsi_{2k}),&\text{if }\alpha=\epsi_{2k}+\epsi_{2k+1}\text{ and }\epsi_{i_k}\succ\epsi_{i'}\succ\epsi_{i_0},\\
&\text{or }\alpha=\epsi_{2k-1}+\epsi_{2k+1}\text{ and }\epsi_{i_0}\succ\epsi_{i'}\succ\epsi_{j_0},\\
\mu(\alpha)&\text{otherwise}
\end{cases}
\end{equation*}
(cf. (\ref{formula:nu_A})). It is clear from formula (\ref{formula:Delta_i_B_n_D_n_pf}) that the value of the $i$th canonical generator of $Y(\unt_n)$ on $\nu$, $k+1\leq i\leq[n/2]$, is equal up to a sign to $\prod_{s=1}^k\mu(\epsi_{2s-1}+\epsi_{2s})\prod_{j=k+1}^ig_j\mu(\epsi_{2j-1}+\epsi_{2j})$ and so can be arbitrary nonzero. Thus, one can complete the proof arguing literally as in the last paragraph of the case $A_{\infty}$ (with (\ref{formula:Delta_i_B_n_D_n_pf}) instead of (\ref{formula:Delta_i_A_n})).

\textsc{Case} $\Phi=B_{\infty}$. The proof for $B_{\infty}$ is completely similar to the case $D_{\infty}$, so we omit it.}



\sst{\label{sst:critetion_non_trivial_ifd}Nontriviality criterion for primitive ideals}\addcontentsline{toc}{subsection}{\ref{sst:critetion_non_trivial_ifd}. Nontriviality criterion for primitive ideals}

In this subsection,  we fix a nil-Dynkin algebra $\nt$ together with the corresponding partial order~$\preceq$ as in Subsection~\ref{sst:nilradicals_inf_dim}. 
For such a Lie algebra $\nt$ we provide a necessary and sufficient condition for a primitive Poisson ideal $I(\lambda)$ to be nontrivial (Theorem~\ref{theo:non_trivial_ifd}). 
There are two ways to state this condition. 
The first one is as follows: there exists a countable collection $\{ \tilde{}\Xi_k\}_{k}$ of countable collections of polynomials from $\Sa(\nt)$ such that $I(\lambda)\ne0$ if and only if there exists $k$ for which $\lambda(\xi)=0$ for all $\xi\in\Xi_k$. 
The advantage of this form is that it is quite believable that the result in this form can be generalized to a larger class of Lie algebras (probably, to the entire class of locally nilpotent Lie algebras). 
The second form is based on the structure of nil-Dynkin Lie algebras and, in particular, on the fact that all $\Xi_k$ can be chosen in such a way that all $\xi\in\Xi_k$ are ``minors'' in an appropriate sense. 
Moreover, the condition $\lambda(\xi)=0$ for all $\xi\in \Xi_k$ is a degeneration condition for a certain submatrix of the infinite matrix defined by $\lambda$, see Example~\ref{exam:nontrivial_ifd}. 
Note that $I(\lambda)\ne0$ if and only if $J(\lambda)\ne0$ and hence the criterion also is a nontriviality criterion for primitive ideals of~$\Ua(\nt)$.  
To proceed, we need more notation related to the above mentioned minors.



First, let $\Phi=A_{\infty}$. Pick two sequences of  positive integers $I=\{i_1,~\ldots,~i_k\}$ and $J=\{j_1,~\ldots,~j_k\}$ so that $\epsi_{i_1}\succ\ldots\succ\epsi_{i_k}\succ\epsi_{j_k}\succ\ldots\succ\epsi_{j_1}$. We denote by $\xi_I^J$ the element of the symmetric algebra $\Sa(\nt)$ defined by the following rule:
\begin{equation}
\xi_I^J=\begin{vmatrix}
e_{i_1,j_k}&\ldots&e_{i_1,j_2}&e_{i_1,j_1}\\
e_{i_2,j_k}&\ldots&e_{i_2,j_2}&e_{i_2,j_1}\\
\vdots&\iddots&\vdots&\vdots\\
e_{i_k,j_k}&\ldots&e_{i_k,j_2}&e_{i_k,j_1}\\
\end{vmatrix}=\begin{vmatrix}
e_{\epsi_{i_1}-\epsi_{j_k}}&\ldots&e_{\epsi_{i_1}-\epsi_{j_2}}&e_{\epsi_{i_1}-\epsi_{j_1}}\\
e_{\epsi_{i_2}-\epsi_{j_k}}&\ldots&e_{\epsi_{i_2}-\epsi_{j_2}}&e_{\epsi_{i_2}-\epsi_{j_1}}\\
\vdots&\iddots&\vdots&\vdots\\
e_{\epsi_{i_k}-\epsi_{j_k}}&\ldots&e_{\epsi_{i_k}-\epsi_{j_2}}&e_{\epsi_{i_k}-\epsi_{j_1}}\\
\end{vmatrix}.\label{formula:xi_I_J_A_infty}
\end{equation}
For instance, let $\beta_s=\epsi_{i_s}-\epsi_{j_s}$, $1\leq s\leq k$, be the first $k$ roots of the Kostant cascade of the Lie algebra~$\nt$.
Then, clearly, $\Delta_{\beta_k}=\sigma(\xi_I^J)$ belongs to the set of generators of~$Z(\nt)$ defined by Theorem~\ref{theo:center_ifd}.

Next, assume $\Phi=C_{\infty}$. Consider $I=\{i_1,~\ldots,~i_k\}$, $J=\{j_1,~\ldots,~j_k\}$ with $\epsi_{i_1}\succ\ldots\succ\epsi_{i_k}$ and $\epsi_{j_1}\succ\ldots\succ\epsi_{j_k}$. If $i\neq j$ for all $i\in I$, $j\in J$, then we put
\begin{equation*}
\xi_I^J=\begin{vmatrix}\label{formula:xi_I_J_C_infty}
e_{\epsi_{i_1}+\epsi_{j_k}}&\ldots&e_{\epsi_{i_1}+\epsi_{j_2}}&e_{\epsi_{i_1}+\epsi_{j_1}}\\
e_{\epsi_{i_2}+\epsi_{j_k}}&\ldots&e_{\epsi_{i_2}+\epsi_{j_2}}&e_{\epsi_{i_2}+\epsi_{j_1}}\\
\vdots&\iddots&\vdots&\vdots\\
e_{\epsi_{i_k}+\epsi_{j_k}}&\ldots&e_{\epsi_{i_2}+\epsi_{j_k}}&e_{\epsi_{i_1}+\epsi_{j_k}}\\
\end{vmatrix}.
\end{equation*}
If $i=j$ for some pairs $i\in I$, $j\in J$ then we use a modification of this formula in which we replace $e_{\epsi_i+\epsi_j}$ by~$2e_{2\epsi_i}$. For instance, let $\beta_s=2\epsi_{i_s}$, $1\leq s\leq k$, be the first $k$ roots from the Kostant cascade of~$\nt$. Then $\Delta_{\beta_k}=\sigma(\xi_I^I)$ belongs to the set of generators of $Z(\nt)$ from Theorem~\ref{theo:center_ifd}.

Finally, for $\Phi=B_{\infty}$ or $D_{\infty}$, given a sequence $I=\{i_1,~\ldots,~i_{2k}\}$ with $\epsi_{i_1}\succ\ldots\succ\epsi_{i_{2k}}$, we denote by $\xi_I$ the unique element of $\Sa(\nt)$ such that
\begin{equation*}
\xi_I^2=\pm\begin{vmatrix}
e_{\epsi_{i_1}+\epsi_{i_{2k}}}&\ldots&e_{\epsi_{i_1}+\epsi_{i_3}}&e_{\epsi_{i_1}+\epsi_{i_2}}&0\\
e_{\epsi_{i_2}+\epsi_{i_{2k}}}&\ldots&e_{\epsi_{i_2}+\epsi_{i_3}}&0&-e_{\epsi_{i_1}+\epsi_{i_2}}\\
e_{\epsi_{i_3}+\epsi_{i_{2k}}}&\ldots&0&-e_{\epsi_{i_2}+\epsi_{i_3}}&-e_{\epsi_{i_1}+\epsi_{i_3}}\\
\vdots&\iddots&\vdots&\vdots&\vdots\\
0&\ldots&-e_{\epsi_{i_3}+\epsi_{i_{2k}}}&-e_{\epsi_{i_2}+\epsi_{i_{2k}}}&-e_{\epsi_{i_1}+\epsi_{i_{2k}}}\\
\end{vmatrix}.
\end{equation*}(After suitable reordering of indices, the matrix in the right-hand side becomes skew-symmetric, so $\xi_I$ is nothing but its Pfaffian. Our normalization is such that the term $e_{\epsi_{i_1}+\epsi_{i_2}}e_{\epsi_{i_3}+\epsi_{i_4}}\ldots e_{\epsi_{i_{2k-1}}+\epsi_{i_{2k}}}$ enters $\xi_I$ with coefficient 1.) For instance, let $\beta_s=\epsi_{i_{2s-1}}+\epsi_{i_{2s}}$, $1\leq s\leq k$, be the first $k$ roots from the Kostant cascade of $\nt$. Then $\Delta_{\beta_k}=\sigma(\xi_I)$ belongs to the set of generators of $Z(\nt)$ from Theorem~\ref{theo:center_ifd}.

Given an upper-right pair $p$ and a positive integer $k$, we define the ideal $I(p,k)$ in the following table. Note that if $\Phi=B_{\infty}$ or $D_{\infty}$ and $i=m\neq-j$ then only $k=1$ is allowed. 
\begin{center}
\begin{longtable}{|l|l|l|}
\hline
Type of $\Phi$ & Upper-right pair $p$ & Generators of $I(p, k)$ as an ideal\\
\hline\hline
$A_{\infty}$ & $(i,j)$ & $\xi_I^J$, $|I|=|J|=k$, $\epsi_{i_k}\succeq\epsi_i$, $\epsi_j\succeq\epsi_{j_k}$\\
\hline
$C_{\infty}$ & $(i,-i)$ & $\xi_I^J$, $|I|=|J|=k$, $\epsi_{i_k}\succeq\epsi_i$, $\epsi_{j_k}\succeq\epsi_i$\\
\hline
$B_{\infty}$ or $D_{\infty}$ & $(i,-i)$ & $\xi_I$, $|I|=2k$, $\epsi_{i_{2k}}\succeq\epsi_i$\\
\hline
$B_{\infty}$ or $D_{\infty}$ & $(i,j)$, $i=m\neq-j$ & $\xi_{\{m,s\}}=e_{\epsi_m+\epsi_s}$, $\epsi_s\succ\epsi_j$\\
\hline
\end{longtable}
\end{center}
The ideals $I(p, k)$ might seems artificial from the first look but their zeroes sets are quite conceptual, see Corollary~\ref{coro:ideal_to_rank}.

Note, however, that $I(p,k)$ may be zero. Namely, for $A_{\infty}$, the ideal $I(p,k)$ is zero if and only if the partial order $\preceq$ has the maximal element $\epsi_{i_0}$ or the minimal element $\epsi_{j_0}$ and at least one of the numbers $|\{s\in\Zp_{>0}\mid\epsi_{i_0}\succeq\epsi_s\succeq\epsi_i\}|$, $|\{s\in\Zp_{>0}\mid\epsi_j\succeq\epsi_s\succeq\epsi_{j_0}\}|$ is finite and less than $k$. Similarly, for~$C_{\infty}$, $I(p,k)=\{0\}$ if and only if $\preceq$ has the maximal element $i_0$ and the number $|\{s\in\Zp_{>0}\mid\epsi_{i_0}\succeq\epsi_s\succeq\epsi_i\}|$ is finite and less than $k$. Finally, for $B_{\infty}$ and $D_{\infty}$, $I(p,k)=\{0\}$ if and only if $\preceq$ has the maximal element, $p=(i,-i)$ and the number $|\{s\in\Zp_{>0}\mid\epsi_{i_0}\succeq\epsi_s\succeq\epsi_i\}|$ is finite and less than $2k$.

\propp{The\label{prop:ideal_p_k_Poisson_ifd} ideal $I(p,k)$ is prime and Poisson.}{Pick a finite set $M\subset\Zp_{>0}$, $|M|=n$, containing all positive entries of $p$. Assume in addition that $n\geq k$ for $\Phi=A_{\infty}$ or $C_{\infty}$, and $n\geq2k$ for $\Phi=B_{\infty}$ or $D_{\infty}$. For $A_{\infty}$ and $C_{\infty}$ (respectively, $B_{\infty}$ and $D_{\infty}$), we denote by~$I_M(p,k)$ the ideal of $\Sa(\unt_n)$ generated by the preimages in $\Sa(\unt_n)$ under the isomorphism $\Sa(\unt_n)\to\Sa(\nt_M)$ of all $\xi_I^J$ with $I,~J\subseteq M$ (respectively, of all $\xi_I$ with $I\subseteq M$).
It is enough to check that the ideal $I_M(p,k)$ of $\Sa(\unt_n)$ is prime and Poisson for all such $M$. 
Indeed, after identifying $I_M(p,k)$ with its image in $\Sa(\nt_M)$ one has $$I(p,k)\cap\Sa(\nt_M)=\bigcup_{M'\supseteq M}(I_{M'}(p,k)\cap\Sa(\nt_M)).$$

Note that $\unt_n$ is a subspace of the space $\ut$ of all strictly upper-triangular matrices with zeroes on the diagonal (in particular, for $\Phi=A_{\infty}$, $\unt_n=\ut$); see Subsection~\ref{sst:nilradicals_fin_dim} for the details. Furthermore, the group $N_n=\Exp{\unt_n}=\{\exp(x),~x\in\unt_n\}$ is a subgroup of the group $\Uu$ of all upper-triangular matrices with units on the diagonal, where $\exp(x)=\sum_{l=0}^{\infty}x^l/l!$ (in particular, for $\Phi=A_{\infty}$, $N_n=\Uu$).

On the other hand, $\unt_n^*$ can be identified with the space $\unt_n^T$ (the superscript~$T$ stands for the transposed matrix) by putting $\mu(x)=\tr(\mu x)$ for $\mu\in\unt_n^T$, $x\in\unt_n$. It is easy to check (see, e.g.,\break \cite[Lecture 7, Section 1]{Kirillov2}) that under this identification the coadjoint action has the form\break $g.\mu=(g\mu g^{-1})_{\low}$, where, given a matrix $x$,
\begin{equation*}
(x_{\low})_{i,j}=\begin{cases}0,&\text{if }i\leq j,\\
x_{i,j},&\text{if }i>j.\\
\end{cases}
\end{equation*}
Recall the basis $\{e_{\alpha}^*,~\alpha\in\Phi_n^+\}$ of $\unt_n^*$ and the numeration of rows and columns of matrices from Subsection~\ref{sst:nilradicals_fin_dim}. Under the identification $\unt_n^*\cong\unt_n^T$,
\begin{equation}
e_{\alpha}^*=\begin{cases}
e_{\alpha}^T&\text{if }\Phi=A_{\infty},\text{ or }\Phi=C_{\infty}\text{ and }\alpha=2\epsi_i,\\
e_{\alpha}^T/4&\text{if }\Phi=B_{\infty}\text{ and }\alpha=\epsi_i,\\
e_{\alpha}^T/2&\text{otherwise}.
\end{cases}\label{formula:n_*_n_T_dual_basis}
\end{equation}
Given a matrix $\mu$ from $\unt_n^T$, we denote by $[\mu]_j^i$ the submatrix of $\mu$ which entries are situated nonstrictly to the South-West from the $(j,i)$th position. Then, for $B_{\infty}$ and $D_{\infty}$ (respectively, for $C_{\infty}$), $[\mu]_{-i}^i$ is antisymmetric (respectively, symmetric) with respect to the antidiagonal.

Now, denote by $\ut_j^i$ the Zariski closed subset of $\ut^T$ consisting of all matrices $f$ from $\ut^T$ such that $\rk[f]_j^i<2k$ (respectively, $\rk[f]_j^i<k$) if $\Phi=B_{\infty}$ or $D_{\infty}$ and $j=-i$ (respectively, otherwise). One can easily check that this subset is invariant under the coadjoint action of~$\Uu$ (see, e.g., the proof of \cite[Lemma 2.2]{Ignatyev2}). Indeed, $\Uu$ is generated as a group by elements of the form $x=\exp(te_{a,b})=1+te_{a,b}$, where $a<b$, $t\in\Cp$, and $1$ denotes the identity matrix. Given $f\in\ut^T$, one can easily deduce that
\begin{equation*}
(x.f)_{r,s}=\begin{cases}f_{a,s}+tf_{b, s},&\text{if
}r=a\text{ and }1\leq s<a,\\
f_{r,b}-tf_{r, a},&\text{if
}s=b\text{ and }b<r\leq n,\\
f_{r,s}&\text{otherwise}.
\end{cases}
\end{equation*}
Hence if $r>a$ and $s<b$, then $[x.f]_r^s=[f]_r^s$. If
$r\leq a$ (and so $s<r\leq a<b$), then the $a$th row of
$[x.f]_r^s$ is obtained from the $a$th row of
$[f]_r^s$ by adding the $b$th row of $[f]_r^s$
multiplied by $t$. Similarly, if $s\geq b$ (and so $r>s\geq
b>a$), then the~$b$th column of $[x.f]_r^s$ is obtained from
the~$b$th column of $[f]_r^s$ by subtracting the $a$th
column of $[f]_r^s$ multiplied by $t$. In both cases,
$\rk[x.f]_r^s=\rk[f]_r^s$.

Next, denote by $Z$ the set of common zeroes of polynomials from $I_M(p,k)$ in $\unt_n^*$. We claim that $Z$ is an $N_n$-invariant subset of $\unt_n^*$ (i.e., that $\sqrt{I_M(p,k)}$ is Poisson). To prove this fact it is enough to note that $Z=\unt_n^T\cap\ut_j^i$ if $p=(i,~j)$. (Here we replace the entries of $p$ by the corresponding integers from $-n$ to $n$). For $A_{\infty}$ and $C_{\infty}$, $\mu\in Z$ if and only if all $k\times k$ minors of~$[\mu]_j^i$ are zero, as required. For $B_{\infty}$ and $D_{\infty}$, if $i=m\neq-j$ then the equality $Z=\unt_n^T\cap\ut_j^i$ is obvious; in other cases, $[\mu]_{-i}^i$ is antisymmetric with respect to the antidiagonal, and all its principal (with respect to the antidiagonal) $k\times k$ minors are zero. The expansion formula for Pfaffians shows that all principal minors of $[\mu]_{-i}^i$ of size at least $k$ are zero. By the well-known Principal Minor Theorem, we see that $\rk[\mu]_{-i}^i<k$.

It remains to show that $I_M(p,k)$ is a prime ideal. For $A_{\infty}$, see, e.g., \cite[Corollary 16.29]{MillerSturmfels}. For $B_{\infty}$ and $D_{\infty}$, if $i=m\neq-j$ then it is evident; otherwise it was proved in \cite{Avramov}. For $C_{\infty}$, see, e.g.,\break \cite[Theorem 1]{Kutz1}.}

Pick a linear form $\lambda\in\nt^*$ and an upper-right pair $p$.
Denote by $[\lambda]_p$ the matrix defined by the following rule. For $A_{\infty}$ and $p=(i,~j)$, the rows and the columns of this matrix are indexed by the numbers $i'\in\Zp_{>0},~\epsi_{i'}\succeq\epsi_i$ and $j'\in\Zp_{>0},~\epsi_{j'}\preceq\epsi_j$ respectively (it is possible that this matrix has infinite sizes). For other root systems and $p=(i,~j)$, the rows and the columns of this matrix are indexed by the numbers $i'\in\Zp_{>0}$,~$\epsi_{i'}\succeq\epsi_i$ and $-j'\in\Zp_{<0}$, $\epsi_{j'}\succeq\epsi_{-j}$ respectively. By definition, if $\Phi=A_{\infty}$ then the~$(i',~j')$th entry of $[\lambda]_p$ equals $\lambda(e_{\epsi_{i'}-\epsi_{j'}})$. For other root systems, the $(i',~-j')$th entry equals~$\lambda(e_{\epsi_{i'}+\epsi_{j'}})$ if $i'\neq j'$. Finally, put
\begin{equation*}
([\lambda]_p)_{i',-i'}=\begin{cases}2\lambda(e_{2\epsi_{i'}})&\text{for }\Phi=C_{\infty},\\
0&\text{for }\Phi=B_{\infty}\text{ and }D_{\infty}.\\
\end{cases}
\end{equation*}
For $A_{\infty}$ (respectively, $C_{\infty}$), $\xi_I^J(\lambda)$ is exactly the minor of the matrix $[\lambda]_p$ with the set of rows $I$ and the set of columns $J$ (respectively, $-J$), while for $B_{\infty}$ and $D_{\infty}$, if $j=-i$ then $\xi_I^2(\lambda)$ equals the minor of $[\lambda]_p$ with the row set $I$ and the column set $-I$.

We define the rank of $[\lambda]_p$ as the cardinality of a maximal set of linearly independent rows (columns), written $\rk[\lambda]_p$. 
It is clear that  $[\lambda]_p<k$ if and only if the rank of each submatrix of~$[\lambda]_p$ of finite sizes is less than $k$. 
If $\rk[\lambda]_p<\infty$, then we say that \emph{the rank of $[\lambda]_p$ is finite}.

\corop{\label{coro:ideal_to_rank}The set of common zeroes of an ideal $I(p,k)$ has the form $\{\lambda\in\nt^*\mid\rk[\lambda]_p<k'\}$, where
\begin{equation*}
k'=\begin{cases}2k&\text{for $B_{\infty}$\textup, $D_{\infty}$ with $p=(i,-i)$},\\
k&\text{otherwise}.
\end{cases}
\end{equation*}
}{A minor modification of the proof of Proposition~\ref{prop:ideal_p_k_Poisson_ifd} provides the desired fact.}
\defi{We say that $\rk[\lambda]_p$ \emph{is not maximal} if $\rk[\lambda]_p$ is finite and either both of the sizes of $[\lambda]_p$ are infinite, or one of the sizes is finite and the least size is greater than $\rk[\lambda]_p$.}

Our next main result is as follows.

\theop{Pick\label{theo:non_trivial_ifd} a linear form $\lambda\in\nt^*$. If $\Bu\neq\varnothing$ then $I(\lambda)\neq0$. If $\Bu=\varnothing$ then the following conditions on $I(\lambda)$ are equivalent\textup{:}
\begin{equation*}
\begin{split}
&\text{\textup{i)} $I(\lambda)\neq0$};\\
&\text{\textup{ii)} $I(\lambda)\supseteq I(p,k)$ for some upper-right pair $p$ and some $k\in\Zp_{>0}$ such that $I(p,k)\neq0$};\\
&\text{\textup{iii)} $\rk[\lambda]_p$ is not maximal for some upper-right pair $p$}.\\
\end{split}
\end{equation*}}{According to Theorem~\ref{theo:center_ifd}, if $\Bu\neq\varnothing$ then $Y(\nt)\neq\Cp$. 
Denote by $${\rm ev}_\lambda\colon \Sa(\nt)\to\Cp,\hspace{10pt} \nt\ni z\mapsto \lambda(z)$$ the evaluation map at $\lambda$.
Then $z-\lambda(z)\in I(\lambda)$ for $z\in Y(\nt)\setminus\Cp$ and hence $I(\lambda)\ne0$. 
So, we will consider the case when $\Bu=\varnothing$ and $Z(\nt)=Y(\nt)=\Cp$ (see the first paragraph of Subsection~\ref{sst:cent_gen_ideals_ifd}). Implication $\mathrm{(ii)}\Longrightarrow\mathrm{(i)}$ is evident. Implication $\mathrm{(iii)}\Longrightarrow\mathrm{(ii)}$ follows immediately from Corollary~\ref{coro:ideal_to_rank}. It remains to check that $\mathrm{(i)}\Longrightarrow\mathrm{(iii)}$.

For $A_{\infty}$ (respectively, for $C_{\infty}$; for $B_{\infty}$ and $D_{\infty}$), denote by $\Fo\subset\Sa(\nt)$ the set of all $\xi_I^J$ (respectively, of all $\xi_I^I$; of all $\xi_I$); clearly, $\Fo$ is countable. Further, a linear form $\mu\in\nt^*$ belongs to the dense open set $A$ from the proof of Theorem~\ref{theo:all_Poisson_cent_gen_ifd} if and only if $f(\mu)\neq0$ for all $f\in\Fo$. Assume that $K$ is a prime Poisson ideal of $\Sa(\nt)$ such that $K\cap\Fo=\varnothing$. The fact that $K$ is prime implies that $\Sa(\nt)/K$ is a domain. Since $K\cap\Fo=\varnothing$, we see that $S'=(\Sa(\nt)/K)[f^{-1},~f\in\Fo]$ is nonzero. Hence, there exists a nontrivial maximal ideal $M'$ of $S'$. It is evident that $S'$ is countable-dimensional and hence $M'$ is of codimension~1 in $S'$ by Corollary~\ref{Lschur}. Let $\mu\in\nt^*$ be the point corresponding to $M'$. Then we have $f(\mu)\neq0$ for all $f\in\Fo$. But Theorem~\ref{theo:all_Poisson_cent_gen_ifd} implies that $I(\mu)=0$. On the other hand we have $K\subset I(\mu)$ and hence $K=0$. In particular, we conclude from $I(\lambda)\neq0$ that $f\in I(\lambda)$ for some $f\in\Fo$.

\textsc{Case} $\Phi=A_{\infty}$. We have $\xi_I^J\in I(\lambda)$ for some sequences $I$, $J$ of the same size $k$. The proof depends on the structure of the partial order $\preceq$ defining $\nt$. Since $\Bu=\varnothing$, we conclude that this partial order can not have both the minimal and the maximal elements.

Consider the case when it has neither minimal nor maximal element. Since $\{e_{\epsi_a-\epsi_b},e_{\epsi_b-\epsi_c}\}=e_{\epsi_a-\epsi_c}$ for distinct indices $a,~b,~c$, the row expansion formula for determinants implies that $\{e_{\epsi_a-\epsi_b},\xi_I^J\}=0$ except the case when either $b\in I$ (and then $a\notin J$) or $a\in J$ (and then $b\notin I$). If $b\in I$ and $a\notin J$ then
\begin{equation}
\{e_{\epsi_a-\epsi_b},\xi_I^J\}=\pm\xi_{I'}^J,\label{formula:commuting_minors}
\end{equation}
where $I'$ is obtained from $I$ by replacing $b$ by $a$ and reordering the indices. Similarly, if $a\in J$ and $b\notin I$ then $\{e_{\epsi_a-\epsi_b},\xi_I^J\}=\pm\xi_I^{J'}$, where $J'$ is obtained from $J$ by replacing $a$ by $b$ and reordering the indices. We conclude that $\xi_{I'}^{J'}\in I(\lambda)$ and so $\xi_{I'}^{J'}(\lambda)=0$ for all $I'$, $J'$ such that $|I'|=|J'|=k$ and $\epsi_{i_k'}\succeq\epsi_{i_1}$, $\epsi_{j_1}\succeq\epsi_{j_k'}$, because $I(\lambda)$ is a Poisson ideal. But this means exactly that the rank of $[\lambda]_p$ is less than $k$, where $p=(i_1,~j_1)$, i.e., $\rk[\lambda]_p$ is finite and so is not maximal.

Next, consider the case when the maximal element exists, and so the minimal element does not. (The case when the minimal element exists can be considered similarly.) We denote the maximal element by $\epsi_{m_1}$. Let $\Mo$ be the subset of $\Zp_{>0}$ consisting of positive integers $m$ such that the set $\{m'\in\Zp_{>0}\mid\epsi_{m_1}\succeq\epsi_{m'}\succeq\epsi_m\}$ is finite. Suppose $\Mo$ is infinite and $m_1,~\ldots,~m_k\in\Mo$ are such that $\epsi_{m_1}\succ\ldots\succ\epsi_{m_k}$ are the largest $k$ elements among all $\epsi_s,~s\in\Zp_{>0}$. It is clear that applying (\ref{formula:commuting_minors}) we can obtain $\xi_{M_0}^J\in I(\lambda)$, where $M_0=\{m_1,~\ldots,~m_k\}$. It follows that $\xi_{M_0}^{J'}\in I(\lambda)$ for all $J'$ such that $\epsi_{j_1}\succeq\epsi_{j_k'}$. This means that the rank of $[\lambda]_p$, $p=(m_k,~j_1)$, is less than $k$ and so is not maximal.

Suppose now that $\Mo=\{m_1,~\ldots,~m_{k'}\}$ is finite and $\epsi_{m_1}\succ\ldots\succ\epsi_{m_{k'}}$. If $k'\geq k$ then the proof is similar to the previous case, so let $k'<k$. Arguing as above, we see that $\xi_{I^{\star}}^{J^{\star}}\in I(\lambda)$ for all $I^{\star}$, $J^{\star}$ of cardinality $k$ such that $\epsi_{i_k^{\star}}\succeq\epsi_{i_1}$, $\epsi_{j_1}\succeq\epsi_{j_k^{\star}}$, and $i_s^{\star}=m_s$ for $1\leq s\leq k'$. If there exists $j_0$ such that $\xi_{\Mo}^{J'}(\lambda)=0$ for all $J'$ of cardinality $k'$ with $\epsi_{j_0}\succeq\epsi_{j_{k'}'}$, then the rank of $[\lambda]_p$, $p=(m_{k'},~j_0)$, is less than~$k'$ and so is not maximal. Thus, we may assume that, given $j_0$, there exists $J'$ of cardinality $k'$ such that $\xi_{\Mo}^{J'}(\lambda)\neq0$ and $\epsi_{j_0}\succeq\epsi_{j_k'}$. But in this case the standard methods of matrix rank calculation show that the rank of a submatrix of $[\lambda]_p$, $p=(i_1,~j_1)$, with the rows $I^{\star}$ and the columns $J^{\star}$ such that $|I^{\star}|\geq k$, $|J^{\star}|=k+k'$, $\epsi_{i_k^{\star}}\succeq\epsi_{i_1}$, $\epsi_{j_1}\succeq\epsi_{j_k^{\star}}$, and $j_{k+s}^{\star}=j_s'$, $i_s^{\star}=m_s$ for $1\leq s\leq k'$, is less than $k$. Thus, $\rk[\lambda]_p$ is less than $k$ and so is not maximal.

\textsc{Case} $\Phi=C_{\infty}$. Let $\xi_I^I\in I(\lambda)$ for some $I$ of cardinality $k$. Pick a root $\alpha=\epsi_i-\epsi_j$. We claim that if $i\in I$, $j\notin I$ then there exists a nonzero scalar $c\in\Cp^{\times}$ such that $\{\xi_I^I,~e_{\alpha}\}=c\xi_{I[i\to j]}^I$, where $I[i\to j]$ is obtained from $(I\setminus\{i\})\cup\{j\}$ by reordering the indices. 
Pick a finite subset $M\subset\Zp_{>0}$ containing $I$ and $j$. Let $\unt_n$, $\ut$, etc. be as in the proof of Proposition~\ref{prop:ideal_p_k_Poisson_ifd}. Similarly, we identify $\unt_n^*$ with $\unt_n^T$ and $\ut^*$ with $\ut^T$; we also identify $\alpha$ and $\xi_I^I$ with their images in $\Phi_n$ and $\Sa(\unt_n)$ respectively. Note that $\unt_n^*$ can be now considered as a subspace of $\ut^*$. Assume for simplicity of notation that $\alpha=\epsi_i-\epsi_j\in\Phi_n$ (i.e., that the indices $i$ and $j$ did not change after passing from $\Phi_M$ to $\Phi_n$).

Let $\spt_{2n}(\Cp)\subset\slt_{2n}(\Cp)$ be as in Subsection~\ref{sst:nilradicals_fin_dim}, so that $\unt_n$ and $\ut$ are the nilradicals of the Borel subalgebras of $\spt_{2n}(\Cp)$ and $\slt_{2n}(\Cp)$ respectively consisting of all upper-triangular matrices from these Lie algebras. The natural embedding $\unt_n^T\hookrightarrow\ut^T$ together with the isomorphisms $(\unt_n^T)^*\cong(\unt_n^*)^*\cong\unt_n$ and $(\ut^T)^*\cong(\ut^*)^*\cong\ut$ define the projection $\ut\to\unt_n$ and, consequently, the surjective morphism of associative algebras $\zeta\colon\Sa(\ut)\to\Sa(\unt_n)$. According to (\ref{formula:n_*_n_T_dual_basis}), given $e_{p,q}\in\ut$, one has
\begin{equation*}
\zeta(e_{p,q})=\begin{cases}
e_{2\epsi_p}&\text{if }p=-q>0,\\
e_{\epsi_p-\epsi_q}/2&\text{if }0<p<q,\\
e_{\epsi_p+\epsi_q}/2&\text{if }0<p<-q.\\
\end{cases}
\end{equation*}
(Recall that we index the rows and the columns of matrices by the numbers $1,~\ldots,~n,~-n,~\ldots,~-1$.) It is well known that the trace form is $\unt_n$-invariant and this provides an $\unt_n$-module isomorphism of $\ut$ and the quotient of $\slt_{2n}(\Cp)$ by the corresponding Borel subalgebra. 
In a similar way we have that $\zeta$~is in fact a morphism of $\unt_n$-modules.

Next, given two sequences $I'$, $J'$ of cardinality $k$, denote by $\kappa_{I'}^{J'}$ the element of
$\Sa(\ut)$ defined by formula (\ref{formula:xi_I_J_A_infty}) with $i_s'$, $-j_s'$ instead of $i_s$,
$j_s$ respectively, $1\leq s\leq k$. It follows from the previous paragraph that $\zeta(\kappa_I^I)=\xi_I^I/2^k$ and $$\zeta(\kappa_{I[i\to j]}^I)=\zeta(\kappa_I^{I[i\to j]})=\xi_{I[i\to j]}^I/2^k.$$ Further, $\zeta$ is a morphism of $\unt_n$-modules, hence $\zeta(\{\kappa_I^I,e_{i,j}-e_{-j,-i}\})=
\{\zeta(\kappa_I^I),\zeta(e_{i,j})\}=\{\xi_I^I,e_{\alpha}\}/2^k$. It remains to note that, thanks to the row expansion formula for determinants, $\{\kappa_I^I,e_{i,j}\}=\pm\kappa_{I[i\to j]}^I$ and $\{\kappa_I^I,e_{-j,-i}\}=\mp\kappa_I^{I[i\to j]}$.

Now we turn again from $\unt_n$ to $\nt$. We have just checked that, up to a nonzero scalar, $\{\xi_I^I,~e_{\alpha}\}$ coincides with $\xi_{I[i\to j]}^I$. It is clear from the row expansion formula for the determinants that, up to a nonzero scalar, $\{\xi_{I[i\to j]}^I,e_{\alpha}\}$ is equal to $\xi_{I[i\to j]}^{I[i\to j]}$. It is easy to deduce from this fact that $\xi_{I'}^{J'}\in I(\lambda)$ for all $I'$, $J'$ of cardinality $k$ such that $\epsi_{i_k'}\succ\epsi_{i_1}$, $\epsi_{j_1'}\succ\epsi_{i_1}$ and either $I'=J'$ or $|I'\setminus J'|=1$. In \cite{FallatMartinezRivera1}, such minors of our matrix (after reordering the indices making the matrix symmetric) are called \emph{quasiprincipal}. In other words, we have proved that all $k\times k$ quasiprincipal minors of the matrix $[\lambda]_p$, $p=(i_1,-i_1)$, are zero. It follows from \cite[Theorem 2.6]{FallatMartinezRivera1} that all quasiprincipal minors of $[\lambda]_p$ of size at least $k$ are zero. By the Principal Minor Theorem, the rank of $[\lambda]_p$ is less than $k$, hence is finite and so is not maximal.

\textsc{Case} $\Phi=D_{\infty}$. (The proof for $B_{\infty}$ is completely similar, so we omit it.) Pick $\xi_I\in I(\lambda)$ for some~$I$ of cardinality $2k$. Such an element $\xi_I\in I(\lambda)$ exists according to the third paragraph of the proof. Recall that $\Bu=\varnothing$. First, assume that the linear order $\preceq$ does not have the maximal element. It follows immediately from the row expansion formula for Pfaffians that if $i\in I$ and $j\notin I$ then $$\{\xi_I,e_{\epsi_i-\epsi_j}\}=\pm\xi_{I[i\to j]}.$$ Hence, arguing as above, we conclude that $\xi_{I'}\in I(\lambda)$ for all $I'$ of cardinality~$2k$ with $\epsi_{i_{2k}'}\succeq\epsi_{i_1}$. Applying again the row expansion formula for Pfaffians and the Principal Minor Theorem, we see that $\rk[\lambda]_p<k<\infty$ is not maximal, where $p=(i_1,~-i_1)$.

Finally, assume that the linear order $\preceq$ has the maximal element $\epsi_m$. Since $\Bu=\varnothing$, the order $\preceq$ does not have the second maximal element. If there exists $j'\in\Zp_{>0}$, $j'\neq m$, such that $$\lambda(e_{\epsi_m+\epsi_{j''}})=0$$ for all $\epsi_{j''}\succeq\epsi_{j'}$ then $\rk[\lambda]_p=0<1$, $p=(m,-j')$, is not maximal. Hence, given $j'\in\Zp_{>0}$, $j'\neq m$, there exists $j''\in\Zp_{>0}$, $j''\neq m$, such that $\epsi_{j''}\succ\epsi_{j'}$ and $\lambda(e_{\epsi_m+\epsi_{j''}})\neq0$.

Arguing as above, we see that $\xi_{I'}\in I(\lambda)$ for all $I'$ of cardinality $2k$ such that $\epsi_{i_{2k}'}\succeq\epsi_{i_1}$ and $i_1'=m$. We will show that the rank of $[\lambda]_p$ is less than $2k$ (and so finite and not maximal), where $p=(i_1,~-i_1)$. It is enough to check that the rank of each its square submatrix with the row set $I'$ and the column set $-I'$ is less than $2k$, where $I'$ is a sequence of cardinality $2k'\geq2k$ with $\epsi_{i_{2k'}'}\succeq\epsi_{i_1}$ and $i_1'\neq m$. Assume to the contrary that its rank is at least $2k$. As we noticed, there exists $j'\in\Zp_{>0}$ such that $\epsi_{j'}\succ\epsi_{i_1'}$ and $$\lambda(e_{\epsi_m+\epsi_{j'}})=\xi_{m,{j'}}(\lambda)\neq0.$$ Clearly, the rank of the submatrix of $[\lambda]_p$ with the row set $I''=\{m,~j'\}\cup I'$ and the column set $-I''$ is at least $2k$. According to \cite[Theorem 6]{Thompson1}, there are sequences $$I_1=\{m,~j'\}\subsetneq I_2\subsetneq\ldots\subsetneq I_k$$ contained in $I''$ such that $|I_s|=2s$ for $1\leq s\leq k$ and $\xi_{I_s}(\lambda)\neq0$. But all these sequences contain $m$, so $\xi_{I_k}(\lambda)=0$. This contradicts our assumption and hence we have finished the proof.}
We would like to discuss two Corollaries of Theorem~\ref{theo:non_trivial_ifd}. The first one is as follows.
\coro{\label{Cndp} Let $\nt$ be the nil-Dynkin Lie algebra. Then there exists a countable collection $\{ \tilde{}\Xi_k\}_{k}$ of countable collections of polynomials from $\Sa(\nt)$ such that $I(\lambda)\ne0$ if and only if there exists $k$ for which $\lambda(\xi)=0$ for all $\xi\in\Xi_k$. 
}
Of course, Corollary~\ref{Cndp} is much weaker then Theorem~\ref{theo:non_trivial_ifd}. 
Nevertheless, it allows to state the following conjecture.
\okr{Conjecture}{\label{Conln} Let $\nt$ be a locally nilpotent Lie algebra with $Y(\nt)=\Cp$. Then there exists a countable collection $\{ \tilde{}\Xi_k\}_{k}$ of countable collections of polynomials from $\Sa(\nt)$ such that $I(\lambda)\ne0$ if and only if there exists $k$ for which $\lambda(\xi)=0$ for all $\xi\in\Xi_k$. 
}
The following corollary provide an interesting enhancement of Theorem~\ref{theo:non_trivial_ifd} under the assumption that $\Bu=\varnothing$.
\corop{Let $\nt$ be a nil-Dynkin Lie algebra and let $I\subset\Sa(\nt)$ be a nonzero prime Poisson ideal. If $\Bu=\varnothing$ then $I\supseteq I(p,k)$ for some upper-right pair $p$ and some $k\in\Zp_{>0}$ such that $I(p,k)\neq\{0\}$.}{
Denote by $A$ the quotient $\Sa(\nt)/I$. 
Assume to the contrary that $I\not\supset I(p, k)$ for all $p$, $k$ with nonzero $I(p,k)$. 
This implies that for each such a pair there exists $f_{p, k}\in I(p, k)\setminus I$. 
It is clear that the collection of such polynomials $f_{p, k}$ is at most countable and we form a sequence $f_1, f_2, \ldots$ out of them. 
Denote by $\bar f_1, \bar f_2,\ldots$ the images of $f_1, f_2, \ldots$ in $A$. Thanks to the definition of $f_1, f_2, \ldots$ we have $\bar f_k\ne0$ for all $k\ge1$. 
This allows us to consider the localization of $A[\bar f_1^{-1}, \bar f_2^{-1},\ldots]$ of $A$ at $\bar f_1, \bar f_2, \ldots$. 
Let $M$ be a maximal ideal of $A[\bar f_1^{-1}, \bar f_2^{-1},\ldots]$ and $\wt M$ be the preimage of $M$ in $\Sa(\nt)$. 
Thanks to Corollary~\ref{Lschurc} we have that $\wt M$ has codimension 1 in $\Sa(\nt)$ and thus $\wt M$ is the maximal ideal attached to a properly chosen $\lambda\in\nt^*$. 
The definition of $\lambda$ and $f_1, f_2, \ldots$ implies that $I(\lambda)\supset I$ and thus $I(\lambda)\ne0$ and that $I(\lambda)\not\supset I(p, k)$ for all $(p, k)$. This contradicts Theorem~\ref{theo:non_trivial_ifd}.
}

\exam{Below we draw schematically matrices $[\lambda]_p$ for different choices of $\lambda$ and $p$.

i) First,\label{exam:nontrivial_ifd} let $\Phi=A_{\infty}$. Taking into account our speculations from the proof of Proposition~\ref{prop:ideal_p_k_Poisson_ifd}, it would be more natural to draw $\nt^*$ as the set of infinite lower-triangular matrices. Assume, for example, that $\epsi_1\succ\epsi_2\succ\ldots$, and $p=(2,~4)$. Below we marked in gray boxes $(j',~i')$ for which $\epsi_{i'}\succeq\epsi_2$ and $\epsi_{j'}\preceq\epsi_4$ (i.e., $i'\leq2$, $j'\geq4$), so that $\lambda(e_{\epsi_{i'}-\epsi_{j'}})$ is the $(i'~,j')$th entry of $[\lambda]_p$. (For $A_{\infty}$, the picture for other choices of $p$ looks essentially like this.)}
\begin{center}
\includegraphics{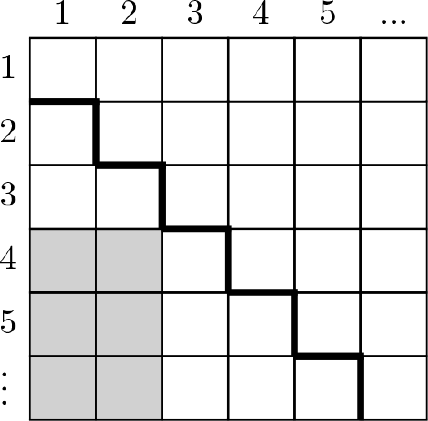}
\end{center}

ii) Another example: let $\Phi=C_{\infty}$ or $D_{\infty}$, $\epsi_{i+1}\succ\epsi_i$ for all $i\geq1$, and $p=(2,~-2)$. Given $i,~j\in\Zp_{>0}$, $\epsi_i\succ\epsi_j$, it is convenient to identify the box $(i,~j)$ (respectively, $(i,-j)$) with the root $\epsi_i-\epsi_j$ (respectively, $\epsi_i+\epsi_j$). Besides, the box $(i,~-i)$ corresponds to the root $2\epsi_i$ or is filled by zero for $\Phi=C_{\infty}$ or $D_{\infty}$ respectively.
Below we marked in gray boxes $(-j',~i')$ for which $\epsi_{i'}\succeq\epsi_2$ and $\epsi_{j'}\succeq\epsi_2$, i.e., boxes corresponding to roots $\alpha$ such that $\lambda(e_{\alpha})$ is involved in $[\lambda]_p$.
\begin{center}
\includegraphics{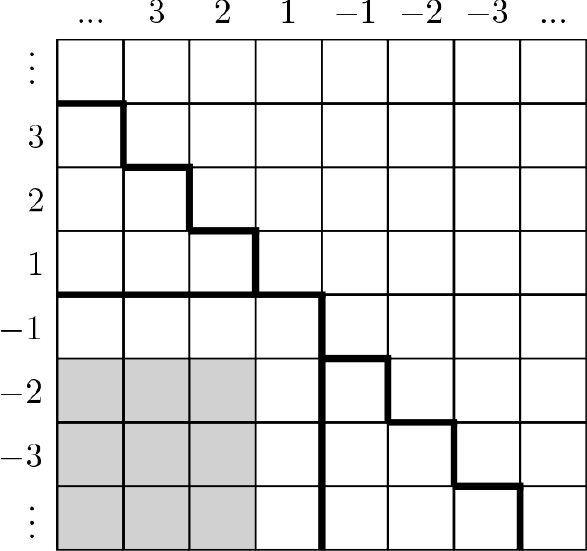}
\end{center}

iii) Finally, let $\Phi=D_{\infty}$, $\epsi_1\succ\epsi_i$ and $\epsi_{i+1}\succ\epsi_i$ for all $i\geq2$, and $p=(1,-3)$. Below we marked in gray the boxes $(-j',~1)$ and $(-1,~j')$ for which $\epsi_{j'}\succeq\epsi_2$, so that $\lambda(e_{\epsi_1+\epsi_{j'}})$ is the $(1,~{-j'})$th entry of~$[\lambda]_p$.
\begin{center}
\includegraphics{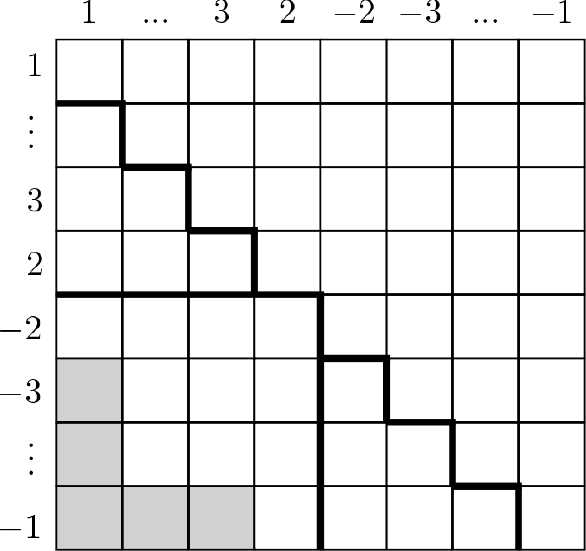}
\end{center}

We hope that these pictures explain why we call $p$ ``an upper-right pair'': in fact, the entries of $[\lambda]_p$ equal $\lambda(e_\alpha)$, where $\alpha$ is the root corresponding to a gray box situated nonstrictly to the South-West from the box corresponding to $p$. (On the last picture we also marked in gray boxes in the bottom row to stress that the picture is antisymmetric with respect to the antidiagonal.)



\medskip\textsc{Mikhail V. Ignatyev: Samara National Research University, Ak. Pavlova 1,\break\indent Samara 443011, Russia}

\emph{E-mail address}: \texttt{mihail.ignatev@gmail.com}

\medskip\textsc{Alexey Petukhov: Institute for Information Transmission Problems, Bolshoy\break\indent Karetniy 19-1, Moscow 127994, Russia}

\emph{E-mail address}: \texttt{alex--2@yandex.ru}

\end{document}